\documentclass[11pt,  a4paper,reqno]{amsart}
\usepackage{amssymb,amsmath,latexsym, amsthm}
\usepackage{color}
\usepackage{comment}
\usepackage{tikz}
\usepackage{hyperref}
\allowdisplaybreaks

\newcommand{\normal}{\color{black}}
\newcommand{\bk}{\color{black}}
\usepackage{mathtools}
\usepackage{amsthm}

\def\rF {{\overline \sF}}
\def\rY {\overline Y}
\def\rP {\overline P}
\def\oR {{\overline{ \mathbb R}}}
\def\1{{\bf 1}}

\def\nn{\nonumber}
\def\bee{\begin{equation}}
\def\eee{\end{equation}}
\def\sA {{\mathcal A}} \def\sB {{\mathcal B}} \def\sC {{\mathcal C}}
 \def\sE {{\mathcal E}} \def\sF {{\mathcal F}}

 \def\sN {{\mathcal N}}

\def\R {{\mathbb R}}

\newtheorem{thm}{Theorem}[section]
\newtheorem{lemma}[thm]{Lemma}
\newtheorem{defn}[thm]{Definition}
\newtheorem{prop}[thm]{Proposition}
\newtheorem{corollary}[thm]{Corollary}
\newtheorem{remark}[thm]{Remark}
\newtheorem{example}[thm]{Example}
\numberwithin{equation}{section}

\def\qed{{\hfill $\Box$ 
\smallskip
}}

\def\NN{{\mathcal N}}

\def\FF{{\mathcal F}}
\def\EE{{\mathcal E}}

\def\R{{\mathbb R}}

\def\E{{\mathbb E}}

\def\P{{\mathbb P}}
\def\N{{\mathbb N}}
\def\eps{\varepsilon}
\def\wh{\widehat}
\def\wt{\widetilde}
\def\pf{\noindent{\bf Proof.} }

\def\e{{\mathbf e}}
\usepackage{hyperref}
\begin{document}
\title
[Dirichlet forms with jump kernels blowing up at the boundary]
{Potential theory of Dirichlet forms with jump kernels blowing up at the boundary}
\author{ Panki Kim \quad Renming Song \quad and \quad Zoran Vondra\v{c}ek}

\address[Kim]{Department of Mathematical Sciences and Research Institute of Mathematics,
	Seoul National University,
	Seoul 08826, Republic of Korea}\thanks{This research is  supported by the National Research Foundation of Korea(NRF) grant funded by the Korea government(MSIP) (No. RS-2023-00270314).
}
\curraddr{}
\email{pkim@snu.ac.kr}

\address[Song]{
Department of Mathematics, University of Illinois, Urbana, IL 61801,
USA}
\curraddr{}\thanks{Research supported in part by a grant from
the Simons Foundation (\#960480, Renming Song)}
\email{rsong@illinois.edu
}

\address[Vondra\v{c}ek]
{
Department of Mathematics, Faculty of Science, University of Zagreb and 
``Dr. Franjo Tu\dj man'' Defense and Security University, Zagreb, Croatia
}
\curraddr{}\thanks{ Research supported in part by the Croatian Science Foundation under the project 4197. 
(Zoran Vondra\v{c}ek)}
\email{vondra@math.hr}

 \date{}

\begin{abstract}
In this paper we study the potential theory of Dirichlet forms 
on the half-space $\R^d_+$ defined by the jump kernel $J(x,y)=|x-y|^{-d-\alpha}\sB(x,y)$ 
and the killing potential $\kappa x_d^{-\alpha}$,
where $\alpha\in (0, 2)$ and $\sB(x,y)$ can blow up to infinity at the boundary.
The jump kernel and the killing potential 
depend
on several parameters. 
 For all admissible values of 
the parameters involved and all $d \ge 1$, 
we prove that the boundary Harnack principle holds, and  establish sharp two-sided estimates on the Green functions of these processes.
\end{abstract}
 \date{}

\maketitle

\noindent {\bf MSC 2020 Mathematical Sciences 
Classification
 System:} Primary 60J45; Secondary 60J46, 60J76.


\section{Introduction}\label{s:intro}

In this paper, we study the potential theory of purely discontinuous symmetric Markov processes in the upper half-space
$\R^d_+:=\{x=(\wt{x}, x_d):\, x_d>0\}$, $d\ge 1$, with  jump kernel of the form $J(x,y)=|x-y|^{-d-\alpha}\sB(x,y)$, $\alpha\in (0,2)$, where $\sB(x,y)$ is degenerate at
the boundary of $\R^d_+$.  In our recent papers  \cite{KSV, KSV20, KSV21}, we have studied the case when  $\sB(x,y)$ decays to zero at the boundary.
In this paper, we study the case when $\sB(x,y)$ blows up at the boundary and    establish the boundary Harnack principle and sharp two-sided estimates on the Green functions.

One of our main motivation to study this problem 
comes from the following natural example of a process with 
jump kernel
blowing up at the boundary. Let $X=(X_t, \P_x)$ be an isotropic $\alpha$-stable process in $\R^d$. Define
$A_t:=\int_0^t \1_{(X_s\in \R^d_+)}ds$
and let $\tau_t:=\inf\{s>0:\, A_s>t\}$ be its right-continuous inverse. The process $Y=(Y_t)_{t\ge 0}$, defined by $Y_t=X_{\tau_t}$,  is a Hunt process 
on $\R^d_+$, called the \emph{trace} process of $X$ on $\R^d_+$ (the name \emph{path-censored} process is also used in some literature, see \cite{KPW14}). 
The part of the process $Y$ until its first hitting time of the boundary $\partial \R^d_+=\{(\wt{x}, 0):\, \wt{x}\in \R^{d-1}\}$ can be described in the following way: 
Let $\tau=\tau_{\R^d_+}=\inf\{t>0:\, X_t\notin \R^d_+\}$ be the exit time of $X$ from $\R^d_+$, $x=X_{\tau-}\in \R^d_+$ the position from which $X$ jumps 
out of $\R^d_+$, and $z=X_{\tau}$ be the position where $X$ lands at the exit from $\R^d_+$. 
Then $z\in \R^d_-$ a.s., where $\R^d_-:=\{x=(\wt{x}, x_d):\, x_d<0\}$. The distribution
of the returning position of $X$ to $\R^d_+$ is given by the Poisson kernel 
of the process $X$ in $\R^d_-$   (i.e., the density of the distribution of $X_{\tau_{\R^d_-}}$ on $\R^d_+$):
\begin{equation}\label{e:poisson-100}
P_{\R^d_-}(z,y)=\int_{\R^d_-}G^X_{\R^d_-}(z,w)j(w, y)\, dw, \quad y\in \R^d_+.
\end{equation}
Here 
$G^{X}_{\R^d_-}(z,w)$
is
the Green function of the process $X$ killed upon 
exiting 
$\R^d_-$, $j(w, y)=\sA(d, \alpha)|w-y|^{-d-\alpha}$ is the jump kernel of $X$ and $\sA(d, \alpha)=2^{\alpha}\pi^{-d/2}\Gamma ((d+\alpha)/2)/|\Gamma(-\alpha/2)|$.

This implies that when $X$ jumps out of $\R^d_+$ from the point $x$, we continue the process by
resurrecting it at $y\in \R^d_+$ according to the kernel
\begin{equation}\label{e:rsnew1}
q(x, y):=\int_{\R^d_-}j(x, z)P_{\R^d_-}(z,y)\, dz, \quad x\in \R^d_+.
\end{equation}
We will call $q(x,y)$ a \emph{resurrection kernel}. Since the Green function $
 G^{X}_{\R^d_-}
(\cdot, \cdot)$ is symmetric, it follows that $q(x,y)=q(y,x)$ for all $x,y\in \R^d_+$. The kernel $q(x, y)$ introduces additional jumps from $x$ to $y$. 
By using Meyer's construction (see \cite{M}), 
one can construct a resurrected process 
on
$\R^d_+$ with jump kernel 
$J(x,y)=j(x, y)+q(x,y)$. 
The resurrected process is equal to the part of the trace process $Y$ until it first hits $\partial \R^d_+$. It follows from \cite[Theorem 6.1]{BGPR} (where $q(x,y)$ is called the \emph{interaction kernel}) that in case $d\ge 3$, 
$$
J(x,y)\asymp q(x,y)\asymp |x-y|^{–d-\alpha}\left(\frac{|x-y|^2}{x_d y_d}\right)^{\alpha/2}, \quad x_d\wedge y_d\le |x-y|.
$$
This asymptotic relation shows that the jump kernel $J(x,y)$ blows up with rate $x_d^{-\alpha/2}$ when $x$ approaches the boundary $\partial \R^d_+$. Here and throughout the paper, the notation $f\asymp g$ for non-negative functions $f$ and $g$ means
that there exists a constant $c\ge 1$ such that $c^{-1}g\le f \le cg$. We also use $a\wedge b:=\min\{a,b\}$ and $a\vee b:=\max\{a,b\}$.

Another motivation for this paper is the process introduced in \cite{DROV, Von21} to study non-local Neumann problems.  
See also \cite{FM} and the references therein.
For the process in \cite{DROV, Von21}, the resurrection kernel $q(x, y)$ 
is given by \eqref{e:rsnew1} with
 the Poisson kernel 
 $P_{\R^d_-}(z,y)$ replaced by
$j(z, y)/\int_{\R^d_+}j(z, w)dw$. The jump kernel of this process also blows up at the boundary, see Remark \ref{r:estimate}(b).

In Section \ref{s:res-ker} we substantially generalize  
these two examples by replacing the Poisson kernel $P_{\R^d_-}(z,y)$ and the kernel 
$j(z, y)/\int_{\R^d_+}j(z, w)dw$
by a very general \emph{return kernel} $p(z,y)$. 
The kernel $p(z,y)$, $z\in \R^d_-$, $y\in \R^d_+$, is chosen so that the corresponding resurrection kernel 
$$
q(x,y)
=
\int_{\R^d_-}j(x, z)p(z,y)\, 
 dz,  \quad x,y\in \R^d_+,
$$
is symmetric. This flexibility in choosing the return kernel allows us to obtain resurrection kernels with various blow-up rates at the boundary. The main result in this direction is Theorem \ref{t:estimate}.  

Note that the jump kernel 
$J(x,y)=j(x, y)+q(x,y)$ 
of the resurrected process may be written in the form 
$J(x,y)=j(x, y)\sB(x,y)$ with 
$\sB(x,y):=1+q(x,y)/j(x, y)$. 
Since the jump kernel $j(x, y)$
is bounded away from the diagonal, the blow up at the boundary comes from the term $\sB(x,y)$.
The estimates in Theorem \ref{t:estimate} contain also the asymptotics of the term $\sB(x,y)$ 
and imply that the resurrected process satisfies \textbf{(A1)}--\textbf{(A4)} below.
The proof of Theorem \ref{t:estimate} is quite long and technical and is therefore postponed to Section \ref{s:proof-res-ker}. Let us mention that Sections \ref{s:res-ker} and \ref{s:proof-res-ker} are logically independent from the rest of the paper, and  
also serve as the motivation for the general set-up that we now introduce.

\smallskip
Let $d \ge 1$, $\alpha\in (0,2)$
and assume that 
$0 \le \beta_1\le  \beta_2< 1\wedge\alpha $. 
Let $\Phi$ be a positive function on $[2, \infty)$ satisfying  
the following weak scaling condition:  
There exist constants  $C_1, C_2>0$ such that 
\begin{equation}\label{e:Phi-infty}
 C_1  (R/r)^{\beta_1}\leq\frac{\Phi(R)}{\Phi(r)}\leq C_2(R/r)^{\beta_2},
\quad 2\le r<R <\infty.
\end{equation}
For notational convenience, we extend the domain of $\Phi$ to $[0, \infty)$ by letting 
$\Phi(t)  \equiv \Phi(2)>0$ on $[0, 2)$. Then for any $\delta>0$,
there exist constants 
$\wt C_1, \wt C_2>0$ 
depending on $\delta$ such that 
$$
 \wt C_1(R/r)^{\beta_1}\leq\frac{\Phi(R)}{\Phi(r)}\leq \wt C_2(R/r)^{\beta_2},
 \quad \delta\le r<R <\infty.
$$
Let 
$\wt{\beta}_2$ be the upper Matuszewska index of $\Phi$ (see \cite[pp.~68-71]{BGT}):
$$
\wt{\beta}_2:=\inf\{\beta>0:\, \exists a\in (0,\infty) \text{ 
such that }\, \Phi(R)/\Phi(r)\le a(R/r)^{\beta} \text{ for }2\le r<R<\infty\}.
$$
Note that the inequality $\Phi(R)/\Phi(r)\le 
a(R/r)^{\wt{\beta}_2}$ may, but need not hold for any $a\in (0,\infty)$. 

Define 
$$
j(x, y)=\sA(d, \alpha)|x-y|^{-\alpha-d} 
\quad \text{and} \quad  J(x,y)=j(x, y)\sB(x,y).
$$
We will assume that $\sB(x,y)$ satisfies the following conditions:

\smallskip
\noindent
\textbf{(A1)} $\sB(x,y)=\sB(y,x)$ for all $x,y\in  \R^d_+$.

\smallskip
\noindent
\textbf{(A2)}   If $\alpha \ge 1$, there exists $\theta>\alpha-1$ such that for every $a>0$ there exists
$C=C(a)>0$ such that 
$$
|\sB(x, x)-\sB(x,y)|\le 
C\left(\frac{|x-y|}{x_d\wedge y_d}\right)^{\theta} \quad \text{ for all }x,y\in \R^d_+ \text{ with }x_d \wedge y_d  \ge a|x-y|.
$$ 

\smallskip
\noindent
\textbf{(A3)}
There exists $C \ge 1$ such that
\begin{equation}\label{e:B7}
C^{-1} \Phi\left(\frac{|x-y|^2}{x_dy_d}\right)\le \sB(x,y)
\le C\Phi\left(\frac{|x-y|^2}{x_dy_d}\right) \quad \text{ for all }x,y\in \R^d_+ .
\end{equation}

\smallskip
\noindent
\textbf{(A4)} 
For all $x,y\in \R^d_+$ and  $a>0$, $\sB(ax,ay)=\sB(x,y)$. 
In case $d\ge 2$, for  all 
$x,y\in \R^d_+$ and $\wt{z}\in \R^{d-1}$, $\sB(x+(\wt{z},0), y+(\wt{z},0))=\sB(x, y)$.

\smallskip
Note that \textbf{(A3)} implies that $\sB(x,y)$ is bounded from below by a positive constant, and \textbf{(A4)} implies that $x\mapsto \sB(x,x)$ is constant.
\smallskip

For $\kappa\in [0,\infty)$ we define the function $\kappa(x):=\kappa x_d^{-\alpha}$ on $\R^d_+$ and set
$$
\EE^\kappa(u,v):=\frac12 \int_{\R^d_+}\int_{\R^d_+}(u(x)-u(y))(v(x)-v(y))J(x,y) dy dx+  \int_{\R^d_+} u(x)v(x)\kappa(x)dx,
$$
where $u,v:\R^d_+\to \R$. 
Let $\FF^0$ be the closure of $C_c^{\infty}(\R^d_+)$ in $L^2(\R^d_+, dx)$ under $\EE^0_1 :=\EE^0+(\cdot, \cdot)_{L^2(\R^d_+, dx)}$. 
Then, due to $\beta_2< 1\wedge\alpha$, 
 $(\EE^0,\FF^0)$ is a regular Dirichlet form on $L^2(\R^d_+, dx)$
 (see Section \ref{s:prelim} below). 
 Let 
$$
\sF^\kappa:= 	\wt \sF^0 \cap L^2(\R^d_+, \kappa(x)dx),\nn
$$
where $\wt \sF^0$ is the family of all $\sE^0_1$-quasi-continuous functions in $\sF^0$. Then $(\sE^\kappa,\sF^\kappa)$ is also a regular Dirichlet form on $L^2(\R^d_+,dx)$.
As we will explain in Section \ref{s:prelim}, under assumptions \textbf{(A1)}-\textbf{(A4)}, there exists a symmetric, scale invariant and horizontally translation invariant Hunt process $Y^{\kappa}=((Y_t^{\kappa})_{t\ge 0}, (\P_x)_{x\in \R^d_+})$ associated with $(\EE^{\kappa}, \FF^{\kappa})$. 
In case $\kappa>0$, the process $Y^{\kappa}$ is transient.
To show these facts we will use results proved in \cite{KSV22b}.

We now associate the constant $\kappa$ from the killing function $\kappa(x)=\kappa x_d^{-\alpha}$ with a 
positive  parameter $p=p_{\kappa}$ which will play a major role in the paper.
Let $\e_d:=(\wt{0}, 1)$. For $q \in (-1, \alpha -\wt{\beta}_2)$, set
$$
	C(\alpha,q,\sB) =
	\begin{cases}
	\displaystyle 	 \int_{\R^{d-1}}  \int_0^1 \frac{(s^q-1)(1-s^{\alpha-q-1})}{(1-s)^{1+\alpha}} \frac{\sB\big( (1-s)\wt u,1) ,s \e_d\big)}{(|\wt u|^2+1)^{(d+\alpha)/2}}	 ds d\wt u
	 & \text{if } d\ge 2 \\[2ex]
	\displaystyle   \int_0^1 \frac{(s^q-1)(1-s^{\alpha-q-1})}{(1-s)^{1+\alpha}} \sB\big( 1 ,s \big)ds  & \text{if } d=1.
	\end{cases}
$$
Then $C(\alpha, 0, \sB)=C(\alpha, \alpha-1, \sB)=0$ and
the function $q\mapsto C(\alpha, q, \sB)$ is strictly increasing 
and continuous
on $[(\alpha-1)_+, \alpha-\wt{\beta}_2)$. Consequently, for every $0\le \kappa<\lim_{q \uparrow \alpha-\wt{\beta}_2} C(\alpha,q,\sB)\le \infty$,  there exists a unique $p_\kappa \in [(\alpha-1)_+, \alpha-\wt{\beta}_2)$ such that
$\kappa=C(\alpha, p_\kappa, \sB)$.
 When $\Phi(r)=r^{\beta}$ with $\beta\in (0,1\wedge \alpha)$, it holds that 
$\lim_{q \uparrow \alpha-\beta} C(\alpha,q,\sB)= \infty$ (see Lemma  \ref{l:LB-on-g0}),
 so 
$\kappa\mapsto p_{\kappa}$ is an increasing bijection from $[0,\infty)$ onto $[(\alpha-1)_+, \alpha-\beta)$. 
In the remainder of this introduction
we will fix $\kappa \in [0,\lim_{q \uparrow \alpha-\wt{\beta}_2} C(\alpha,q,\sB))$,
and assume $\alpha>1$ if  $\kappa=0$ 
so that $p_0=\alpha-1>0$.
We will show in Section \ref{s:scaling} that $Y^0$ is transient when $\alpha\in (1,2)$.
For notational simplicity, in the remainder of this introduction, we omit the superscript $\kappa$ from the notation: 
For example, we write $Y$ instead of $Y^{\kappa}$, and $p$ instead of $p_\kappa$ in \eqref{e:killing-potential}.

The role of the parameter $p$ and its connection to $C(\alpha, p, \sB)$ 
can be seen from the following observation. Let 
$$
L^{\sB}f(x)=\textrm{p.v.}\int_{\R^d_+}(f(y)-f(x))J(x,y)\, dy- C(\alpha, p, \sB)x_d^{-\alpha} f(x), \quad x\in \R^d_+,
$$
whenever the principal value integral makes sense. If $g_p(x)=x_d^p$, then $L^{\sB}g_p\equiv 0$, see Lemma \ref{l:LB-on-g}. Hence the operator $L^{\sB}$ annihilates the $p$-th power of the distance to the boundary.

The first main result of the paper is the scale invariant boundary
Harnack principle with  
exact decay rate: If a non-negative harmonic function
vanishes continuously at a part of the boundary $\partial \R^d_+$, then the decay rate is equal to the $p$-th power of the distance to the boundary.

For an open subset $D$ of  $\R_+^d$, let
$\tau_D:=\inf\{t>0:\, Y_t\notin 
D\}$ be the first exit time of the process $Y$ from $D$. 
\begin{defn}\label{D:1.1}   
A non-negative Borel function defined on
$\R_+^d$ is said to be  {harmonic}
in an open set $V\subset \R_+^d$ with respect 
to $Y$ 
if for every bounded open set $
D\subset\overline{
D}\subset V$,
$$
f(x)= \E_x \left[ f(Y_{\tau_{D}}): \tau_{D}<\infty \, \right] \qquad \hbox{for all } x\in D.
$$
A non-negative Borel function $f$ defined on $\R_+^d$ is said to be \emph{regular harmonic} in an open set $V\subset \R_+^d$ if 
$$
f(x)= \E_x \left[ f(Y_{\tau_{V}}): \tau_{V}<\infty \, \right] 
\qquad \hbox{for all } x\in V.
$$
\end{defn}

When $d\ge 2$, for $a,b>0$ and $\wt{w} \in \R^{d-1}$, we define 
\begin{equation}\label{e:box}
D_{\wt{w}}(a,b):=\{x=(\wt{x}, x_d)\in \R^d:\, |\wt{x}-\wt{w}|<a, 0<x_d<b\}.
\end{equation}
By abusing notation,  in case $d=1$, we will use
$D_{\wt{w}}(a,b)$ to stand for the open interval $(0,b)=\{y\in \R_+:\, 0<y<b\}$.

\begin{thm}\label{t:BHP}
Suppose $p\in (0, \alpha-\wt{\beta}_2) \cap [(\alpha-1)_+, \alpha-\wt{\beta}_2)$.
 Assume that  $\sB$ satisfies  \textbf{(A1)}-\textbf{(A4)}. Then there exists $C \ge 1$ such that for all $r>0$, $\wt{w} \in \R^{d-1}$,  and any non-negative function $f$ in $\R^d_+$ which is harmonic in $D_{\wt{w}}(2r, 2r)$ 
with respect to $Y$ and 
vanishes continuously on $B(({\wt{w}},0), 2r)\cap \partial \R^d_+$,  we have 
\begin{equation}\label{e:TAMSe1.8new}
\frac{f(x)}{x_d^{ p}}\le C\frac{f(y)}{y_d^{p}}, \quad x, y\in D_{\wt{w}}(r/2, r/2) .
\end{equation}
\end{thm}

\smallskip
The second main result is on sharp two-sided estimates for the Green function of the process $Y$. We recall in Section \ref{s:prelim} the definition of the Green function $G(x,y)$, $x,y\in \R^d_+$, and comment on its existence.

\begin{thm}\label{t:gfe}
Suppose that
$p\in (0, \alpha-\wt{\beta}_2) \cap [(\alpha-1)_+, \alpha-\wt{\beta}_2)$ and that $\sB$ satisfies  \textbf{(A1)}-\textbf{(A4)}.
Then the process $Y$ admits a Green function $G:\R^d_+\times \R^d_+\to [0,\infty]$
 such that $G(x, \cdot)$ is 
 continuous in $\R^d_+\setminus \{x\}$ and regular harmonic with respect to $Y$ in $\R^d_+\setminus B(x, \epsilon)$ for any $\epsilon>0$. 
Moreover, $G(x,y)$ has the following estimates:  for all $x,y \in \R^d_+$, 
\begin{align}\label{e:FGE}
G(x,y)\asymp 
\begin{cases}
\displaystyle 
\left(\frac{x_d}{|x-y|}  \wedge 1 \right)^p\left(\frac{y_d}{|x-y|}  \wedge 1 \right)^p \frac{1}{|x-y|^{d-\alpha}},
&\alpha<d;\\[8pt]
\displaystyle \left(\frac{x\wedge y}{|x-y|} \wedge 1\right)^p\log\left(e+ \frac{ x \vee  y}{|x-y|}\right),
&\alpha=1=d;\\[8pt]
\displaystyle \left(\frac{x\wedge y}{|x-y|} \wedge 1\right)^p 
\left(x \vee y \vee |x-y|\right)^{\alpha-1},
&\alpha>1=d.
\end{cases}
\end{align} 
\end{thm}

Let us emphasize here that in case $\kappa=0$ and $\alpha>1$, by using several cutting-edge techniques developed here as well as in our previous papers \cite{KSV, KSV20, KSV21}, we succeeded to establish that, regardless of the blow-up rate of the function $\sB$, the decay rate of harmonic functions as well as the Green function is given by $p=\alpha-1$. We have shown in \cite{KSV21} that the same phenomenon also occurs in case when $\sB$ decays to zero at the boundary. In view of the fact that this is the same decay rate as for the censored $\alpha$-stable process (or, equivalently, the regional fractional Laplacian), this can be regarded as a stability result even for degenerate non-local operators.

Our strategy for proving the two main results above consists of several steps.

The first step is to show certain interior potential-theoretic results for the process $Y$. 
This is done in \cite{KSV22b} in a more general setting than 
that of the current paper. One of the key difficulties is the fact that 
$Y$ need not have the Feller property. 
Despite this obstacle we established a Dynkin-type formula  
on relatively compact open subsets $D$ of $\R^d_+$
for functions  in $C^2(\overline{D})$ defined on $\R^d_+$, 
see Theorem \ref{t:dynkin}.
Another important result coming from \cite{KSV22b} 
is the Harnack inequality, see Theorem \ref{t:uhp}. These and some other results are described in the preliminary Section \ref{s:prelim}.

The second step consists of studying the action of the operator $L^{\sB}$ on the powers of the distance to the boundary. This allows an extension of the Dynkin-type formula to \emph{not} relatively compact open sets $D(r,r)$
for functions $x_d^p\1_{D(R,R)}$ for $2r<R$, see Proposition \ref{p:dynkin-hp}. This extension together with Theorem \ref{t:dynkin}  plays a major role throughout this paper.

The third step is to establish certain exit probability estimates, 
see Lemma \ref{l:exit-probability-estimate} and Proposition \ref{p:allcomparable}.  
The key ingredient in proving these lemmas is to find suitable test functions (barriers) and to estimate the action of the operator $L^{\sB}$ on them. This is done in Lemma \ref{l:key-lemma}. 
The proof of this lemma is quite involved and relies on 
some rather delicate estimates of certain integrals due to the general nature of $\Phi$, see Lemma \ref{l:1n}.

The fourth step is the Carleson estimate, Theorem \ref{t:carleson}, for non-negative harmonic functions vanishing on a part of the boundary. The proof, although  standard, requires several modifications due to the blow up of the jump kernel at the boundary.

The next step consists of showing interior estimates for the Green function $G(x,y)$, see Propositions \ref{p:green-int-est} and \ref{p:green-int-est-2}. By interior we mean that the distance between $x$ and $y$ is small comparable to the distance of these points to the boundary. Here we distinguish two cases: $d>\alpha$ and $d=1\le \alpha$.  The proof of the upper bound of the former case uses the Hardy inequality, while the proof of the lower bound employs 
a capacity argument. In the latter case, we use the capacity estimates of the one-dimensional killed isotropic stable process and a version of the capacity argument for the process $Y$.

Next, we obtain
the preliminary upper bound of the Green function with correct boundary decay. We first show, see Theorem \ref{t:green-function-decay}, that the Green function decays at the boundary. This allows us to use the Carleson estimate and extend the upper interior estimate of $G(x,y)$ to all points $x,y\in \R^d_+$, cf.~Proposition \ref{p:gfcnub}. In 
Lemma \ref{l:prelub} we insert in the upper estimate the boundary part $\big(\frac{x_d\wedge y_d}{|x-y|}\wedge 1\big)^p$. The proof depends on delicate estimates of the jump kernel, and again, on the powerful Lemma \ref{l:1n}. 
As an application, in Proposition \ref{p:bound-for-integral-new} we give some upper estimates 
on the Green potentials of powers of the distance to the boundary.
These upper estimates, together with exit probability estimates, the Harnack inequality and the Carleson estimate, lead to a rather straightforward proof of Theorem \ref{t:BHP}.

Finally, we use the interior Green function estimates, the boundary Harnack principle and scaling to obtain the sharp two-sided Green function estimates.

We end this introduction with a few comments on the assumptions \textbf{(A1)}-\textbf{(A4)} and their relation to the assumptions in 
\cite{KSV, KSV20, KSV21}, where the jump kernel decays at the boundary. Assumption \textbf{(A1)} ensures the symmetry of the jump kernel and hence the process $Y$. Assumption \textbf{(A2)} is used in the analysis of the generator $L^{\sB}$, and allows to establish a Dynkin-type formula. Assumption \textbf{(A4)} is natural in the context of the half-space $\R^d_+$ and, in particular, ensures the scaling property of the process $Y$. These three assumptions were also postulated in \cite{KSV, KSV20, KSV21}. The main difference with those papers is in assumption \textbf{(A3)} which provides the blow-up of jump kernel at the boundary
and is motivated by Section \ref{s:res-ker}. 
In case when $\Phi(t)=t^{\beta}$, $t\ge 2$, for $0\le \beta <\alpha\wedge 1$, \textbf{(A3)} 
is equivalent to the condition
\begin{equation}\label{e:negative-beta}
 \sB(x,y) \asymp \left(\frac{x_d\wedge y_d}{|x-y|}\wedge 1\right)^{-\beta}\left(\frac{x_d\vee y_d}{|x-y|}\wedge 1\right)^{-\beta}.
\end{equation}
In \cite{KSV, KSV20, KSV21}, the assumptions on $\sB(x,y)$ included the case when
\begin{equation}\label{e:positive-beta}
 \sB(x,y) \asymp \left(\frac{x_d\wedge y_d}{|x-y|}\wedge 1\right)^{\beta}\left(\frac{x_d\vee y_d}{|x-y|}\wedge 1\right)^{\beta},
\end{equation}
with $\beta\ge 0$. In case $\beta>0$, this implies the decay of jump kernel at the boundary. Thus we can regard \eqref{e:negative-beta} as an extension of \eqref{e:positive-beta} from $\beta \in [0,\infty)$ to $\beta \in (-(\alpha\wedge 1), 0]$. Of course, \textbf{(A3)} is much more general than \eqref{e:negative-beta}.

It is instructive to look at the effect of blow-up and the decay of $\sB$ 
determined by $\beta\in (-(\alpha\wedge 1), \infty)$ on the range of possible values of the parameter $p$. By using  \cite{KSV}, we see that $p\in (0,\alpha+\beta)\cap [(\alpha-1)_+, \alpha+\beta)$.  By increasing 
the parameter $\beta$ from 0 to $\infty$ 
(and thus making the decay of $\sB$ sharper), the upper boundary of the range of $p$ also increases from $\alpha$ to $\infty$. On the other hand, by decreasing 
the parameter $\beta$ from 0 to $-(1\wedge \alpha)$
(and thus making the blow-up higher), the upper boundary of the range of $p$ decreases from $\alpha$ to $(\alpha-1)_+$. Therefore, the larger the blow-up at the boundary, the smaller the effect of the killing function.

\textbf{Notation:} Throughout this paper, 
capital $C$, with or without subscript, is used only for assumptions or the statements of results, 
while lower case  $c$  and $c_i$, $i=1, 2, \dots$, are used in the proofs. The value of $c$ may change from one appearance to another,
but the value of $c_i$ stays fixed in the same proof.
The notation $C=C(a,b,\ldots)$ indicates that the constant $C$
depends on $a, b, \ldots$.
We will use ``$:=$" to denote a
definition, which is read as ``is defined to be".
We will use notations $\log^b a=(\log a)^b$, $a_+:=a \vee 0$ and $a_-:=(-a) \vee 0$.
For any $x\in \R^d$ and $r>0$, we use $B(x, r)$ to denote the open ball of radius $r$ centered at $x$.
For a Borel subset $V$ in $\R^d$, $|V|$ denotes  the Lebesgue measure of $V$ in $\R^d$,
we use the superscript instead of the subscript for the coordinate of processes as $Y=(Y^1, \dots, Y^d)$.

\section{Resurrection kernel}\label{s:res-ker}

Let $p:\R^d_-\times \R^d_+ \to [0, \infty)$ be a function such that, for each $z\in \R^d_-$, $p(z, \cdot)$ is a probability density on $\R^d_+$, that is,  $\int_{\R^d_+}p(z,y)\, dy=1$. 
Recall that $j(x,z)=\sA(d, \alpha)|x-z|^{-d-\alpha}$, $\alpha\in (0,2)$. 
Let
$$
q(x,y):=\int_{\R^d_-}j(x,z)p(z,y)\, dz\,, \quad x,y\in \R^d_+,
$$
and define a resurrected process 
on
$\R^d_+$ with jump kernel 
$J(x,y)=j(x, y)+q(x,y)$. 
The idea is that when an isotropic $\alpha$-stable process exits $\R^d_+$ by jumping to $z\in \R^d_-$, it is immediately returned to $y\in \R^d_+$ according to the probability distribution $p(z,y)dy$. Therefore we call $p(z,y)$ a \emph{return kernel}. The kernel $q(x,y)$, which we call a \emph{resurrection kernel}, introduces additional jumps from $x$ to $y$, 
thus,  the jump kernel of the resurrected process should be $J(x,y)=j(x,y) + q(x,y)$. 
The process can be constructed via 
Meyer's construction in  \cite{M}, 
or, in case of symmetric $q(x,y)$, by using Dirichlet form theory.
Since $p(z,\cdot)$ is a probability density, an application of Fubini's theorem gives that
$
\int_{\R^d_+}q(x,y)dy 
=\int_{\R^d_-}j(x,z)dz <\infty$, $x\in \R^d_+$.

We would like the resurrected process to be symmetric, to have the scaling property and to be invariant with respect to horizontal translation.
Since  $j(x,z)=\sA(d, \alpha)|x-z|^{-d-\alpha}$, 
the above properties will follow from the symmetry of $q$, 
the homogeneity of $q$:
\begin{align}\label{e:whsc}
q(\lambda x, \lambda y)=\lambda^{-d-\alpha} q(x,y)\, , \qquad \lambda >0, x,y\in \R^d_{ + },
\end{align}
and the horizontal translation invariance (in case $d\ge 2$) of $q$:
\begin{align}\label{e:hti}
q(x+(\wt{u},0), y+(\wt{u},0))=q(x,y), \quad \wt{u}\in \R^{d-1}.
\end{align}
This will depend on properties of the probability kernel $p(z,y)$.
We now recall the examples from the introduction. 

\begin{example}\label{ex:special-return}{\rm
(a) For the trace process of an isotopic $\alpha$-stable process on $\R^d_+$, 
\begin{align}
\label{e:Pst}
p(z,y)
=c\frac{|z_d|^{\alpha/2}}{y_d^{\alpha/2}}|z-y|^{-d}=
c
{|z_d|^{\alpha}}
\left(\frac{|y-z|^2}{y_d|z_d|}\right)^{\alpha/2}
|y-z|^{-d-\alpha}, 
\end{align}
is the Poisson kernel for $\R^d_-$. 
The formula \eqref{e:Pst} can be derived from the Poisson kernel for balls, 
see \cite{BY, B}.
From \eqref{e:poisson-100} and \eqref{e:rsnew1} we see that the corresponding resurrection kernel $q(x,y)$ is symmetric, and from \eqref{e:Pst} that it satisfies \eqref{e:whsc} and \eqref{e:hti}.

\noindent
(b) 
For the process studied in \cite{DROV, Von21},
$$
p(z,y)=\frac{j(z,y)}{\mu(z)}, \quad \text{where } \mu(z)=\int_{\R^d_+}j(z,y)dy.
$$
Clearly, the corresponding resurrection kernel  $q(x,y)$ is symmetric. Since $\mu(z)=c^{-1}|z_d|^{-\alpha}$, we get that
$$
p(z,y)=c |z_d|^{\alpha}|z-y|^{-d-\alpha}.
$$
So $q(x,y)$ satisfies \eqref{e:whsc} and \eqref{e:hti}.
}
\end{example}

Motivated by these two examples, we now introduce a very general return kernel $p(z,y)$.
Let $\gamma_1, \gamma_2$ be two constants such that
$-\infty<\gamma_1\le \gamma_2< 1\wedge\alpha$.
Let $\Psi$ be a positive function on $[2, \infty)$ 
satisfying  the following weak scaling condition: 
There exist  constants $C_1, C_2>0$ such that
$$
C_1(R/r)^{\gamma_1}\leq\frac{\Psi(R)}{\Psi(r)}\leq C_2
   (R/r)^{\gamma_2},\quad 
2\le r<R <\infty.
$$
For notational convenience, we extend the domain of $\Psi$ to $[0, \infty)$ by letting 
$\Psi(t)  \equiv \Psi(2)>0$ on $[0, 2)$.
In particular, 
for any $\delta>0$, there exist 
$\wt{C}_1, \wt{C}_2$
depending on $\delta$ such that
$$
   \wt{C}_1(R/r)^{\gamma_1}\leq\frac{\Psi(R)}{\Psi(r)}\leq \wt{C}_2(R/r)^{\gamma_2},\quad \delta\le r<R <\infty,
$$
and 
\begin{equation}\label{e:Psi-infty0}
    \wt{C}_1(R/r)^{-{\gamma_1}_-}\leq\frac{\Psi(R)}{\Psi(r)}\leq \wt{C}_2(R/r)^{{\gamma_2}_+},\quad \delta\le r<R <\infty.
\end{equation}

Observe that 
after the change of variables
$\wt u =u_d \wt v$ (when $d\ge 2$), 
we see that
\begin{align*}
A&:=
\int_{\R^d_+}
\frac{
\Psi({( |\wt{u}|^2+(u_d+1)^2})/{u_d})
}{
\left( |\wt{u}|^2+(u_d+1)^2\right)^{(d+\alpha)/2}}du 
\nn\\
&\le c \Psi(2) 
\int_{\R^{d-1}}\int_0^\infty
\frac{
{u_d}^{-{\gamma_2}_+}}
{( |\wt{u}|^2+(u_d+1)^2)^{(d+\alpha)/2-
{\gamma_2}_+}
}du_dd\wt{u}\nn\\
&\le c 
\int_0^1
\frac{du_d}{{u_d}^{{\gamma_2}_+}}
\int_{\R^{d-1}}\frac{d\wt{u}}
{(|\wt{u}|+1)^{d+\alpha-2{\gamma_2}_+}}+c
\int_{\R^{d-1}}\int_1^\infty
\frac{{u_d}^{-{\gamma_2}_+}du_dd\wt{u}}
{( |\wt{u}|^2+u_d^2)^{(d+\alpha)/2-{\gamma_2}_+}}
\nn\\
&\le c 
+c \int_1^\infty 
u_d^{{\gamma_2}_+-\alpha-1}
du_d
\int_{\R^{d-1}}\frac{d\wt{v}}
{(|\wt{v}|+1)^{d+\alpha-2{\gamma_2}_+}}
<\infty.
\end{align*}
In the second line we used \eqref{e:Psi-infty0} and, in the last inequality we used
the fact $0 \le{\gamma_2}_+< 1\wedge\alpha$.

Note that for $y\in \R^d_+$ and $z\in \R^d_-$,
${|y-z|^2}/(y_d|z_d|) \ge(y_d+|z_d|)^2/(y_d|z_d|) \ge 2.$
For $y\in \R^d_+$ and $z\in \R^d_-$, define
\begin{equation}\label{e:return-kernel}
\wt p(z, y)
:=
{|z_d|^{\alpha}}
\Psi\left(\frac{|y-z|^2}{y_d|z_d|}\right)
|y-z|^{-d-\alpha}.
\end{equation}
It is easy to see that 
\noindent
{\rm (a)} $\wt p(\lambda z, \lambda y)=\lambda^{-d}\wt p(z,y)$ for all $\lambda >0$, $z\in \R^d_-$, $y\in \R^d_+$;
\noindent
{\rm (b)} $\wt p(z+(\wt{u},0), y+(\wt{u},0))=\wt p(z,y)$ for all $\wt{u}\in \R^{d-1}$,  
$z\in \R^d_-$, $y\in \R^d_+$;
\noindent
{\rm (c)} 
There exists $c\ge 1$ such that for all $r>0$ and $y_0\in \R^{d}_+$  with $B(y_0, 2r)\subset \R^{d}_+$ and $ y_1,y_2\in B(y_0,r),$
\begin{align}
 \label{e:ccc}
 c^{-1}\wt p(w,y_1)\le \wt p(w, y_2) \le c \wt p(w,y_1) \quad \text{for all } w\in \R^d_- .
 \end{align}
Moreover, by the change of variables
$u=|z_d|^{-1}(\wt{y}-\wt{z}, y_d)$ we also have the property:
(d) 
$
\int_{\R^d_+}\wt p(z,y)dy=A$ for all $z \in \R^d_-$.

Thus 
$p(z,\cdot):=A^{-1}\wt p(z,\cdot)$
is a probability density. When $\Psi(t)=t^{\alpha/2}$, $t\ge 2$, we recover the return kernel from Example \ref{ex:special-return}(a), while $\Psi(t)=1$ gives the return kernel in Example \ref{ex:special-return}(b).

With the $\wt p(z, y)$  defined in \eqref{e:return-kernel}, $q(x, y)$ can be written as
\begin{align}
\label{e:rerq}
q(x,y)
=
 \sC \int_{\R^d_-} 
\Psi\left(\frac{|y-z|^2}{y_d|z_d|}\right)\frac{
|z_d|^{\alpha}
}{ |x-z|^{d+\alpha}|y-z|^{d+\alpha}}dz, \quad 
x, y\in \R^d_+ ,
\end{align}
where $\sC:=\sA(d, \alpha)A^{-1}$.
From the properties (a)-(b) above 
we have \eqref{e:whsc} and 
\eqref{e:hti}.

In the next result, we show that the kernel $q$ in   \eqref{e:rerq} is symmetric.

\begin{prop}\label{p:symm}
The resurrection kernel $q$ is symmetric.
\end{prop}
\pf Assume that $d\ge 2$, the proof for $d=1$ being much easier.
Let $x$ and $y$ be any two points in $\R^d_+$.
If $x_d= y_d$, by the change of variables 
$\wt x-\wt z=\wt w-\wt y$ and $w_d=z_d$
we see that 
\begin{align*}
&q(x,y) 
= 
\sC
 \int_{\R^{d-1}} \int_{-\infty}^0\Psi\left(\frac{|\wt y-\wt z|^2+|y_d-z_d|^2}{y_d|z_d|}\right)
\frac{{|z_d|^{\alpha}d\wt{z}\, dz_d}}{ |x-z|^{d+\alpha}|y-z|^{d+\alpha}} \\
&= 
\sC
\int_{\R^{d-1}} \int_{-\infty}^0\Psi\left(\frac{
|\wt x-\wt w|^2+|x_d-w_d|^2}{x_d|w_d|}\right)
\frac{
{|w_d|^{\alpha}d \wt w\,  dw_d}
}{ |y-w|^{d+\alpha}|x-w|^{d+\alpha}}=q(y, x).
\end{align*}
For the remainder of the proof, we assume $x_d\neq y_d$. Without loss of generality, 
we assume that $x=(x_1, \wh 0, x_d)$, $y=(y_1, \wh 0, y_d)$,
and that the line connecting $x$ and $y$ intersects the hyperplane $z_d=0$ at the origin. 
Then
$$
0=\frac{y_d x_1 -x_d y_1}{y_d-x_d},  
\quad x_1=\frac{y_1-x_1}{y_d-x_d}x_d,
\quad y_1=\frac{y_1-x_1}{y_d-x_d}y_d.
$$

For $r>0$, we define
$Tz=r^2z/|z|^2$.
We choose $r$ so that $Tx=y$ and $Ty=x$, i.e.,
$
\frac{r^2}{|x|^2}(x_1, \wh 0, x_d)=(y_1, \wh 0, y_d).
$
Thus
\begin{align}
\frac{|x|^2}{r^2}
=\frac{x_d}{y_d} \quad \text{and}  \quad\frac{|y|^2}{r^2}
=\frac{r^2}{|x|^4} (x_1^2+x_d^2)=\frac{r^2}{|x|^2} 
=\frac{y_d}{x_d}.
\label{e:ss1}
\end{align}

We now fix this $r$. We also write $Tz$ as $z^*$. Then $Tx=y$ and  $Ty=x$.
We have $z_d=r^2z^*_d/|z^*|^2$ and
\begin{equation}\label{e:ss2}
|x-z|=|Ty-Tz^*|=\frac{r^2|y-z^*|}{|y||z^*|}, \quad 
|y-z|=|Tx-Tz^*|=\frac{r^2|x-z^*|}{|x||z^*|}.
\end{equation}
Hence, by \eqref{e:ss1}, \eqref{e:ss2} and the fact $z_d=r^2z^*_d/|z^*|^2$,
\begin{align*}
\frac{|y-z|^2}{y_d|z_d|}=\frac{r^2|x-z^*|^2}{|x|^2y_d|z^*_d|}=
\frac{r^2x_d}{|x|^2y_d} \frac{|x-z^*|^2}{x_d|z^*_d|}=\frac{|x-z^*|^2}{x_d|z^*_d|},
\end{align*}
and  
\begin{align*}
&|x-z|^{-d-\alpha}{|z_d|^{\alpha}}|z-y|^{-d-\alpha}
=
|y-z^*|^{-d-\alpha}|z^*_d|^{\alpha}|z^*-x|^{-d-\alpha}r^{-2d}
|z^*|^{2d}.
\end{align*}
Note that
$
|{\rm det}(JTz)|={r^{2d}}/{|z|^{2d}}.
$
Consequently
\begin{align*}
q(x, y)
&= 
\sC
\int_{\R^d_-}|y-z^*|^{-d-\alpha}
\Psi\left(\frac{|x-z^*|^2}{x_d|z_d^*|}\right)
\frac{
{|z_d^*|^{\alpha}}
}{|x-z^*|^{d+\alpha}}
r^{-2d}
|z^*|^{2d}dz\\
&= 
\sC
\int_{\R^d_-}|y-z^*|^{-d-\alpha}
\Psi\left(\frac{|x-z^*|^2}{x_d|z_d^*|}\right)
\frac{
{|z_d^*|^{\alpha}}
}{|x-z^*|^{d+\alpha}}dz^*=q(y, x).
\quad \ \ \,  \Box
\end{align*}

Define
$$
\Psi_1(u):=\int^u_1\frac{\Psi(v)}{v}dv, \quad u\ge 2.
$$
\begin{lemma}\label{l:phsnew}
{\rm (a)} $\Psi_1 \asymp 1$
when $\gamma_2<0$. {\rm (b)} When $\gamma_1>0$, we have $\Psi_1\asymp \Psi$.
{\rm (c)} When $\gamma_2 \ge 0$, 
there exists a constant $C>0$ such that
$$
1\leq\frac{\Psi_1(R)}{\Psi_1(r)}\leq 
 C(R/r)^{\gamma_2}\log(R/r),
 \quad
2\le r<R <\infty.
$$
\end{lemma}
\pf
Since $\Psi_1$ is an increasing function by definition,  clearly, $1\leq\frac{\Psi_1(R)}{\Psi_1(r)}$ for $2\le r<R <\infty$. When $\gamma_2 \ge 0$,
for any $u\ge 2$ and $\lambda \ge 1$, since $\Psi(w)=\Psi(2)$ for all $w\in [0, 1]$, we have
\begin{align*}
&\Psi_1(\lambda u)
=\int^{u}_{1/\lambda}\frac{\Psi(\lambda w)}{w}dw\le c\lambda^{\gamma_2}\int^{u}_{1/\lambda}\frac{\Psi(w)}{w}dw
\le c \lambda^{\gamma_2}(\Psi_1(u)
+\int^{1}_{1/\lambda}\frac{1}{w}dw)\\
&= c \lambda^{\gamma_2}(\Psi_1(u)
+\log\lambda)
\le  c \lambda^{\gamma_2}(\Psi_1(u)
+(\Psi_1(u)/\Psi_1(2))\log\lambda)\le c\Psi_1(u)\lambda^{\gamma_2}\log\lambda.
\end{align*}
When $\gamma_2<0$, for any $u>2$, since $\Psi(v)=\Psi(2)$ for all $v\in [0, 1]$, we have
\begin{align*}
 \Psi_1(2) \le \Psi_1(u)=\int^u_1\frac{\Psi(v)}{v}dv\le c\Psi_1(2)\int^{u}_1\frac{dv}{v^{1-\gamma_2}}\le c\Psi_1(2)\int^{\infty}_1\frac{dv}{v^{1-\gamma_2}}\le c\Psi_1(2).
\end{align*}
 If $\gamma_1>0$ we have that for $u\ge 2$,
  \begin{align*}
 &\Psi (u) 
\asymp \Psi(u) u^{-\gamma_2} \int_{1}^{u}
 v^{-1+\gamma_2} dv \le c \Psi(u) \int_{1}^{u}
 v^{-1} \frac{\Psi(v)}{\Psi(u)}
  dv =c \Psi_1(u) \\
\qquad \  
   &= c \Psi(u) \int_{1}^{u}
 v^{-1} \frac{\Psi(v)}{\Psi(u)}
  dv \le c \Psi(u) u^{-\gamma_1} \int_{1}^{u}
 v^{-1+\gamma_1} dv 
   \asymp \Psi (u) .
   \qquad \quad \Box
 \end{align*}

We now state the main result of this section -- sharp two-sided estimates for the resurrection kernel $q(x,y)$ and the jump kernel $J(x,y)$. Since the proof of this result is quite technical and long, we postpone it to Section \ref{s:proof-res-ker}.

\begin{thm}\label{t:estimate}
Let $x, y \in \R^d_+$.

\noindent
{\rm (a)} For $x_d \wedge y_d >|x-y|$, it holds that
\begin{align}\label{e:whks1}
q(x,y) &\asymp 
(x_d \wedge y_d)^{-d-\alpha}\asymp(x_d \vee y_d)^{-d-\alpha} 
\end{align}
and
\begin{align}\label{e:whks2}
\sB(x,y)-1
&  \asymp
\left(\frac{|x-y|}{x_d \wedge y_d} \right)^{d+\alpha}.
\end{align}

\noindent
{\rm (b)} For  $x_d \wedge y_d  \le |x-y|$,  it holds that
\begin{align}\label{e:whks3}
q(x,y) &\asymp J(x,y) \asymp |x-y|^{-d-\alpha}
\Psi_1\left(\frac{|x-y|^2}{x_dy_d}\right).
\end{align}
In particular, if $\gamma_1>0$,  then
  \begin{align}\label{e:whks4}
q(x,y) \asymp J(x,y)  \asymp |x-y|^{-d-\alpha} \Psi\big(\frac{|x-y|^2}{x_dy_d}\big) \quad \text{ for } x_d \wedge y_d  \le |x-y|,\end{align}
and, if $\gamma_2<0$, then
  \begin{align}\label{e:whks5}
q(x,y) \asymp J(x,y) \asymp |x-y|^{-d-\alpha}  \quad\text{ for } x_d \wedge y_d  \le|x-y|.
\end{align}
\end{thm}

Recall that we can write 
$J(x,y)=j(x,y)\sB(x,y)$ with $\sB(x,y):=1+q(x,y)/j(x,y)$. 
We have already shown in \eqref{e:whsc} and 
\eqref{e:hti} that the function $\sB(x, y)$ satisfies \textbf{(A4)}.
Note that, by Proposition \ref{p:symm}, \eqref{e:ccc}  and  \cite[Lemma 7.2]{KSV22b}, \textbf{(A1)}-\textbf{(A2)} hold.
Moreover, combining Theorem \ref{t:estimate} with Lemma \ref{l:phsnew}(c), we now see \textbf{(A3)} holds too. Therefore, the resurrected process with the resurrection kernel \eqref{e:rerq} satisfies \textbf{(A1)}--\textbf{(A4)}.

 From Theorem \ref{t:estimate}(b), we also see that if $\gamma_1>0$,  then the function $\Psi$ and the function $\Phi$ from \eqref{e:Phi-infty} can be taken to be the same. The next corollary, which  is an immediate consequence of  the theorem above, 
shows that the functions $\Psi$ and $\Phi$ may not be the same in general.

\begin{corollary}\label{c:estimate}
Let $\gamma\in (-\infty, 1\wedge\alpha)$ 
and $\delta \in \R$.
Suppose  $\Psi(t)=t^\gamma\log^\delta t$,  $t\ge 2$, that is, up to a multiplicative constant,
$$
p(z, y)=\frac{|z_d|^{\alpha-\gamma}}{y_d^{\gamma}}\frac{
\log^\delta\left(\frac{|y-z|^2}{y_d|z_d|}\right)
}{|y-z|^{d+\alpha-2\gamma}},
\quad z\in \R^d_-, y\in \R^d_+.
$$
Then for any $x, y \in \R^d_+$ with
$x_d \wedge y_d >|x-y|$, it holds that
$$
q(x,y) \asymp 
(x_d \wedge y_d)^{-d-\alpha} \asymp (x_d \vee y_d)^{-d-\alpha}, \quad \sB(x,y) -1  \asymp
\left(\frac{|x-y|}{x_d \wedge y_d} \right)^{d+\alpha}
$$
and for $x, y \in \R^d_+$ with
$x_d \wedge y_d  \le |x-y|$, it holds that
\begin{align}\label{e:whkss}
&q(x,y) \asymp J(x,y)\nn\\
&\asymp |x-y|^{-d-\alpha}
\begin{cases}
\left(\frac{|x-y|^2}{x_dy_d } \right)^{\gamma}
\log^\delta\Big(\frac{|x-y|^2}{x_dy_d}\Big) &\text{when }\gamma>0;\\  
\log^{\delta+1}\Big(\frac{|x-y|^2}{x_dy_d}\Big)
&\text{when }\delta>-1, \gamma=0; \\ 
\log\Big(e+\log\big(\frac{|x-y|^2}{x_dy_d}\big)\Big)
&\text{when }\delta=-1, \gamma=0;\\
1&\text{when } \delta<-1, \gamma=0;\\
1  &\text{when }\gamma<0.
\end{cases}
\end{align}
\end{corollary}

\begin{remark}\label{r:estimate}
{\rm
(a) When
 $\Psi(t)=t^{\alpha/2}$, $t\ge2$ (so $\gamma=\alpha/2$, $\delta=0$), which corresponds to the trace process, see Example \ref{ex:special-return}(a), we get 
from  \eqref{e:whkss} that 
for $x, y \in \R^d_+$ with
$x_d\wedge y_d\le |x-y|$, 
$$
J(x,y)\asymp q(x,y)
 \asymp|x-y|^{-d-\alpha} \left(\frac{|x-y|^2}{x_d y_d}\right)^{\alpha/2}=
|x-y|^{-d} x_d^{-\alpha/2}y_d^{-\alpha/2} .
$$
This generalizes \cite[Theorem 6.1]{BGPR} to dimensions $1$ and $2$.

\smallskip
\noindent
(b) When  $\Psi(t)=1$ (so $\gamma=\delta=0$), which corresponds to Example \ref{ex:special-return}(b), we get 
from  \eqref{e:whkss} that 
for $x, y \in \R^d_+$ with
$x_d\wedge y_d\le |x-y|$, 
$$
J(x,y)\asymp q(x,y)\asymp 
|x-y|^{-d-\alpha} \log\left(\frac{|x-y|^2}{x_d y_d}\right)\asymp
|x-y|^{-d-\alpha}\log\left(e+\frac{|x-y|}{x_d\wedge y_d}\right).
$$
}
\end{remark}


\section{Consequences of main results of \cite{KSV22b}}\label{s:prelim}

In this section, we recall the main results of \cite{KSV22b} and apply them to our setting. 
The paper \cite{KSV22b} deals with a general proper open set $D\subset \R^d$ with weaker assumptions. 
For readers' convenience, we restate some results in \cite{KSV22b},  that will be needed in this paper, in the present setting. 

Let $d\ge 1$, $\alpha\in (0,2)$, and $J(x,y)=j(x,y)\sB(x,y)$, $x,y\in \R^d_+$, where 
$j(x,y)=j(|x-y|)=\sA(d, \alpha)|x-y|^{-d-\alpha}$ 
and $\sB(x,y)$ satisfies \textbf{(A1)}-\textbf{(A4)}.
We recall the assumptions \textbf{(H1)}-\textbf{(H5)} imposed in
 \cite{KSV22b} in case $D=\R^d_+$. 
The assumption \textbf{(H1)}, respectively \textbf{(H4)}, are precisely \textbf{(A1)}, respectively \textbf{(A2)}. The other three assumptions are:

\smallskip
\noindent
\textbf{(H2)} For any $a\in (0,1)$ there exists 
$C=C(a)\ge 1$ such that for all $x,y\in \R^d_+$ satisfying
 $x_d\wedge y_d\ge a|x-y|$,
 it holds that
$
C^{-1}\le \sB(x,y)\le C.
$ 

\smallskip
\noindent
\textbf{(H3)} For any $a >0$ there exists 
$C=C(a)>0$
such that
$$
\int_{\R^d_+, |y-x|>a x_d} J(x,y)dy \le  
Cx_d^{-\alpha}
\, ,\quad x\in \R^d_+.
$$

\smallskip
\noindent
\textbf{(H5)} For any  $\epsilon \in (0,1)$ there exists 
$C=C(\epsilon)\ge 1$ with the following property: 
For all $x_0\in \R^d_+$ and $r>0$ with $B(x_0, (1+\epsilon)r)\subset \R^d_+$, we have
$$
C^{-1}\sB(x_1,z)\le \sB(x_2,z) \le C\sB(x_1, z)
$$
$\text{for all }x_1,x_2\in B(x_0,r), \, \,z\in \R^d_+\setminus B(x_0, (1+\epsilon)r)\, .$
\smallskip

It is shown in \cite[Section 7]{KSV22b} that \textbf{(A1)}-\textbf{(A4)} imply the assumptions \textbf{(H1)}-\textbf{(H5)}. This allows us to use here all the results proved in \cite{KSV22b}. Note that \textbf{(H5)} immediately implies that for any  $\epsilon \in (0,1)$ there exists $C=C(\epsilon)\ge 1$ with the following property: For all $x_0\in \R^d_+$ and $r>0$ with $B(x_0, (1+\epsilon)r)\subset \R^d_+$, 
it holds that
\begin{equation}\label{e:estimate-J}
C^{-1}J(x_1,z)\le J(x_2,z) \le C J(x_1, z)
\end{equation}
$\text{for all }x_1,x_2\in B(x_0,r), \, \,z\in \R^d_+\setminus B(x_0, (1+\epsilon)r)$, 
see \cite[(1.8)]{KSV22b}.

Recall from Section \ref{s:intro}  that for $\kappa(x)=\kappa x_d^{-\alpha}$, $\kappa\in [0,\infty)$, 
we introduced 
$$
\EE^\kappa(u,v):=\frac12 \int_{\R^d_+}\int_{\R^d_+} (u(x)-u(y))(v(x)-v(y))J(x,y)dydx+ 
 \int_{\R^d_+} u(x)v(x)\kappa(x)dx,
$$
where $u,v:\R^d_+\to \R$. 
Let $\FF^0$ be the closure of $C_c^{\infty}(\R^d_+)$ in $L^2(\R^d_+, dx)$ 
under $\EE^0_1=\EE^0+(\cdot, \cdot)_{L^2(\R^d_+, dx)}$ and let $\sF^\kappa:= 	\wt \sF^0 \cap L^2(\R^d_+, \kappa(x)dx)$,
where $\wt \sF^0$ is the family of all $\sE^0_1$-quasi-continuous functions in $\sF^0$. Then 
$(\EE^0,\FF^0)$ and $(\sE^\kappa,\sF^\kappa)$ are Dirichlet forms on $L^2(\R^d_+,dx)$.
By \cite[Proposition 3.3]{KSV22b} 
there exists a symmetric Hunt process $Y^{\kappa}=((Y_t^{\kappa})_{t\ge 0}, (\P_x)_{x\in \R^d_+})$ associated with $(\EE^{\kappa}, \FF^{\kappa})$ which can start from every point $x\in \R^d$.
By $\zeta^\kappa$ we denote the lifetime of $Y^\kappa$ and define $Y^\kappa_t=\partial$ for $t \ge \zeta^\kappa$, where $\partial$ is a cemetery point added to $\R^d_+$. 

If $D\subset \R^d$ is an open set, 
let $\tau_D:=\inf\{t>0:\, Y^{\kappa}\notin D\}$ be the exit time of $Y^{\kappa}$ from $D$. 
The part process $Y^{\kappa, D}$ is defined by $Y^{\kappa, D}_t=Y^{\kappa}_t$ if $t<\tau_D$, 
and is equal to $\partial$ otherwise. The Dirichlet form of $Y^{\kappa, D}$ is $(\EE^{\kappa}, \FF^{\kappa}_D)$, 
where $\FF^{\kappa}_D=\{u\in \FF^{\kappa}: \, u=0 \textrm{ quasi-everywhere on  } \R^d_+\setminus D\}$. 
Here quasi-everywhere means that the equality holds everywhere except on a set of capacity zero with respect to $Y^{\kappa}$. 

Let $\mathrm{Cap}^{Y^{\kappa,D}}$ and $\mathrm{Cap}^{X^D}$ denote the capacities with respect to the killed processes 
$Y^{\kappa,D}$, and killed isotropic stable process $X^D$ respectively. The following result is proved in 
\cite[Lemma 3.2]{KSV22b}.
Set $d_D:=\mathrm{diam}(D)$ and $\delta_D:=\mathrm{dist}(D, \partial \R^d_+)$.

\begin{lemma}\label{l:df-comparison}
\cite[Lemma 3.2]{KSV22b}
For every $a>0$, there exists $C=C(a)>0$ such that for all relatively compact open subset 
$D$ of $\R^d_+$ with $d_D \le a \delta_D$, and  for any Borel $A\subset D$, 
\begin{equation}\label{e:comp-cap}
C^{-1}\mathrm{Cap}^{Y^{\kappa,D}}(A)\le \mathrm{Cap}^{X^D}(A)\le C \mathrm{Cap}^{Y^{\kappa, D}}(A).
\end{equation}
\end{lemma}

We will also need the following mean exit time estimates.
\begin{prop}\label{p:exit-time-estimate}
\cite[Proposition 5.3]{KSV22b}
{\rm (a)} There exists a constant $C>0$ such that for all $x_0\in \R^d_+$ and   $r>0$
 with $B(x_0,r) \subset \R^d_+$, it holds that
$$
\E_x \tau_{B(x_0,r)}\ge Cr^{\alpha}\, ,
\quad x\in B(x_0,r/3).
$$

\noindent
{\rm (b)} 
For every $\varepsilon >0$,  there exists  
$C=C(\varepsilon)>0$  
such that for all $x_0\in  \R^d_+$ and all $r >0$ satisfying 
$B(x_0, (1+\eps)r)\subset \R^d_+$, it holds that
$$
\E_x \tau_{B(x_0,r)}\le Cr^{\alpha}\, ,\quad x\in B(x_0,r)\, .
$$
\end{prop}

\smallskip
For $f:\R^d\to \R$ and $x\in \R^d_+$,  set
\begin{align}\label{d:Lalpha}
L_{\alpha}^\sB f(x):=
\textrm{p.v.}\int_{\R^d_+}(f(y)-f(x))J(x,y)\, dy\, , 
\end{align}
whenever the principal value integral on the right-hand side makes sense.
Define 
$$
L^\sB f(x):=L_{\alpha}^\sB f(x)- \kappa (x)f(x)\, ,
\quad x\in \R^d_+\, .
$$
By \cite[Proposition 4.2(a)]{KSV22b}, if $f\in C_c^2(\R^d)$, then $L^\sB f$ and $L_\alpha^\sB f$ are well defined for all  $x \in \R^d_+$. 
Moreover, by \cite[Proposition 4.2(a)]{KSV22b} and using the same argument as in 
\cite[Section 8.2]{KSV} (or derive directly from \eqref{d:Lalpha}), we see that for $u\in C_c^2(\R^d)$ with $u\equiv 0$ on $\R^d_-$ and any $r>0$,
\begin{align}\label{e:operator-interpretation}
L^\sB_\alpha f(x)
&=\int_{\R^d_+}(u(y)-u(x)-\nabla u(x) {\bf 1}_{\{ |y-x|<r\}}\cdot (y-x))J (y, x)dy
\nn\\
&\quad+\int_{\R^d_+}\nabla u(x) {\bf 1}_{\{ |y-x|<r\}}\cdot (y-x) j(y,x)(\sB (y, x)-\sB (x, x))dy
\nn\\
&\quad-\sB (x, x)\int_{\R^d_-}\nabla u(x) {\bf 1}_{\{ |y-x|<r\}}\cdot(y-x) j(y,x)dy.
\end{align}
The expression 
\eqref{e:operator-interpretation} was crucially used in the proof of  \cite[Lemma 5.8(a)]{KSV}.
In this paper we also use \eqref{e:operator-interpretation}  to estimate the action of the operator $L^{\sB}$ on suitable test functions (barriers), see the proof of Lemma \ref{l:key-lemma}.

The following Dynkin-type formula is one of the main results of \cite{KSV22b} and will be extremely important in 
this paper.

\begin{thm}\label{t:dynkin}\cite[Theorem 4.8]{KSV22b}
Suppose that $
D\subset \R^d_+$ is a relatively compact open set. 
For any non-negative function $f$ on $\R^d_+$ with
$f\in C^2(\overline{D})$ and any $x\in D$,
$$
\E_x[f(Y^{\kappa}_{\tau_D})]=f(x)+\E_x\int_0^{\tau_D}L^\sB f(Y^{\kappa}_s)ds\, .
$$
\end{thm}

We also need the following Krylov-Safonov type estimate.
Let $T_A$  be the first hitting time to $A$ for $Y^{\kappa}$.
 
\begin{lemma}\label{l:krylov-safonov}\cite[Lemma 5.4]{KSV22b}
For every $\epsilon\in (0,1)$ there exists $C=C(\epsilon)>0$ such that for 
all $x\in \R^d_+$ and $r>0$ with $B(x,(1+3\epsilon)r)\subset  \R^d_+$,
 and any Borel set $A\subset B(x,r)$,
$$
\P_y(T_A<\tau_{B(x,(1+2\epsilon)r)})\ge C\frac{|A|}{|B(x,r)|}\, , \qquad y\in B(x,(1+\epsilon)r).
$$
\end{lemma}

The following scale invariant Harnack inequality is one of the main results in \cite{KSV22b}.

\begin{thm}\label{t:uhp}\cite[Theorem 1.1]{KSV22b}
\noindent
{\rm (a)}
There exists a constant $C>0$ such that 
for any $r\in (0,1]$, $B(x_0, r) \subset \R^d_+$ and any non-negative function $f$ in $\R^d_+$ which is 
harmonic in $B(x_0, r)$ with respect to $Y^{\kappa}$, we have
$$
f(x)\le C f(y), \qquad \text{ for all } x, y\in 
B(x_0, r/2).
$$

\noindent
{\rm (b)}
For any $L>0$, there exists a constant $C=C(L)>0$ such that for 
any $r \in (0, 1]$, 
all $x_1,x_2 \in  \R^d_+ $  with $|x_1-x_2|<Lr$ and $B(x_1,r)\cup B(x_2,r)  \R^d_+ $  
and any non-negative function $f$ in  $\R^d_+$ 
which is  harmonic in $B(x_1,r)\cup B(x_2,r)$ with respect to 
 $Y^{\kappa}$, we have
$$
f(x_2)\le 
Cf(x_1)\, .
$$
\end{thm}

\smallskip
For a Borel function $f:\R^d_+\to \R$, let
$$
Gf(x):=\E_x\int_0^{\zeta^{\kappa}}f(Y_t^{\kappa})\, dt, \quad x\in \R^d_+,
$$
be the Green potential of $f$. It is shown in \cite[Proposition 6.2]{KSV22b} that if $Y^{\kappa}$ is transient, then there exists a symmetric function $G:\R^d_+\times \R^d_+\to [0,\infty]$ which is lower semi-continuous in each variable and finite 
off the diagonal such that for every non-negative Borel $f$,
$$
Gf(x)=\int_{\R^d_+}G(x,y)f(y)\, dy\, .
$$
Moreover, $G(x, \cdot)$ is harmonic with respect to $Y$ in $\R^d_+\setminus \{x\}$ and regular harmonic with respect to 
$Y^{\kappa}$ in $\R^d_+\setminus B(x, \epsilon)$ for any $\epsilon>0$. The function $G(\cdot, \cdot)$ 
is called the Green function of $Y^\kappa$.
Transience of the process $Y^{\kappa}$ is clear in case $\kappa>0$, see \cite[Lemma 6.1]{KSV22b}. For the case $\kappa=0$, see Lemma \ref{l:alpha>1}.

\section{Scaling and consequences}\label{s:scaling}
In this section we discuss scaling, transience in case $\alpha\in (1,2)$ and $\kappa=0$, and the role of the constant $\kappa=C(\alpha, p, \sB)$.

Let $\rF$ be the closure of $C_c^\infty(\oR^d_+)$ in $L^2(\R^d_+,dx)$ under the norm  $\EE^0_1:=\EE^0+(\cdot,\cdot)_{L^2(\R^d_+,dx)}$ .  Then $(\EE^0, \rF)$ is a regular Dirichlet form on $L^2(\R^d_+,dx)$. Let $((\rY_t)_{t \ge 0}, (\P_x)_{\overline\R^d_+\setminus \sN_0})$ be the Hunt process associated with $(\EE^0, \rF)$, where $\sN_0$ is an exceptional set.
We write $(\rP_t)_{t\ge0}$ and  $(P_t^\kappa)_{t\ge 0}$  for the semigroups of $\rY$ and  $Y^\kappa$  respectively.

Let $(\EE^0, \overline{\FF}_{\R^d_+})$ be the part form of 
 $(\EE^0, \overline{\FF})$ on $\R^d_+$.
i.e.,  the form corresponding to the process $\overline{Y}$ killed at 
the exit time
$\tau_{\R^d_+}:=\inf\{t>0:\, \overline{Y}_t\notin \R^d_+\}$. It follows from \cite[Theorem 4.4.3(i)]{FOT} that 
$(\EE^0, \overline{\FF}_{\R^d_+})$
is a regular Dirichlet form on $L^2({\R^d_+},dx)$ and that $C^{\infty}_c(\R^d_+)$ is its core. Hence 
$\overline{\FF}_{\R^d_+}=\FF^0$, 
implying that $\overline{Y}$ killed upon exiting $\R^d_+$ is equal to $Y^0$. Thus   we conclude that $Y^0$ is a subprocess of $\overline{Y}$, 
that the exceptional set 
$\NN_0$ can be taken to be a subset of $\partial \R^d_+,$ 
and
that the lifetime $\zeta^0$ of $Y^0$ can be identified with $\tau_{\R^d_+}$. 
Suppose that for all $x\in \R^d_+$ it holds that $\P_x(\tau_{\R^d_+}=\infty)=1$. Then $(Y^0_t, \P_x, x\in \R^d_+)\stackrel{d}{=} (\overline{Y}_t, \P_x, x\in \R^d_+)$, implying that $\FF^0=\overline{\FF}_{\R^d_+}=\overline{\FF}$. 

For any $r>0$, define processes $\rY^{(r)}$ and  $Y^{\kappa, (r)}$ by $\rY_t^{(r)}:=r\rY_{r^{-\alpha}t}$ and $Y_t^{\kappa, (r)}:=rY^{\kappa}_{r^{-\alpha}t}$. We have the following scaling and horizontal translation invariance properties of $\rY$ and $Y^\kappa$.

\begin{lemma}\label{l:scaling}
{\rm (a)}	For any $\kappa \ge0$, $r>0$ and $x \in \R^d_+$, $(\rY^{ (r)}, \P_{x/r})$  and  $(Y^{\kappa, (r)}, \P_{x/r})$ have the same laws as $(\rY,\P_x)$ and $(Y^\kappa,\P_x)$ respectively.

\noindent 
{\rm (b)}
In case $d\ge 2$, for any $\kappa \ge 0$, 
$\wt z \in \R^{d-1}$ and $x \in \R^d_+$, 	 $(\rY+(\wt z,0), \P_{x-(\wt z,0)})$ and 
$(Y^\kappa+(\wt z,0), \P_{x-(\wt z,0)})$ 
have the same laws as $(\rY,\P_x)$ and $(Y^\kappa,\P_x)$ respectively.

\noindent
{\rm (c)} If $Y^{\kappa}$ is transient, then for all $x,y\in  \R^d_+$, $x\neq y$, and all $r>0$, 
\begin{equation}\label{e:scaling-of-G}
G(x,y)=G\left(\frac{x}{r}, \frac{y}{r}\right)r^{\alpha-d}\, .
\end{equation}
\end{lemma}
\pf Part (a) follows in the same way as in \cite[Lemma 5.1]{KSV} and \cite[Lemma 2.1]{KSV21}, while part (b) is an immediate consequence of \textbf{(A4)}. Part (c) is a direct consequence of  
part (a),
see the proof in \cite[Proposition 2.4]{KSV20}. \qed

\smallskip
The following two results address the case when $\kappa=0$.
\begin{lemma}\label{l:alpha>1}
Suppose $ \alpha \in (1,2)$ and $\kappa=0$. Then
$\FF ^0\not=\overline \FF$ and $\P_x(\zeta^0<\infty)=1$ for all $x \in \R^d_+$.
\end{lemma}
\pf
Take  $u\in  C_{ c }^{\infty}(\overline \R^d_+)$ 
such that $u \ge 1$ on $B(0, 1) \cap \R^d_+$, then 
 $u \notin \FF^0$. In fact, if $u \in \FF^0$, then by 
Hardy's  inequality for censored $\alpha$-stable processes (see \cite{CS, D}), 
$$
\infty > \EE(u,u) \ge c 
 \int_{\R^d_+}\int_{\R^d_+} \frac{(u(x)-u(y))^2}
 {|x-y|^{d+\alpha}}dxdy
 \ge c
\int_{\R^d_+ }  \frac{u(x)^2}{x_d^\alpha}dx \ge c 
\int_{B(0, 1) \cap \R^d_+ }  |x|^{-\alpha}dx=\infty, 
$$
which gives a contradiction. 

The fact that  $\FF^0 \not=\overline \FF$ implies 
that there is a point $x_0  \in \R^d_+$ such that 
$\P_{x_0}(\zeta^0<\infty)>0$. 
Then by the scaling property of $Y^0$ in Lemma \ref{l:scaling}(a), we have that 
$\P_x(\zeta^0<\infty)=\P_{x_0}(\zeta^0<\infty)>0$ for all $x \in \R^d_+$.
Now, by the same argument as in the proof of \cite[Proposition 4.2]{BBC}, we have that  
$\P_x(\zeta^0<\infty)=1$ for all $x \in \R^d_+$.
\qed

Consequently, in case $\alpha>1$, the process $Y^0$ is transient. This is the reason why in the sequel we consider only $\alpha>1$ when there is no killing. The next result says that $Y^0$ dies at the boundary $\partial \R^d_+$ at its lifetime.

\begin{corollary}\label{c:lemma4-1}
Suppose $\alpha\in (1,2)$ and  $\kappa=0$.
\noindent
{\rm (a)} For any $x\in \R^d_+$, 
 $\P_x(
 Y^0_{\zeta^0-}\in \partial \R^d_+)=1$.	
	\noindent
{\rm (b)} There exists a constant
	$n_0 \ge 2$ such that for any $x\in \R^d_+$,  
$\P_x\left(  \tau_{B(x,n_0 x_d)}=\zeta^0     \right)> 1/2$. 
\end{corollary}
\pf Using Lemma \ref{l:scaling}(a), we see that 
$$
\P_x\big(  \tau_{B(x,n x_d)}=\zeta^0\big)=  \P_{(\wt 0, 1)}\big(  \tau_{B((\wt 0, 1) ,n )} =\zeta^0  \big), \quad x \in \R^d_+.
$$
The sequence of events $(\{ \tau_{B((\wt 0, 1) ,n )}=\zeta^0 \})_{n\ge 1}$ is increasing in $n$ and 
\begin{align}\label{e:zetaf}
\cup_{n=1}^\infty\big\{  \tau_{B((\wt 0, 1) ,n )}
=\zeta^0  \big \} =\big\{\zeta^0<\infty\big\}.
\end{align}
Thus, by Lemma  \ref{l:alpha>1} we have
\begin{equation}\label{e:c-lemma4-1}
\lim_{n \to \infty} 
\P_{(\wt 0, 1)}\big(  \tau_{B((\wt 0, 1) ,n )}=\zeta^0  \big)=
\P_{(\wt 0, 1)}\big(\zeta^0<\infty\big)=1.
\end{equation} 
Moreover, since there is no killing inside $\R^d_+$, it holds that  
$\{ \tau_{B((\wt 0, 1) ,n )}=\zeta^0\}\subset \{Y^0_{\zeta^0-}\in \partial \R^d_+\}$ for each $n\ge 1$. 
Thus it follows from \eqref{e:zetaf} and \eqref{e:c-lemma4-1} that 
$\P_{(\wt 0, 1)}(
Y^0_{\zeta^0-}\in \partial \R^d_+)=1$. 
The claim (a) now follows by Lemma \ref{l:scaling} (a) and (b).

For (b), note that by \eqref{e:c-lemma4-1} there exists $n_0 \ge 2$ such that 
$\P_{(\wt 0, 1)}\big(  \tau_{B((\wt 0, 1) ,n_0 )} =\zeta^0  \big)>1/2$.
Therefore,  
$$ 
\qquad 
\P_x\big(  \tau_{B(x,n_0 x_d)}=\zeta^0    \big)=\P_{(\wt 0, 1)}\big(  \tau_{B((\wt 0, 1) ,n_0 )} =\zeta^0
\big)>1/2, \quad x \in \R^d_+.
\qquad \  \Box
$$

\smallskip
Recall the constant $C(\alpha, q, \sB)$ 
from the introduction: Let $\e_d:=(\wt{0}, 1)$. 
For $q \in (-1, \alpha)$,
\begin{align*}
	C(\alpha,q,\sB) =
	\begin{cases}	 
	\int_{\R^{d-1}}  \int_0^1 \frac{(s^q-1)(1-s^{\alpha-q-1})}{(1-s)^{1+\alpha}} \frac{\sB\big( (1-s)\wt u,1) ,s \e_d\big)}{(|\wt u|^2+1)^{(d+\alpha)/2}}	 ds d\wt u,
	 &  d\ge 2 \\
	\int_0^1 \frac{(s^q-1)(1-s^{\alpha-q-1})}{(1-s)^{1+\alpha}} \sB\big( 1 ,s \big)ds,  & d=1.
	\end{cases}
\end{align*}

\begin{lemma}\label{l:LB-on-g0}
\noindent
{\rm (a)}
 For any $q \in (-1+\wt{\beta_2}, \alpha -\wt{\beta_2})$, 
$C(\alpha,q,\sB) \in (-\infty, \infty)$ is well defined. Further,
$C(\alpha,q,\sB)=0$ if and only if $q\in\{0,\alpha-1\}$.
{\rm (b)}
For any $q\in [\alpha-\beta_1, \alpha)$ it holds that $C(\alpha,q,\sB)= \infty$. Moreover,
$\lim_{q \uparrow \alpha-\beta_1} C(\alpha,q,\sB)= \infty$.
\end{lemma}
\pf 
We only give the proof for $d\ge 2$. The case $d=1$ is simpler.

\noindent
(a)
We first choose $\beta_2$ such that
 \eqref{e:Phi-infty} holds and $q\in (-1+\beta_2, \alpha-\beta_2)$.
Due to the sign of $(s^q-1)(1-s^{\alpha-q-1})$, we see that 
\begin{align*}
\label{e:C}
C(\alpha, q, \sB) 
\begin{cases}
\in (0, \infty] &  q\in ( \beta_2-1, (\alpha-1)\wedge 0) \cup((\alpha-1)_+,  \alpha-\beta_2);\\
=0 &  q=0, \alpha-1;\\
\in [-\infty, 0) &  q\in (\alpha-1, 0) \cup (0, \alpha-1).
\end{cases}
\end{align*}
In the rest 
of the proof we assume that $q\not=0$ and $q\not=\alpha-1$.
By \textbf{(A3)},
\begin{eqnarray}
\sB\big( (1-s)\wt u,1) ,s \e_d\big)& \le& c
{\bf 1}_{0<s<\frac{(|\wt u|^2+1)^{1/2}}{(|\wt u|^2+1)^{1/2}+1}} \left(\frac{(1-s)^2(|\wt u|^2+1)}{s}\right)^{\beta_2}+c{\bf 1}_{\frac{(|\wt u|^2+1)^{1/2}}{(|\wt u|^2+1)^{1/2}+1}<s<1}.\label{e:sbb1}
\end{eqnarray}
Note that for $0<s<1/2$,
\begin{align*}
 \frac{(s^q-1)(1-s^{\alpha-q-1})}{s^{\beta_2}(1-s)^{1+\alpha-2\beta_2}}
\asymp \begin{cases}
 s^{-\beta_2+\alpha-q-1},
 & (\alpha-1) \vee 0<q<\alpha-\beta_2;\\
- s^{-\beta_2},
 & 0<q < \alpha-1;\\
 s^{-\beta_2+q}, &\beta_2-1<q < (\alpha-1)\wedge 0;\\
 - s^{-\beta_2+\alpha-1}, & \alpha-1< q< 0,
 \end{cases}
 \end{align*}
 and, for $1/2<s<1$, 
 \begin{align*}
& \frac{(s^q-1)(1-s^{\alpha-q-1})}{s
^{\beta_2}(1-s)^{1+\alpha-2\beta_2}}
 \asymp \begin{cases} (1-s)^{1-\alpha+2\beta_2}
 & q\in ( \beta_2-1, (\alpha-1)\wedge 0) \cup((\alpha-1)_+,  \alpha-\beta_2);\\
- (1-s)^{1-\alpha+2\beta_2} & q\in (\alpha-1, 0) \cup (0, \alpha-1).
 \end{cases}
 \end{align*}
Thus, 
 for $q \in (-1+\beta_2, \alpha -\beta_2)$, 
\begin{align*}
&\int_{\R^{d-1}} \frac{1}{(|\wt u|^2+1)^{(d+\alpha)/2}}  \int_0^1 \frac{|(s^q-1)(1-s^{\alpha-q-1})|}{(1-s)^{1+\alpha}} \sB\big( (1-s)\wt u,1) ,s \e_d\big)ds \,d\wt u\\
&\le c\int_{\R^{d-1}} \frac{1}{(|\wt u|^2+1)^{(d+\alpha-2\beta_2)/2}}
 \int_0^{\frac{(|\wt u|^2+1)^{1/2}}{(|\wt u|^2+1)^{1/2}+1}} \frac{|(s^q-1)(1-s^{\alpha-q-1})|}{s
^{\beta_2}(1-s)^{1+\alpha-2\beta_2}} ds \,d\wt u
 \\
 &\quad +c\int_{\R^{d-1}} \frac{1}{(|\wt u|^2+1)^{(d+\alpha)/2}}
  \int_{\frac{(|\wt u|^2+1)^{1/2}}{(|\wt u|^2+1)^{1/2}+1}}^1 \frac{|(s^q-1)(1-s^{\alpha-q-1})|}{(1-s)^{1+\alpha}} ds \,d\wt u<\infty,
\end{align*}
which implies that 
\begin{align*}
C(\alpha, q, \sB) 
\begin{cases}
\in (0, \infty) &  q\in ( \beta_2-1, (\alpha-1)\wedge 0) \cup((\alpha-1)_+,  \alpha-\beta_2);\\
=0 &  q=0, \alpha-1;\\
\in (-\infty, 0) &  q\in (\alpha-1, 0) \cup (0, \alpha-1).
\end{cases}
\end{align*}
\noindent
(b) 
By \textbf{(A3)}, for  $(\alpha-1) \vee 0<q<\alpha-\beta_1$ and $0<s<1/2$,
$$
 \frac{(s^q-1)(1-s^{\alpha-q-1})}{(1-s)^{1+\alpha}}\sB\big( (1-s)\wt u,1) ,s \e_d\big) \ge c
 s^{-\beta_1+\alpha-q-1}(|\wt u|^2+1)^{\beta_1}.
$$
Thus, 
$$
C(\alpha,q,\sB)\ge c \int_{\R^{d-1}} \frac{d\wt u}{(|\wt u|^2+1)^{(d+\alpha-2\beta_1)/2}}\int_0^{1/2} s^{-\beta_1+\alpha-q-1}ds,
$$
which implies the claim.
\qed

As already mentioned in the introduction, the function $q\mapsto C(\alpha, q, \sB)$ is strictly increasing 
and continuous
on $[(\alpha-1)_+, \alpha-\wt{\beta}_2)$. Consequently, for every $0\le \kappa<\lim_{q \uparrow \alpha-\wt{\beta}_2} C(\alpha,q,\sB)\le \infty$,  there exists a unique $p_\kappa \in [(\alpha-1)_+, \alpha-\wt{\beta}_2)$ such that
\begin{equation}\label{e:killing-potential}
	\kappa=C(\alpha, p_\kappa, \sB).
\end{equation}
{\it
In the rest of this paper, 
unless explicitly mentioned otherwise, 
we will fix 
$$
\kappa \in [0,\lim_{q \uparrow \alpha-\beta_2} C(\alpha,q,\sB)),
$$	
and assume $\alpha>1$ if  $\kappa=0$. 
Moreover, we omit the superscript $\kappa$ from the notation, i.e., write 
$Y^D$, $\tau_D$  and $\zeta$ instead of 
$Y^{\kappa,D}$, $\tau^\kappa_D$   and $\zeta^\kappa$ 
respectively.
Also, we denote by $p$ the constant $p_\kappa$ in \eqref{e:killing-potential}.
}

The connection between $p$ and $C(\alpha, p, \sB)$ is explained in the following result which is an analog of \cite[(5.4)]{BBC}.
For $q>0$, let $g_q(y):=y_d^q=\delta_{\R^d_+}(y)^q$.

\begin{lemma}\label{l:LB-on-g}
Let $p\in (\wt{\beta}_2-1, \alpha-\wt{\beta}_2)$. Then
$$
L_{\alpha}^\sB g_p(x)= C(\alpha,p, \sB)x_d^{p-\alpha}, \quad x\in \R^d_+.
$$
\end{lemma}
\pf 
We only give the proof for $d \ge 2$. The case $d=1$ is simpler.
Recall  $\mathbf{e}_d=(\tilde{0}, 1)$. By  \textbf{(A4)},  we can for simplicity take $x=(\wt{0}, x_d)$.
Fix $x=(\wt{0}, x_d) \in \R^d_+$ and let  
$\eps \in  (0, (x_d \wedge 1)/2]$.
Let 
$$
I_1(\eps)
:=\int_{\R^d_+, |\wt{z}|^2+| z_d-1|^2> (\varepsilon/x_d)^2}\frac{z_d^p-1}{| (\wt{z}, z_d)-\mathbf{e}_d|^{d+\alpha}}\sB(\mathbf{e}_d,  (\wt{z}, z_d))\, 
dz_d d \wt{z}.$$
We see, by the change of variables $y=x_d z$ and \textbf{(A4)}, that
$
L_{\alpha}^\sB g_p(x)=x_d^{p-\alpha} \lim_{\eps \to 0} I_1(\eps)\, .
$
Using the change of variables $\wt{z}=(z_d-1)\wt{u}$, we get
\begin{align*}
I_1(\eps)&=
 \int_{\R^d_+, | z_d-1|^2 |\wt{u}|^2+| z_d-1|^2> (\varepsilon/x_d)^2}
\frac{z_d^p-1}{|z_d-1|^{1+\alpha}} 
 \, \frac{\sB\big(\mathbf{e}_d, ((z_d-1)\wt{u}, z_d)\big)}{(|\wt{u}|^2+1)^{(d+\alpha)/2}}dz_dd\wt{u}\\
&=
 \int_{\R^{d-1}}
 I_2(\eps, \wt{u})
 \, \frac{d\wt{u}}{(|\wt{u}|^2+1)^{(d+\alpha)/2}},
\end{align*}
where
$$I_2(\eps, \wt{u})=
 \left(\int_0^{1-(\eps/x_d)(|\wt{u}|^2+1)^{-1/2} }+\int_{1+(\eps/x_d)(|\wt{u}|^2+1)^{-1/2}}^{\infty}\right) 
\frac{z_d^p-1}{|z_d-1|^{1+\alpha}} 
\sB\big(\mathbf{e}_d, ((z_d-1)\wt{u}, z_d)\big)dz_d .
$$
Fix $\wt{u}$ and let $\epsilon_0=(\eps/x_d)(|\wt{u}|^2+1)^{-1/2} $. 
Using the change of variables $s=1/z_d$, \textbf{(A1)} and \textbf{(A4)}, by 
the same argument as that in  the proof of \cite[Lemma 5.4]{KSV} ($\wt{z}=|z_d-1|\wt{u}$ there should be $\wt{z}=(z_d-1)\wt{u}$), 
we have that
$
I_2(\eps, \wt{u})=I_{21}(\eps, \wt{u})+I_{22}(\eps, \wt{u})
$
where
\begin{align*}
I_{21}(\eps, \wt{u})&:=\int_0^{1-\epsilon_0}\frac{(s^p-1)+(s^{\alpha-1-p}-s^{\alpha-1})}{(1-s)^{1+\alpha}}\sB\big(((1-s)\wt{u}, 1), s\mathbf{e}_d \big)\, ds, \\
I_{22}(\eps, \wt{u}) &:=\int_{1-\epsilon_0}^{\frac{1}{1+\epsilon_0}}\frac{
s^{\alpha-1-p}(1-s^p)
}{(1-s)^{1+\alpha}}\sB\big(((1-s)\wt{u}, 1), s\mathbf{e}_d \big)\, ds.
\end{align*}
From the proof of Lemma \ref{l:LB-on-g0}, 
we see that $I_{21}(\eps, \wt{u})$ is bounded and 
\begin{align}
\label{e:I21n}
\lim_{\eps\to 0}I_{21}(\eps, \wt{u})=\int_0^1 \frac{(s^p-1)(1-s^{\alpha-p-1})}{(1-s)^{1+\alpha}}\sB\big(((1-s)\wt{u}, 1), s\mathbf{e}_d \big)\, ds\,  .
\end{align}
On the other hand, by \eqref{e:sbb1},  $\sB\big((1-s)\wt{u}, 1), s\mathbf{e}_d \big)$ is bounded by a positive constant  when $\frac{(|\wt u|^2+1)^{1/2}}{(|\wt u|^2+1)^{1/2}+1}<s<1$. 
Since $x_d/\eps \ge 2 \ge 1+(|\wt{u}|^2+1)^{-1/2}$ for $\eps \in  (0, x_d/2]$, we have that for $\eps \in  (0, x_d/2]$, 
$$1-\epsilon_0=1-\frac{\eps/x_d}{(|\wt{u}|^2+1)^{1/2}} \ge
1-\frac{1}{(|\wt u|^2+1)^{1/2}+1}= \frac{(|\wt u|^2+1)^{1/2}}{(|\wt u|^2+1)^{1/2}+1}.$$
Therefore, 
using the facts that
$\epsilon_0 \le 1/2$ and $\frac{1}{1+\epsilon_0} \le 1-\epsilon_0+\epsilon_0^2<1$,
we have
$$\left| I_{22}(\eps, \wt{u})\right| \le 
c\int_{1-\epsilon_0}^{1-\epsilon_0+\epsilon_0^2} \frac{
1-s^p
}{(1-s)^{1+\alpha}}ds \le c \epsilon_0^{2-\alpha}, \quad \eps \in  (0, (x_d \wedge 1)/2], 
$$
(cf., \cite[p.121]{BBC}) 
which implies that $\lim_{\eps\to 0}I_{22}(\eps, \wt{u})=0.$
Therefore, $I_2(\eps, \wt{u})$  is bounded on $(0, (x_d \wedge 1)/2]$ and $\lim_{\eps\to 0}I_{2}(\eps, \wt{u})=\lim_{\eps\to 0}I_{21}(\eps, \wt{u})$. We conclude that 
\begin{align*}
\lim_{\eps\to 0} I_1(\eps)=
 \int_{\R^{d-1}}  \int_0^1 \frac{(s^p-1)(1-s^{\alpha-p-1})}{(1-s)^{1+\alpha}} \frac{\sB\big( (1-s)\wt u,1) ,s \e_d\big)}{(|\wt u|^2+1)^{(d+\alpha)/2}}	 ds d\wt u
=C(\alpha, p, \sB).
\end{align*}
\qed

An immediate, but important, consequence of this lemma is the fact that for $p\in (\wt{\beta}_2-1, \alpha-\wt{\beta}_2)$,
$$
L^\sB g_p(x)= L_{\alpha}^\sB g_p(x)-\kappa(x)x_d^p= L_{\alpha}^\sB g_p(x)-C(\alpha,p, \sB)x_d^{p-\alpha}=0,
$$
for all $x\in \R^d_+$. Thus, the operator $L^\sB$ annihilates the function $x_d^p$.

\section{Dynkin's formula and some estimates}\label{s:dynkin}
Recall that $D_{\wt{w}}(a,b)$ was defined in \eqref{e:box}.
Without loss of generality, we will mostly deal with the case $\wt w=\wt{0}$. 
We will write $D(a,b)$ for $D_{\wt{0}}(a,b)$ and
$U(r)=D_{\wt{0}}(\frac{r}2, \frac{r}2).$
Further  we use  $U$ for $U(1)$.
In case $d=1$, $U(r)=(0,r/2)$.

In the rest of the paper (except Subsection \ref{ss:int1d} and Proposition \ref{p:bound-for-integral-new}, that exclusively deal with the case $d=1$), all the proofs, and even the statements of some lemmas, are given for $d\ge 2$ only. The case $d=1$ is much simpler.

We first recall an important consequence of the L\'evy system formula that will be used repeatedly in this paper, see e.g.~\cite[(3.2), (3.3)]{KSV20}: Let $f:\R^d_+\to [0,\infty)$ be a Borel function, and let $V,W$ be two Borel subsets of $\R^d_+$ with disjoint closures. Then for all $x\in \R^d_+$,
\begin{equation}\label{e:levy-system}
\E_x [f(Y_{\tau_V}), Y_{\tau_V}\in W]=\E_x \int_0^{\tau_V}\int_W f(y)J(Y_s,y)\, dy\, ds\, .
\end{equation}

The next lemma will be used several times in this paper.
\begin{lemma}\label{l:new-two-estimates}
Let $\beta_2$ be the constant in  \eqref{e:Phi-infty}
and let $q\in [0, \alpha-{\beta}_2)$. There exists
 $C=C(q, \beta_2)>0$ such that for all $0<r\le R<\infty$ and all $y\in U(r)$, 
\begin{align*}
&\int_{\R_+^d, z_d>R/2}\Phi\left(\frac{|z|^2}{y_dz_d}\right)\frac{z_d^q dz}{ |z|^{d+\alpha}}  
+\int_0^R \int_{\R^{d-1}, |\wt z|>R/2}\Phi\left(\frac{|z|^2}{y_dz_d}\right)\frac{z_d^q }{ |z|^{d+\alpha}}\, d\wt zd z_d \le C\Phi\left(\frac{r}{y_d}\right)\frac{ R^{q-\alpha+\beta_2}}{r^{\beta_2}}. 
\end{align*}
\end{lemma}
\pf 
Let $y\in U(r)$. 
By the change of variables $\wt z=z_d \wt u$ and the 
facts that 
$\alpha-2\beta_2>-1$, $\beta_2-\alpha<0$ and $q+\beta_2-\alpha<0$,
\begin{align*}
&\int_{\R_+^d, z_d>R/2}
\Phi\left(\frac{|z|^2}{y_dz_d}\right)\frac{z_d^qdz}{ |z|^{d+\alpha}}
\nn\\&=
c_1 \int_{R/2}^\infty
z_d^{q-1-\alpha}
\int_{\R^{d-1}}
\Phi\left(\frac{(|\wt u|^2+1)z_d}{y_d}\right)
\frac{d \wt u}{(|\wt u|^2+1)^{(d+\alpha)/2}} dz_d\\
&\le 
c_2 \Phi\left(\frac{r}{y_d}\right)
r^{-\beta_2}
\int_{R/2}^\infty
z_d^{q-1-\alpha+\beta_2}
\int_{\R^{d-1}}
\frac{d \wt u}{(|\wt u|^2+1)^{(d+\alpha-2\beta_2)/2}} dz_d
\nn\\&=c_3 \Phi\left(\frac{r}{y_d}\right)
r^{-\beta_2} R^{q-\alpha+\beta_2}.
\end{align*}
On the other hand, using the fact that $q-\beta_2 \ge -\beta_2 >-1$,
\begin{align*}
&\int_0^R \int_{\R^{d-1}, |\wt z|>R/2}
\Phi\left(\frac{|z|^2}{y_dz_d}\right)\frac{z_d^q}{ |z|^{d+\alpha}}\, d\wt zdz_d 
\nn\\&\asymp
\int_0^R \int_{\R^{d-1}, |\wt z|>R/2}
\Phi\left(\frac{|\wt z|^2}{y_dz_d}\right)\frac{z_d^q}{|\wt z|^{d+\alpha}}\, d\wt zdz_d \\
& \le c_4 \Phi\left(\frac{r}{y_d}\right) 
r^{-\beta_2} 
\int_0^R z_d^{q-\beta_2}dz_d  \int_{\R^{d-1}, |\wt z|>R}
\frac{d\wt z}{|\wt z|^{d+\alpha-2\beta_2}}\\
& \le c_5 \Phi\left(\frac{r}{y_d}\right) 
r^{-\beta_2} 
  R^{q-\beta_2+1} \int_{R}^\infty t^{-\alpha+2\beta_2-2} dt \le  c_6 \Phi\left(\frac{r}{y_d}\right) 
r^{-\beta_2}  R^{q-\alpha+\beta_2}.
\end{align*}
This completes the proof of the lemma. \qed

For $q, R>0$, let $h_{q, R}(x)=x_d^q {\bf 1}_{D(R,R)}(x)$, $x \in \R^d_+$.

\begin{lemma}\label{l:estimate-of-L-hat-B}
Suppose that $p\in (0, \alpha-\wt{\beta}_2) \cap [(\alpha-1)_+, \alpha-\wt{\beta}_2)$.
There exists $C>0$ such that for 
any $R >0$,
$$
0>L^\sB  h_{p, R}(z) \ge 
-C  R^{p-\alpha}\Phi(R/z_d), \quad z\in U(R).
$$
\end{lemma}

\pf 
We first choose $\beta_2$ such that \eqref{e:Phi-infty} holds and  $p\in (0, \alpha-\beta_2) \cap [(\alpha-1)_+, \alpha-\beta_2)$. 
Then, since
$\sB (z,y) \asymp \Phi\big( \frac{|y|^2}{y_dz_d} \big)$ for $y\in D(R,R)^c\cap \R^d_+$ and $z\in U(R)$, 
 using Lemma \ref{l:new-two-estimates}, we have that
for $z\in U(R)$,
\begin{align}\label{e:neww1}
\int_{D(R,R)^c\cap \R^d_+}\frac{y_d^{ p} }{|y-z|^{d+\alpha}}\sB (z,y)\, dy
 \le  c(p)  R^{p-\alpha}\Phi(R/z_d).
\end{align}
Let $z\in U (R)$. 
By  Lemma \ref{l:LB-on-g}, $L^\sB  g_p(x)=0$. 
Thus, by \eqref{e:neww1},
\begin{align*}
0>L^\sB  h_{p, R}(z)=-\int_{D(R,R)^c\cap \R^d_+}\frac{y_d^{ p} }{|y-z|^{d+\alpha}}\sB (z,y)\, dy \ge 
-c(p)
R^{p-\alpha}\Phi(R/z_d).
\end{align*}
\qed

Next we extend the Dynkin-type formula in Theorem \ref{t:dynkin} to some not relatively compact open sets.

\begin{prop}\label{p:dynkin-hp}
Let $p\in (0, \alpha-\wt{\beta}_2) \cap [(\alpha-1)_+, \alpha-\wt{\beta}_2)$, 
$R \ge 1$ and $r\le R$. For any $x\in U(r)$ it holds that
\begin{equation}\label{e:dynkin-hp}
\E_x [h_{p, R}(Y_{ \tau_{U(r)}})]=h_{p, R}(x)+\E_x \int_0^{\tau_{U(r)}}  L^\sB h_{p, R}(Y_s)\, ds \, .
\end{equation}
\end{prop}
\pf 
We first choose $\beta_2$ such that \eqref{e:Phi-infty} holds and $p\in (0, \alpha-\beta_2) \cap [(\alpha-1)_+, \alpha-\beta_2)$. For $k\in \N$ let $U_k:=\{w \in U(r): w_d>2^{-k}\}$.
 Then $U_k$ is a relatively open compact subset of $\R^d_+$ and 
$h_{p,R}\in C^2(\overline{U_k})$. 
By Theorem \ref{t:dynkin}, for every $k\in \N$, it holds that
\begin{equation}\label{e:dynkin-hp2}
\E_x [h_{p,R}(Y_{\tau_{U_k}})]=h_{p,R}(x)+\E_x \int_0^{\tau_{U_k}}L^{\sB} h_{p,R}(Y_s)\, ds.
\end{equation}
Since $\tau_{U_k}\to \tau_{U(r)}$, the left-hand side converges to $\E_x [h_{q, R}(Y_{ \tau_{U(r)}})]$ by the dominated convergence theorem.
By Lemma \ref{l:estimate-of-L-hat-B}, 
$ L^\sB h_{p, R}(z)\le 0$ for all $z\in U(r)$. 
Thus we can use use the monotone convergence theorem in the right-hand side of \eqref{e:dynkin-hp2}  and obtain \eqref{e:dynkin-hp}. 
\qed

\begin{lemma}\label{l:upper-bound-for-integral}
Let $p\in (0, \alpha-\wt{\beta}_2) \cap [(\alpha-1)_+, \alpha-\wt{\beta}_2)$.
There exists a constant $C>0$ such that for all $R>0$
and all $x\in U(R)$,
$$
\E_x \int_0^{\tau_{U(R)}}
\Phi(R/Y_t^{ d})
\, dt \le C R^{\alpha-p} x_d^{p}\, .
$$
\end{lemma}
\pf
We first choose $\beta_2$ such that \eqref{e:Phi-infty} holds and $p\in (0, \alpha-\beta_2) \cap [(\alpha-1)_+, \alpha-\beta_2)$.
For $R>0$, let 
$C(R):=(D(R,R)\setminus D(3R/4, 3R/4) )  \cap \{y\in \R^d_+: y_d\ge |\wt{y}|\}$.
Note that for $z\in U(R)$ and $y\in C(R)$ we have $|z-y|\le 2|y|\le 2\sqrt{2}y_d$, 
$|z|\le R/\sqrt{2}\le 4|y|/(3\sqrt{2})$, 
and therefore $(\sqrt{2}- (4/3))y_d\le (\sqrt{2}- (4/3))|y|\le  \sqrt{2}(|y|-|z|)  \le \sqrt{2}|z-y| \le 4y_d$. 
Thus, 
\begin{eqnarray*}
\sB (z,y) \asymp
\Phi\left(\frac{|z-y|^2}{z_dy_d} \right) 
\ge c_1 \Phi\left(\frac{|y|}{z_d}\right) 
\ge c_2\Phi\left(\frac{R}{z_d}\right) \left(\frac{|y|}{R}\right)^{\beta_1}.
\end{eqnarray*}
Using 
\eqref{e:levy-system}
we get for  $x\in U(R)$
(with constants $c_3$, $c_4$ independent of $R$),
\begin{align*}
&\E_x[
h_{p, R}
(Y_{\tau_{U(R)}})]\ge \E_x[
h_{p, R}
(Y_{\tau_{U(R)}}), Y_{\tau_{U(R)}}\in C(R)]\\
&\ge  c_3  \E_x\int_{C(R)} \int_0^{\tau_{U(R)}} 
|Y_t-y|^{-d-\alpha} \
 \left(\frac{|y|}{R}\right)^{\beta_1}
\Phi\left(\frac{R}{Y^{d}_t}\right)  y_d^p \, dy \, dt\\
&\ge  c_4 R^{-\beta_1}\int_{C(R)}{y_d^p}{|y|^{-d-\alpha+\beta_1}}\, dy
\left(\E_x \int_0^{\tau_{U(R)}} 
\Phi\left(\frac{R}{Y_t^{d}}\right)   \, dt\right)\, .
\end{align*}
Note that, 
$$ R^{-\beta_1}\int_{C(R)}y_d^p |y|^{-d-\alpha+\beta_1}\, dy
\asymp  
R^{p-\alpha} \int_{C(1)}z_d^p |z|^{-d-\alpha+\beta_1}\, dz 
\asymp
R^{p-\alpha}.$$ 
Thus, by Proposition \ref{p:dynkin-hp} and  
Lemma \ref{l:estimate-of-L-hat-B}, for all $R>0$,
\begin{align}\label{e:wow}
c R^{p-\alpha} \E_x \int_0^{\tau_{U(R)}}
\Phi\left(\frac{R}{Y_t^{d}}\right) \, dt \le
 \E_x[h_{p, R}(Y_{\tau_{U(R)}})]=x_d^p +
 \E_x\int_0^{\tau_{U(R)}}L^\sB  h_{p, R}
(Y_s)\, ds \le x_d^p\, .
\end{align}
\qed

\begin{corollary}\label{c:exltlow}
Let $p\in (0,\alpha-\wt{\beta}_2) \cap [(\alpha-1)_+, \alpha-\wt{\beta}_2)$. 
Then there exists $C>0$ such that 
$$
\E_x \int_0^{\tau_{U}} \Phi\left(\frac{1}{Y^{ d}_s}\right)ds\le C \P_x\left(Y_{\tau_{U}}\in  D(1, 1)  \right)  \quad \text{for all } x \in U.
$$
\end{corollary}
\pf
This corollary follows from \eqref{e:wow} and  the fact that 
$h_{p, 1}$ is bounded by 1 and supported 
on $D(1, 1)$ so that 
$$
\qquad \, 
\P_x\left(Y_{\tau_{U}}\in D(1, 1)\right)   
\ge \E_x [h_{p, 1}(Y_{ \tau_{U}})] 
\ge c \E_x \int_0^{\tau_{U}} \Phi\left(\frac{1}{Y^{d}_s}\right)ds.
\qquad \qquad \Box
$$

\begin{corollary}\label{c:exit-away}
Let $p\in (0, \alpha-\wt{\beta}_2) \cap [(\alpha-1)_+, \alpha-\wt{\beta}_2)$. There exists $C>0$ such that, for all $r>0$ and $x\in U(r)$, it holds that
$$
\P_x(Y_{\tau_{U(r)}}\notin D(r,r))
\le C\left(\frac{x_d}{r}\right)^p.
$$
\end{corollary}
\pf We first choose $\beta_2$ such that \eqref{e:Phi-infty} holds and $p\in (0, \alpha-\beta_2) \cap [(\alpha-1)_+, \alpha-\beta_2)$.
By scaling in Lemma \ref{l:scaling} (a), it suffices to prove the claim for $r=1/2$.  Let $D=D(1,1)$. For $z\in U$ and $w\in \R^d_+\setminus D$, it holds that $|z-w|\asymp |w|$.
Hence, by \eqref{e:levy-system},
\begin{align*}
\P_x(Y_{\tau_{U}}\notin D)
&=\E_x \int_0^{\tau_U} \int_{\R^d_+\setminus D} 
J(w,Y_t)dw\, dt \le 
c_1 \E_x \int_0^{\tau_U}\int_{\R^d_+\setminus D}|w|^{-d-\alpha} 
\Phi\left(\frac{|w|^2}{w_d Y_t^{d}}\right)dw\, dt.
\end{align*}
It follows from Lemma \ref{l:new-two-estimates} (with $r=1/2$ and $R=2$) that
$$
\int_{\R^d_+\setminus D}|w|^{-d-\alpha} 
\Phi\left(\frac{|w|^2}{w_d Y_t^{d}}\right)dw \le 
c_2 \Phi\left(\frac{1}{Y_t^{d}}\right).
$$
Therefore, by using Lemma \ref{l:upper-bound-for-integral} we get 
\begin{align*}
\qquad \qquad \qquad
\P_x(Y_{\tau_{U}}\notin D)\le 
c_3 \E_x \int_0^{\tau_U} 
\Phi\left(\frac{1}{Y_t^{d}}\right)dt\le c_4 x_d^p. 
\qquad \qquad \qquad  \Box
\end{align*}

\begin{prop}\label{p:POTAe7.14}
Let $\beta_2$ be such that \eqref{e:Phi-infty} holds and let $p\in (0, \alpha-\beta_2) \cap [(\alpha-1)_+, \alpha-\beta_2)$. Then
there exists $C=C(\beta_2)>0$ such that for all 
$0<4r\le R<\infty$ and $w\in D(r,r)$, 
$$
\P_w\Big(Y_{\tau_{B(w, r)\cap \R^d_+}}
\in \R^d_+ \setminus B(w, R)\Big)\le C
\frac{r^{\alpha-\beta_2}}{R^{\alpha-\beta_2}}\frac{w_d^{p}}{r^{p}}.
$$
\end{prop}
\pf
Let $\wt w= \wt0$, 
 $0<4r\le R<\infty$, 
$w\in D(r,r)$  and $y\in B(w,r)\cap \R^d_+$ and $z\in A(w, R, 4)\cap \R^d_+$. 
Then 
 $|z-y| \asymp |z| \asymp |z-w|>R >y_d$.
Thus, 
$$
J (y,z)\asymp \frac{1}{ |y-z|^{d+\alpha}}
\Phi\left(\frac{|z-y|^2}{y_dz_d}\right)
\asymp
\frac{1}{ |z|^{d+\alpha}}
\Phi\left(\frac{|z|^2}{y_dz_d}\right).
$$
Thus by using \eqref{e:levy-system} in the first inequality below and Lemma \ref{l:new-two-estimates} in the last inequality, we get
\begin{align*}
&\P_w\Big(Y_{\tau_{B(w, r)\cap \R^d_+}}\in \R^d_+ \setminus B(w, R)\Big)\le 
c\, \E_w\int^{\tau_{ B(w, r)\cap \R_d^+}}_0   \int_{\R^d_+ \setminus B(w, R)}
\Phi\left(\frac{|z|^2}{Y_t^{d}z_d}\right)
\frac{dzdt}{ |z|^{d+\alpha}}\\
&\le c\, \E_w\int^{\tau_{ B(w, r)\cap \R_d^+}}_0  \int_{\R^d_+\setminus D(R/2, R/2)}\Phi\left(\frac{|z|^2}{Y_t^{d}z_d}\right)\frac{dzdt}{ |z|^{d+\alpha}}\\
&\le c\, r^{-\beta_2} R^{-\alpha+\beta_2}\E_w\int^{\tau_{ B(w, r)\cap \R_d^+}}_0
\Phi\left(\frac{r}{Y_t^{ d}}\right)  dt.
\end{align*}
Since $B(w, r)\cap \R_d^+\subset D(2r,2r)$,
applying
Lemma \ref{l:upper-bound-for-integral}, we get that 
for all  
 $0<4r\le R<\infty$ and $w\in D(r,r)$,
\begin{align*}
\P_w\Big(Y_{\tau_{B(w, r)\cap \R^d_+}}
\in \R^d_+ \setminus B(w, R)\Big)
& \le  c r^{-\beta_2} R^{-\alpha+\beta_2}\E_w\int^{\tau_{ D(2r,2r)}}_0
\Phi\left(\frac{r}{Y_t^{d}}\right)  dt\nn\\
\qquad \qquad \qquad \qquad \qquad  \, 
&\le c \frac{r^{\alpha-\beta_2}}{R^{\alpha-\beta_2}}\frac{w_d^p}{r^p}.
\qquad \qquad \qquad  \qquad \qquad \quad \Box
\end{align*}

\section{The key technical result and exit probability estimates}
\subsection{The key lemma and exit probability estimates}
The following lemma is the key technical result of the paper. It will allow us to obtain 
exit probability estimates essential 
for the proof of
Theorem \ref{t:BHP}.
\begin{lemma}\label{l:key-lemma}
Let $p\in ((\alpha-1)_+, \alpha-\wt{\beta}_2)$.
\noindent
{\rm (a)} There exist a $C^2$-function $\psi:\R^d\to [0,\infty)$ with compact support, and 
a constant $C_{1}>0$ such that
$$
L^\sB  \psi(x)\le C_{1}  \Phi(1/x_d),\quad x\in U,
$$
and the following assertions hold: 

\noindent
{\rm (b)} 
The function 
$\phi(x):=h_{p, 1}(x)-\psi(x)$, $x\in \R^d_+$, 
satisfies the following properties: 
	\begin{itemize}
	\item[(b1)]  $\phi(x)=x_d^p$ for all $x=(\wt{0}, x_d)\in U$ with $0<x_d<1/4$;
	\item[(b2)]  $\phi(x)\le 0$ for all $x\in U^c\cap \R^d_+$;
	\item[(b3)]  There exists $C_{2}>0$ such that $L^\sB \phi(x)\ge -C_{2} \Phi(1/x_d)$ for all $x\in U$.
	\end{itemize}
\end{lemma}

Note that Lemma
\ref{l:key-lemma} has the  stronger assumption  
$p>(\alpha-1)_+$, which requires
the killing function to be strictly positive.
The proof of this lemma is  long and involved, therefore we postpone it to the end of this section. We now prove 
several consequences of Lemma \ref{l:key-lemma}.

\begin{lemma}\label{l:lower-bound-for-integral}
Let $p\in ((\alpha-1)_+, \alpha-\wt{\beta}_2)$.
For any $x=(\wt{0}, x_d)$ with $0<x_d<1/4$, it holds that
\begin{equation}\label{e:lower-bound-for-integral}
\E_x \int_0^{\tau_{U}}  \Phi(1/Y_t^{d})\,  dt \ge C_{2}^{-1} x_d^p\, ,
\end{equation}
where $C_{2}$ is the constant from Lemma \ref{l:key-lemma}.
\end{lemma}
\pf
We first choose $\beta_2$ such that \eqref{e:Phi-infty} holds 
and $p\in ((\alpha-1)_+, \alpha-\beta_2)$. Recall that $\phi=h_p-\psi$.
For $k\in \N$ let 
$U_k:=\{y\in U:\, y_d>2^{-k}\}$. 
Then $U_k$ is a relatively open compact subset of $\R^d_+$ and by Lemma \ref{l:key-lemma}  $\phi\in C^2(\overline{U_k})$. Let $x=(\wt{0}, x_d)$ with $0<x_d<1/4$. By Theorem \ref{t:dynkin} (applied separately to 
$h_{p, 1}$ 
and $\psi$, and then taking the difference), 
for every $k\in \N$ with $2^{-k}<x_d$, it holds that
$$
\E_x [\phi(Y_{\tau_{U_k}})]= \phi(x)+\E_x \int_0^{\tau_{U_k}}L^{\sB} \phi(Y_s)\, ds. 
$$
From Lemma \ref{l:key-lemma} (b3), we know that $
L^\sB \phi(z)\ge -C_{2}\Phi(1/z_d)$  for all $z\in U$. 
Therefore,
\begin{equation}\label{e:new-nn0}
\E_x [\phi(Y_{\tau_{U_k}})]- \phi(x)\ge -C_{2} \E_x \int_0^{\tau_{U_k}}\Phi(1/Y_s^d)\, ds.
\end{equation}
Since $\tau_{U_k}\to \tau_{U}$, by letting $k\to \infty$, 
and using the monotone convergence theorem, 
the right-hand side converges to $-C_{2}\int_0^{\tau_U}\Phi(1/Y_s^d)\, ds$. Since $\psi$ and $h_{p, 1}$ are bounded, 
by the dominated convergence theorem, we have
$\E_x [\psi(Y_{\tau_{U_k}})]\to \E_x [\psi(Y_{\tau_{U}})]$ and
$
\E_x [
h_{p, 1}
(Y_{\tau_{U_k}})] \to \E_x [
h_{p, 1}
(Y_{\tau_{U}})].
$
Hence, by letting $k\to \infty$ in \eqref{e:new-nn0}, and using Lemma \ref{l:key-lemma} (b1)--(b2),  we get 
$$
-x_d^p\ge \E_x[
 \phi(Y_{\tau_U})]-\phi(x)\ge -C_{2}\int_0^{\tau_U}\Phi(1/Y_s^d)\, ds.
$$
This proves \eqref{e:lower-bound-for-integral}. \qed

\begin{lemma}\label{l:exit-probability-estimate}
If $p\in ((\alpha-1)_+, \alpha-\wt{\beta}_2)$
 then 
there exists  $C>0$ such that  for 
$x=(\wt{0}, x_d)\in D(1/8,1/8)$,
$$
\P_x(Y_{\tau_{D(1/4,1/4)}}
\in D(1/4,1)\setminus D(1/4,3/4))
\ge C  x_d^p \, .
$$
\end{lemma}
\pf  
For $y\in D(1/4,1/4)$ 
and $z\in D(1/4,1)\setminus D(1/4,3/4)$, it holds that $y_d<z_d$, $|z|\asymp |y-z|\asymp z_d \asymp 1$ and $y_d<2|y-z|$. Hence, 
$
\sB (y,z) \asymp \Phi\left({1†}/{y_d}\right)
$ 
and, by using \eqref{e:levy-system}
and Lemma \ref{l:scaling}(a), we get that for 
$0<x_d<1/8$, 
\begin{align*}
&\P_{(\wt{0}, x_d)}\left(
Y_{\tau_{D(1/4,1/4)}}\in 
D(1/4, 1)\setminus D(1/4, 3/4)\right)\nonumber\\
&\ge c\,\E_{(\wt{0}, x_d)}
\int^{\tau_{D(1/4,1/4)}}_0 
\Phi(1/Y_t^{d})\, 
\int_{D(1/4, 1)\setminus D(1/4, 3/4)}\frac{dz}{|z|^{d+\alpha}}\, dt\nonumber\\
&\ge c\,\E_{(\wt{0}, x_d)}
\int^{\tau_{D(1/4,1/4)}}_0 
\Phi(1/Y_t^{d})\,\, dt \asymp
\E_{(\wt{0}, 2x_d)}\int^{\tau_{U}}_0 
\Phi(1/Y_t^{d})dt\, .
\end{align*}
The claim now follows from Lemma  \ref{l:lower-bound-for-integral}. \qed

Note that in the next two results, we allow $p=(\alpha-1)_+$.

\begin{lemma}\label{l:POTAl7.4}
Suppose $p\in (0, \alpha-\wt{\beta}_2) \cap [(\alpha-1)_+, \alpha-\wt{\beta}_2)$.
There exists $C>0$  such that for any 
$x \in U(2^{-4}) $,
$$
\P_x\left(Y_{\tau_{U}}\in D(1, 1)\right)\le C
\P_x\left(Y_{\tau_{U}}\in
 D(1/2, 1)\setminus D(1/2, 3/4)\right).
$$
\end{lemma}
\pf
We first choose $\beta_2$ such that
\eqref{e:Phi-infty} holds and $p\in (0, \alpha-\beta_2) \cap [(\alpha-1)_+, \alpha-\beta_2)$.
Let
$$
H_2:=\{Y_{\tau_{U}}\in D(1, 1)\}, \quad H_1:=\{Y_{\tau_{U}}\in D(1/2, 1)\setminus D(1/2, 3/4)\}.
$$
We first note that, by Lemma \ref{l:scaling}(b),
\begin{align}\label{e:POTAe7.151}
\P_w(H_1)&\ge
\P_w(Y_{\tau_{D_{\wt w}(1/4, 1/4)}
}\in D_{\wt w}(1/4, 1)\setminus D_{\wt w}(1/4, 3/4))\nn\\
&=\P_{(\wt 0, w_d)}
(Y_{\tau_{D(1/4, 1/4)}}
\in D(1/4, 1)\setminus D(1/4, 3/4)).
\end{align}
When $p=\alpha-1>0$, we 
choose a $q_* \in (\alpha-1, \alpha-\beta_2)$ and let
$\kappa^*(x)=C(\alpha, q_*,  \sB)x_d^{-\alpha}$.
Let $Y^{\kappa^*}$ be the process associated with Dirichlet form
$\EE(u,v)+\int_{\R^d_+} u(x)v(x)\kappa^*(x) dx$.
 By Lemma \ref{l:exit-probability-estimate}, we get that, when $p=\alpha-1$,
$$
\P_{(\wt 0, w_d)}
(Y^{\kappa^*}_{\tau_{D(1/4, 1/4)}}
\in D(1/4, 1)\setminus D(1/4, 3/4))
 \ge c \, w^{q_*}_d, 
\quad w\in U(1/4).
$$
Thus, by this and \eqref{e:POTAe7.151},
\begin{equation*}
\P_w(H_1)\ge 
\P_{(\wt 0, w_d)}
(Y^{\kappa^*}_{\tau_{D(1/4, 1/4)}}
\in D(1/4, 1)\setminus D(1/4, 3/4))
 \ge  c \,  w^{q_*}_d, \quad 
w\in U(1/4).
\end{equation*}
When $p\in ((\alpha-1)_+, \alpha-\beta_2)$, we just use  \eqref{e:POTAe7.151} and Lemmas \ref{l:exit-probability-estimate} directly to obtain that 
$\P_w(H_1)\ge  c\, w^p_d$,
for $w\in U(1/4)$.

Therefore, we see that for all $p\in (0, \alpha-\beta_2) \cap [(\alpha-1)_+, \alpha-\beta_2)$, there exists 
$q \in ((\alpha-1)_+, \alpha-\beta_2)$ with $q \ge p$ such that 
$
\P_w(H_1)\ge  c w^q_d$ for
$w\in U(1/4).$
Using this and  Proposition \ref{p:POTAe7.14},  the remaining part of the proof closely follows that of \cite[Lemma 6.2]{KSV20} and \cite[Lemma 5.5]{KSV21} (Proposition \ref{p:POTAe7.14} is used in the proof). Therefore we omit the rest of the proof. \qed

The next comparability result 
summarizes the exit probability estimates obtained so far and will play a crucial role in the remainder of  this paper.

\begin{prop}\label{p:allcomparable}
Let $p\in (0,\alpha-\wt{\beta_2}) \cap [(\alpha-1)_+, \alpha-\wt{\beta}_2)$.
For all $r>0$,
$$
\P_x\left(Y_{\tau_{U(r)}}\in  D(r, r)  \right)
\asymp 
\P_x\left(Y_{\tau_{U(r)}}\in \R^d_+\right)\asymp 
\left(\frac{x_d}{r}\right)^p
 \quad  
\text{for all } x\in U(2^{-4}r).
$$
\end{prop}
\pf 
By scaling in Lemma \ref{l:scaling}(a)  it suffices to prove both results for $r=1$.
By Proposition \ref{p:dynkin-hp}, Lemma
\ref{l:estimate-of-L-hat-B} and  the fact that 
$h_{p, 1}$
is bounded by 1 and supported 
on $D(1, 1)$, we have that for every $x\in U(2^{-4})$,
\begin{align*}
&\P_x\left(Y_{\tau_{U}}\in D(1, 1)\right) 
\ge \E_x [h_{p, 1}(Y_{ \tau_{U}})] \\
&=
x_d^p +\E_x\int_0^{\tau_{U}}L^\sB  
h_{p, 1}(Y_s)\, ds 
\ge
x_d^p -c_1\E_x \int_0^{\tau_{U}} 
\Phi\left(\frac{1}{Y^{d}_s}\right)ds.
\end{align*}
Thus, by Corollary \ref{c:exltlow},
\begin{align}\label{e:xdp}
x_d^p \le c_1\E_x \int_0^{\tau_{U}} 
\Phi\left(\frac{1}{Y^{ d}_s}\right)ds +
\P_x\left(Y_{\tau_{U}}\in D(1, 1)\right) 
\le c_2 \P_x\left(Y_{\tau_{U}}\in D(1, 1)\right).
\end{align}

On the other hand, 
for $y\in U$ 
and $z\in D(1/2,1)\setminus D(1/2,3/4)$, it holds that $y_d<z_d$, $|z|\asymp |y-z|\asymp z_d$ and $y_d<2|y-z|$. Hence, 
$
\sB (y,z) \asymp \Phi\left({|z|}/{y_d}\right)
$
and, by using \eqref{e:levy-system} 
and Lemma \ref{l:upper-bound-for-integral}, we get that for every $x\in U(2^{-4})$,
\begin{align*}
\P_x\left(Y_{\tau_{U}}\in 
D(1/2, 1)\setminus D(1/2, 3/4)\right)&\le 
c_3\,\E_x\int^{\tau_{U}}_0 
\Phi(1/Y_t^{d})\, 
\int_{D(1/2, 1)\setminus D(1/2, 3/4)}\frac{dz dt}{|z|^{d+\alpha-\beta_2}}\nonumber\\
&\le c_4\,\E_x\int^{\tau_{U}}_0 
\Phi(1/Y_t^{d})\,\, dt \le c_5 x_d^p.
\end{align*}
Thus, by  Lemma \ref{l:POTAl7.4},
$\P_x\left(Y_{\tau_{U}}\in D(1, 1)\right)\le   c_6 x_d^p$ for every $x\in U(2^{-4})$.
Combining this with \eqref{e:xdp} and Corollary \ref{c:exit-away}, 
we get that for every $x\in U(2^{-4})$,
\begin{align*}
c_2^{-1} x_d^p&\le 
\P_x\left(Y_{\tau_{U}}\in D(1, 1)\right) \le 
\P_x(Y_{\tau_{U}}\in \R^d_+) \nn\\
\qquad &=
\P_x(Y_{\tau_{U}}\in D(1,1))+
\P_x(Y_{\tau_{U}}\in \R^d_+ \setminus D(1,1))\le  c_7 x_d^p.  \qquad \quad \ \Box
\end{align*}

\subsection{
Auxiliary
 lemmas}
In this subsection we give two
lemmas needed in the proof of Lemma \ref{l:key-lemma}. The index $\beta_2$ below is such that $\beta_2<1\wedge \alpha$ and \eqref{e:Phi-infty} holds.

The next lemma is one of the key technical results in this paper.
\begin{lemma}\label{l:1n} 
{\rm (a)} For any $k \in \R$,  
 there exists $C>0$ such that 
for $0<x_d \le R/2$,
\begin{align}\label{e:ln}
& \int_{D_{\wt x}(R, R)\cap \{|y-x| \ge x_d/2\}}\Phi\left(\frac{|x-y|^2}{x_dy_d}\right) \frac{dy}{|x-y|^{d+\alpha-k}} \nn \\
 \le &
C\begin{cases}
R^{k-\alpha}\Phi(R/x_d)(1+{\bf 1}_{\{  k+\beta_1=\alpha  \}}\log (R/x_d)),
& k+\beta_1 \ge \alpha; \\
[R^{-\beta_1}\Phi(R/x_d)  x_d^{k-\alpha+\beta_1}] \wedge\1 [R^{k+\beta_2-\alpha}  x_d^{-\beta_2} ],
& k+\beta_1 < \alpha<k+\beta_2;\\
  x_d^{k-\alpha} (1+{\bf 1}_{\{  k+\beta_2=\alpha  \}} \log (R/x_d)),
& k+\beta_2 \le \alpha\, .
\end{cases}
\end{align}
{\rm (b)} For any $k >\alpha $, there exists 
$C>0$ such that for $0<x_d \le R/2$,
$$
\int_{D_{\wt x}(R, R)}\frac{\sB(x, y)}{|x-y|^{d+\alpha-k}}dy \le   
CR^{k-\alpha} (1\vee \Phi(R/x_d)).
$$
\end{lemma}
\pf
Without loss of generality, 
we assume $\wt x=\wt 0$.
Let $R>0$ and $x_d \le R/2$.  

\noindent 
(a) Define
\begin{align*}
&I(k) :=\int_{D(R, R)\cap\{|x-y|\ge x_d/2\}}\Phi\left(\frac{|x-y|^2}{x_dy_d}\right) \frac{dy}{|x-y|^{d+\alpha-k}}\\
\normal
&=\int_{D(R, x_d/2)} +\int_{D(R, R)\setminus D(R, 3x_d/2)}+\int_{(D(R,3 x_d/2)\setminus D(R, x_d/2)) \cap \{|y-x|\ge x_d/2\}}\\
&=:I_1(k)+I_2(k)+I_3(k).
\end{align*}

\noindent 
(a-i)
Clearly, $y_d<x_d$ for $y_d \in  D(R, x_d/2)$.
Using the change of 
variables $y_d=x_d h$ and $r=x_d s$ in the second line below, we get
\begin{align*}
&I_1(k)
\asymp  \int_0^R r^{d-2} \int_0^{x_d/2} \frac{1}{((x_d-y_d)+r)^{d+\alpha-k}}  
\Phi\big(\frac{((x_d-y_d)+r)^2}{x_dy_d}\big)
dy_d\, dr\\
&= x_d^{-\alpha+k} \int_0^{R/x_d} s^{d-2} \int_0^{1/2} \frac{1}{[(1-h)+s]^{d+\alpha-k}} 
\Phi\big(\frac{[(1-h)+s]^2}{h}\big)
 dh\, ds, 
\end{align*}
which is, 
using $1-h\asymp 1$,
comparable to 
$$
x_d^{-\alpha+k} \int_0^{R/x_d} \frac{s^{d-2}}{(1+s)^{d+\alpha-k}} \int_0^{1/2} \Phi\big(\frac{(1+s)^2}{h}\big)dhds. 
$$
Since 
$
 \int_0^{1/2} \Phi({(1+s)^2}/{h})dh
\le c \Phi((1+s)^2)  \int_0^{1/2} h^{-\beta_2}dh, 
$
we have 
\begin{align}\label{e:l1bd}
I_1(k)
&\le c 
x_d^{-\alpha+k} 
\int_0^{R/x_d} \frac{s^{d-2}  \Phi\big((1+s)^2\big)}{(1+s)^{d+\alpha-k}} ds \le c  x_d^{-\alpha+k} 
\Big(
1+
\int_{1}^{R/x_d} \frac{  \Phi(s^2)}{s^{2+\alpha-k}} ds \Big). 
\end{align}

In order to estimate $I_3(k)$, for $a>0$ we define 
$K_a:=\{y\in \R^d:\, |\wt{y}| < ax_d/2, |y_d-x_d| < ax_d/2\}$. Then $K_{1/\sqrt{d}}\subset B(x, x_d/2)\subset K_1$, hence
\begin{align*}
I_3(k)& \le
 \int_{(D(R,3 x_d/2)\setminus D(R, x_d/2))\setminus K_1}
 |x-y|^{-d-\alpha+k} \Phi\big(\frac{|x-y|^2}{x_d^2}\big) dy\\
 & + \int_{K_1\setminus K_{1/\sqrt{d}}} 
|x-y|^{-d-\alpha+k} \Phi\big(\frac{|x-y|^2}{x_d^2}\big) dy =:I_{31}(k)+I_{32}(k).
\end{align*}
For  $y_d \in (D(R,3 x_d/2)\setminus D(R, x_d/2))\setminus K_1$, 
we have $y_d\asymp x_d$ and $x_d\le 2 |x-y|$. Thus, using the change of 
variables $y_d=rt+x_d$ in the second line below, we get
\begin{align*}
I_{31}(k)
&= 
c \int_{x_d/2}^R r^{d-2}\int_{x_d/2}^{3x_d/2}   {(|x_d-y_d|+r)^{-d-\alpha+k}}\Phi\big(\frac{(|x_d-y_d|+r)^2}{x_d^2}\big) 
 dy_d\, dr\\
&=  
c\int_{x_d/2}^R r^{-\alpha+(k-1)}\int_{-\frac{x_d}{2r}}^{\frac{x_d}{2r}}(|t|+1)^{-d-\alpha+k}
\Phi\big(\frac{(|t|+1)^2r^2}{x_d^2}\big) 
dt\, dr,
\end{align*}
which is, by the change of variables $r=x_d s$, comparable to
\begin{align}\label{e:lnew22}
x_d^{k-\alpha}  
\int_{1/2}^{R/x_d}  s^{-\alpha+(k-1)}\int_0^{1/s} 
\Phi\big((t+1)^2s^2\big) 
{(t+1)^{-d-\alpha+k}}\,  dt \,ds.
\end{align}
Since 
$
\int_0^{1/s} 
\Phi\big((t+1)^2s^2\big) 
{(t+1)^{-d-\alpha+k}} dt\asymp \Phi(s^2)/s$ for 
$s >1/2$, from \eqref{e:lnew22} we get 
\begin{align}
\label{e:I_31(k)}
I_{31}(k)
&\le c 
x_d^{-\alpha+k} \Big(1+\int_{1}^{R/x_d} \frac{  \Phi(s^2)}{s^{2+\alpha-k}} ds  \Big).
\end{align}
For $y\in K_1\setminus K_{1/\sqrt{d}}$ 
it holds that $x_d/(2\sqrt{d})\le |y-x|\le \sqrt{d}x_d$. 
By using that the volume $|K_1\setminus K_{1/\sqrt{d}}|\asymp x_d^d$, we get
$$
I_{32}(k)\asymp 
 \int_{K_1\setminus K_{1/\sqrt{d}}}
x_d^{-d-\alpha+k}\Phi(1)\ dy \asymp x_d^{-\alpha+k}.
$$
Together with 
\eqref{e:I_31(k)} 
this gives
$$
I_{3}(k)
\le c 
x_d^{-\alpha+k} \Big(2+\int_{1}^{R/x_d} \frac{  \Phi(s^2)}{s^{2+\alpha-k}} ds  \Big).
$$ 
Therefore, combining this inequality, \eqref{e:l1bd} and the fact that
$\int_{1}^{R/x_d} \frac{  \Phi(s^2)}{s^{2+\alpha-k}} ds \ge  \int_{1}^{3/2} \frac{  \Phi(s^2)}{s^{2+\alpha-k}} ds \ge c>0$,
we conclude that 
\begin{align}\label{e:e:13bd}
I_1(k)+I_3(k)\le c 
x_d^{-\alpha+k} \int_{1}^{R/x_d} \frac{  \Phi(s^2)}{s^{2+\alpha-k}} ds.
\end{align}

\noindent 
(a-ii)
Clearly,  $y_d>x_d$ for  $y_d \in D(R, R)\setminus D(R, 3x_d/2)$.
Thus, using the change of variables $y_d=x_d h$ and $r=x_d s$ in the second line below, 
we get
\begin{align}
&I_2(k)\le c  \int_0^{2R} r^{d-2} \int_{(3x_d/2)}^{R}
\Phi\big(\frac{((y_d-x_d)+r)^2}{x_dy_d}\big)\frac{dy_d\, dr }
{((y_d-x_d)+r)^{d+\alpha-k}} 
\nn\\
&= x_d^{-\alpha+k}   \int_{3/2}^{R/x_d} 
\int_0^{2R/x_d} 
\frac{s^{d-2}}{[(h-1)+s]^{d+\alpha-k}} 
\Phi\big(\frac{[(h-1)+s]^2}{h}\big)ds\, dh,
\label{e:RS1}
\end{align}
which is, by the change of variables $s=(h-1)t$ and using $(h-1)/h\asymp 1$, equal to 
\begin{align}
&x_d^{-\alpha+k} \int_{3/2}^{R/x_d} 
\int_0^{\frac{2R}{(h-1)x_d}}
\frac{t^{d-2}}{(h-1)^{1+\alpha-k}
(1+t)^{d+\alpha-k}} \Phi\big(\frac{(h-1)^2}{h} (1+t)^2\big)
dt\, dh\nonumber\\
&\asymp x_d^{-\alpha+k} \int_{3/2}^{R/x_d}
h^{-1-\alpha+k}
  \int_0^{\frac{2R}{(h-1)x_d}}
\frac{ \Phi\big(h(1+t)^2\big)}{
(1+t)^{d+\alpha-k}}
t^{d-2}dt\, dh. \label{e:RS2}
\end{align}
Then, for  $3/2\le h \le R/x_d$, we have
$\frac{2R}{(h-1)x_d} \ge \frac{2R}{R-x_d} \ge 2$. 
In particular, 
$$
 \int_{1}^{\frac{2R}{(h-1)x_d}}
\frac{\Phi(ht^2)dt}
{t^{2+\alpha-k}}  \ge  \int_{1}^{2}\frac{\Phi(ht^2)dt}
{t^{2+\alpha-k}} \asymp  \Phi(h), \quad 
3/2\le h \le R/x_d,
$$
Thus, 
\begin{align}
& 
 \int_0^{\frac{2R}{(h-1)x_d}}
\frac{ \Phi\big(h(1+t)^2\big)}{
(1+t)^{d+\alpha-k}}
t^{d-2}dt\asymp\Phi(h)+
\int_{1}^{\frac{2R}{(h-1)x_d}}
\frac{\Phi(ht^2)dt}
{t^{2+\alpha-k}} \asymp 
\int_{1}^{\frac{2R}{(h-1)x_d}}
\frac{\Phi(ht^2)dt}
{t^{2+\alpha-k}}. \label{e:RS31}
\end{align}
Combining \eqref{e:RS1}--\eqref{e:RS31}, we get
\begin{align}\label{e:RS32}
I_2(k) &
\le cx_d^{-\alpha+k}  \int_{3/2}^{R/x_d}\int_{1}^{\frac{2R}{(h-1)x_d}}
\frac{\Phi(ht^2)dt}
{t^{2+\alpha-k}} \frac{dh}{h^{1+\alpha-k}}.
\end{align}
Since $x_d \le R/2$, we have
\begin{align}\label{e:RS33}
 \int_{3/2}^{R/x_d}
\int_{1}^{\frac{R}{(h-1)x_d}}\frac{\Phi(ht^2)dt}
{t^{2+\alpha-k}} \frac{dh}{h^{1+\alpha-k}} \ge 
\int_{3/2}^{2}
\int_{1}^{\frac{R}{x_d}}\frac{\Phi(ht^2)dt}
{t^{2+\alpha-k}} \frac{dh}{h^{1+\alpha-k}} \asymp \int_{1}^{\frac{R}{x_d}}\frac{\Phi(t^2)dt}
{t^{2+\alpha-k}} .
\end{align}
Combining \eqref{e:e:13bd},\eqref{e:RS32} and \eqref{e:RS33}, we conclude that 
\begin{align}\label{e:RS34}
&I(k) \le c  x_d^{-\alpha+k}  \int_{3/2}^{R/x_d}
 \int_{1}^{\frac{2R}{(h-1)x_d}}\frac{\Phi(ht^2)dt}
{t^{2+\alpha-k}} \frac{dh}{h^{1+\alpha-k}}\nn\\
&\le x_d^{-\alpha+k} 
 \int_{1}^{6R/x_d}\int_{1}^{\frac{6R}{x_dh}}
\frac{\Phi(ht^2)dt}
{t^{2+\alpha-k}} \frac{dh}{h^{1+\alpha-k}}=: 
 x_d^{-\alpha+k}II(k, 6R/x_d)
\end{align}
where, by Fubini's theorem (for $a \ge 4$),
\begin{align}
&II(k, a):=\int_{1}^{a}
\int_{1}^{a/h}\frac{\Phi(ht^2)dt}
{t^{2+\alpha-k}} \frac{dh}{h^{1+\alpha-k}}=
\int_{1}^{a}
\int_{1}^{a/t}\frac{\Phi(ht^2)dh}{h^{1+\alpha-k}}
\frac{dt}{t^{2+\alpha-k}} \label{e:RS35}\\ &\le c 
\int_{1}^{a} \Phi(at)\left(\frac{t}{a}\right)^{\beta_1}
\int_{1}^{a/t}\frac{dh}{h^{-\beta_1+1+\alpha-k}}
\frac{dt}{t^{2+\alpha-k}}.
\label{e:RS351}
\end{align}
When $k>\alpha-\beta_1$, from \eqref{e:RS351} we have 
\begin{align}\label{e:RS36}
& II(k, a)\le c 
\int_{1}^{a} \Phi(at)\left(\frac{a}{t}\right)^{-\alpha+k}
\frac{dt}{t^{2+\alpha-k}}=ca^{k-\alpha}\Phi(a/2)\int_{1}^{a}\frac{ \Phi(at)}{\Phi(a/2)}
\frac{dt}{t^2}\nn\\
&\le ca^{k-\alpha}\Phi(a)\int_{1}^{a}\
\frac{dt}{t^{2-\beta_2}}  \le ca^{k-\alpha}\Phi(a)\int_{1}^{\infty}\
\frac{dt}{t^{2-\beta_2}} \asymp a^{k-\alpha}\Phi(a).
\end{align}
If $k=\alpha-\beta_1$, from \eqref{e:RS351} we have
\begin{align}\label{e:RS370}
& II(k, a) 
\le ca^{k-\alpha}\Phi(a/2)\int_{1}^{a}\frac{ \Phi(at)}{\Phi(a/2)}
\log({a/t})\frac{dt}{t^2}\le ca^{k-\alpha}\Phi(a)\int_{1}^{a} \log(a)
\frac{dt}{t^{2-\beta_2}} \nn\\
& \le c\log(a)a^{k-\alpha}\Phi(a)\int_{1}^{\infty}\
\frac{dt}{t^{2-\beta_2}} \asymp \log(a) a^{k-\alpha}\Phi(a).
\end{align}
If $k<\alpha-\beta_1$, using $2+\alpha-k-\beta_1-\beta_2 >1+\alpha-k-\beta_1>1$, from \eqref{e:RS351} we have
\begin{align}\label{e:RS37}
& II(k, a) 
\le c
 a^{-\beta_1}\Phi(a/2)\int_{1}^{a}\frac{ \Phi(at)}{\Phi(a/2)}
\frac{dt}{t^{2+\alpha-k-\beta_1}}
\int_{1}^{\infty}\frac{dh}{h^{1+\alpha-k-\beta_1}}
\nn\\
&\le ca^{-\beta_1}\Phi(a/2)\int_{1}^{\infty}
\frac{dt}{t^{2+\alpha-k-\beta_1-\beta_2}} \asymp a^{-\beta_1}\Phi(a/2).
\end{align}
If $k>\alpha-\beta_2$, from \eqref{e:RS35} we have
\begin{align}\label{e:RS380}
& II(k, a) \le c
\int_{1}^{a}  h^{\beta_2} t^{2\beta_2}
\int_{1}^{a/t}\frac{dh}{h^{1+\alpha-k}}
\frac{dt}{t^{2+\alpha-k}}\asymp
\int_{1}^{a}  
(a/t)^{k+\beta_2-\alpha}
\frac{dt}{t^{2-k+\alpha-2\beta_2}}\nn\\
&\asymp a^{k+\beta_2-\alpha}
\int_{1}^{a} 
\frac{dt}{t^{2-\beta_2}}\le a^{k+\beta_2-\alpha}
\int_{1}^{\infty} 
\frac{dt}{t^{2-\beta_2}} \asymp a^{k+\beta_2-\alpha}.
\end{align}
If $k=\alpha-\beta_2$, from \eqref{e:RS35} we have
\begin{align}\label{e:RS38}
 II(k, a) &\le c
\int_{1}^{a}  \int_{1}^{a/t}t^{2\beta_2}h^{\beta_2}
\frac{dh}{h^{1-\beta_2}}
\frac{dt}{t^{1+\alpha}}=
c
\int_{1}^{a}  t^{2\beta_2}
\int_{1}^{a/t}\frac{dh}{h}
\frac{dt}{t^{1+\alpha}}
\nn\\&\le c  \log a
\int_{1}^{\infty}
\frac{dt}{t^{2-\beta_2}}
\le c \log a.
\end{align}
If $k<\alpha-\beta_2$, using $2+\alpha-k-2\beta_2 >1+\alpha-k-\beta_2>1$, from \eqref{e:RS35} we get
\begin{align}\label{e:RS39}
 II(k, a) &\le c
\int_{1}^{\infty}\frac{dt}{t^{2+\alpha-k-{2\beta_2}}}
\int_{1}^{\infty}\frac{dh}{h^{1+\alpha-k-\beta_2}}<\infty.
\end{align}
Therefore, combining \eqref{e:RS34}--\eqref{e:RS39},
we conclude that \eqref{e:ln} holds.

\smallskip
\noindent
(b) If we further assume that $k>\alpha$, then for $x_d \le R/2$,
\begin{align*}
&\int_{D(R, R)}\frac{\sB(x, y)}{|x-y|^{d+\alpha-k}}dy
 \le
\int_{\{|y-x|<x_d/2\}}\frac{dy}{|y-x|^{d+\alpha-k}}
 +
R^{k-\alpha} \Phi(R/x_d) \\
\le& c
x_d^{k-\alpha}+c R^{k-\alpha} \Phi(R/x_d)
\le c R^{k-\alpha}+c R^{k-\alpha} \Phi(R/x_d)
 \asymp R^{k-\alpha} 
  (1\vee \Phi(R/x_d)
 ).\end{align*}
\qed

\begin{lemma}\label{l:3}
For every $\alpha\in [1, 2)$,  there exists $C=C(\alpha)>0$ such that  for all $z\in U$,
\begin{align*}
\int_{D_{\wt z}(7, 7)}\frac{|\sB (y, z)-\sB (z, z)|}{|y-z|^{d+\alpha-1}}dy\le    
C\begin{cases}
\Phi(1/z_d)(1+{\bf 1}_{\{  1+\beta_1=\alpha  \}} |\log z_d|)
 &\text{ if } 1+\beta_1 \ge \alpha;\\
 [\Phi(1/z_d)  z_d^{1-\alpha+\beta_1}]\wedge z_d^{-\beta_2} 
 &\text{ if } 1+\beta_1 < \alpha<1+\beta_2;\\
  z_d^{1-\alpha} (1+{\bf 1}_{\{  1+\beta_2=\alpha  \}} |\log z_d|)
 &\text{ if } 1+\beta_2 \le \alpha .
\end{cases}
\end{align*}
\end{lemma}
\pf 
Since $\sB (z, z)\le c\sB(y,z)$ for all $y,z\in \R^d_+$, we have
\begin{align*}
\int_{D(7, 7)}\frac{|\sB (y, z)-\sB (z, z)|}{|y-z|^{d+\alpha-1}}dy &\le
\int_{D(7, 7)\cap \{|y-z|<z_d/2\}}\frac{|\sB (y, z)-\sB (z, z)|}{|y-z|^{d+\alpha-1}}dy\\
&  +c\int_{D(7, 7)\cap \{|y-z|\ge z_d/2\}}
\frac{\sB (y, z)}{|y-z|^{d+\alpha-1}}dy
=:I+II.
\end{align*}
If $y\in B(z, 2^{-1}z_d)$, then $|y-z|\le z_d/2\le y_d$ and $y_d\asymp z_d$, hence by \textbf{(A2)}, 
we have that 
\begin{align}\label{e:dfge1}
I\le cz_d^{-\theta}\int_{|y-z|<z_d/2}|y-z|^{\theta-d-\alpha+1}dy
=cz_d^{-\theta}\int_0^{z_d/2}r^{\theta-\alpha}dr \le cz_d^{1-\alpha}.
\end{align}
Since
\begin{align}\label{e:dfge01}
z_d^{1-\alpha} \le c 
\begin{cases}
z_d^{-\beta_1} \le c \Phi(1/z_d)
 &\text{ if } 1+\beta_1 \ge \alpha;\\
 [\Phi(1/z_d)  z_d^{1-\alpha+\beta_1}]\wedge 
 z_d^{-\beta_2} 
 &\text{ if } 1+\beta_1 < \alpha<1+\beta_2,
 \end{cases}
\end{align}
combining \eqref{e:dfge1} with  by Lemma \ref{l:1n}(a),
we get the lemma.
\qed

\subsection{Proof of Lemma \ref{l:key-lemma}}
Let $\psi$ be a non-negative  $C^{\gamma}$ function in $\R^d_+$ with bounded support and bounded derivatives such that
 $$
\psi(y)=\begin{cases}|\widetilde{y}|^\gamma, 
&y\in D(2^{-2}, 2^{-2});\\
                                1,&y\in D(2, 2)\setminus U;\\
                                 0, & y\in D(3, 3)^c,
\end{cases}
$$
where $\gamma \ge 2$ will be chosen later, and $\psi(y)\ge 4^{-\gamma}$ 
for $y\in U\setminus D(2^{-2}, 2^{-2})$. 
The function $\psi$ in 
$\R^d_+$ can be constructed
so that, for $y=(\widetilde{y}, y_d)$ with
$y_d\in (0, \frac18)$, $\psi(y)$ depends on $\widetilde{y}$ only. 
We  extend  $\psi$ to be identically zero in 
$\R^d_-$.

Note that  \textbf{(A4)} implies that $x\mapsto \sB(x, x)$ is a 
constant on $\R^d_+$. Without loss
of generality, we assume that $\sB(x, x)\equiv 
1$ and for simplicity, in the remainder of this subsection, we neglect the constant $\sA(d, \alpha)$ in $j(x, y)$.

For $z \in U$ and $|y-z|>6$,  $|y| \ge |y-z| -|z| >5$. Thus by \eqref{e:operator-interpretation} (with $r=6$), 
for $\alpha\in [1, 2)$, 
we have that for $z \in U$,
\begin{align}\label{e:b}
&L^\sB_ \alpha\psi(z)=\int_{\R^d_+ \cap \{ |y-z|<6\}}\frac{\psi(y)-\psi(z)-\nabla \psi(z)\cdot (y-z)}{|y-z|^{d+\alpha}}\sB (y, z)dy\nn\\
&-\psi(z)\int_{\R^d_+ \cap \{ |y-z|>6\}}\frac{1}{|y-z|^{d+\alpha}}\sB (y, z)dy
\nn\\ 
&+\int_{\R^d_+ \cap \{ |y-z|<6\}}\frac{\nabla \psi(z)\cdot (y-z)}{|y-z|^{d+\alpha}}(\sB (y, z)-1)dy-\int_{\R^d_- \cap \{ |y-z|<6\}}\frac{\nabla \psi(z) \cdot(y-z)}{|y-z|^{d+\alpha}}dy\nn\\
& \le c_1 \int_{\R^d_+ \cap \{ |y-z|<6\}}\frac{\sB (y, z)}{|y-z|^{d+\alpha-2}}dy
+\int_{\R^d_+ \cap \{ |y-z|<6\}}\frac{|\nabla \psi(z) \cdot (y-z)|}{|y-z|^{d+\alpha}}|\sB (y, z)-1|dy
\nn\\
&\qquad \qquad +\int_{B(z, 6) \setminus B(z, z_d)}\frac{|\nabla \psi(z)\cdot(y-z)|}{|y-z|^{d+\alpha}}dy.
\end{align}

\noindent
{\bf (a)}
When $\alpha\in (0, 1)$, $L^\sB_ \alpha\psi(z)$ is not really a principal value integral and, 
since  $\alpha < 1$,
by Lemma \ref{l:1n}(b) with $k=1$,
\begin{align*}
&L^\sB_\alpha\psi(z)=
\int_{\R^d_+ } \frac{\psi(y)-\psi(z)}{|y-z|^{d+\alpha}}\sB(y, z)dy 
\le 
\int_{\R^d_+ \cap \{ |y-z|<6\}} \frac{\psi(y)-\psi(z)}{|y-z|^{d+\alpha}}\sB(y, z)dy \\
 &\le  c_2\int_{D_{\wt z}(7,7)}
\frac{\sB(y, z)}
{|y-z|^{d+\alpha-1}}dy \le c_3 \Phi(1/z_d).
\end{align*}
Thus, the conclusion of (a) follows for $\alpha\in (0, 1)$.
For the remainder of the proof of (a), we assume that $\alpha\in [1,2)$.

\noindent
{\bf (a1) $\alpha\in [1,2)$ and 
$z\in  D(2^{-2}, 2^{-2})$:}
Since $z\in D(2^{-2}, 2^{-2})$, we have that
 $\psi(z)=\psi(\wt{z})=|\wt{z}|^{\gamma}$ and 
$ |\nabla  \psi(\wt{z})\cdot(\wt{y}-\wt{z})| \le
c_4 |\wt{z}|^{\gamma-1} |\wt{y}-\wt{z}| $.
 We use \eqref{e:b}  and get 
\begin{align*}
L^\sB_\alpha\psi(z)\le& 
c_1 \int_{D_{\wt z}(7,7)}\frac{\sB (y, z)}{|y-z|^{d+\alpha-2}}dy
 + \int_{D_{\wt z}(7,7)}\frac{|\nabla \psi(\wt{z})\cdot(\wt{y}-\wt{z})|}{|y-z|^{d+\alpha}}|\sB(y,z)-1|\, dy \\
&+\int_{B(z, 6) \setminus B(z, z_d)}\frac{|\nabla  \psi(\wt{z})\cdot(\wt{y}-\wt{z})|}{|y-z|^{d+\alpha}}  dy=:I+II+III.
\end{align*}
Since  $2>\alpha$, by  Lemma \ref{l:1n}(b)  with $k=2$ we get that
$I \le c_5  \Phi(1/z_d).$ 
Estimating $II$  by  Lemma \ref{l:3} and using \eqref{e:dfge01}, we get that
\begin{align*}
&II+ III \nn\\
&\le
c_6|\wt{z}|^{\gamma-1}
\int_{D_{\wt z}(7,7)}\frac{|\sB(y, z)-\sB(z, z)|}{|y-z|^{d+\alpha-1}}dy+
c_6 |\wt{z}|^{\gamma-1} \int_{B(z, 6) \setminus B(z, z_d)}\frac{dy}{|y-z|^{d+\alpha-1}} \\
&\le c_7 |\wt{z}|^{\gamma-1}\begin{cases}
\Phi(1/z_d)(1+{\bf 1}_{\{  1+\beta_1=\alpha \text{ or } \alpha=1  \}} |\log z_d|)
 &\text{ if } 1+\beta_1 \ge \alpha; \\
[\Phi(1/z_d)  z_d^{1-\alpha+\beta_1}]\wedge z_d^{-\beta_2} 
 &\text{ if } 1+\beta_1 < \alpha<1+\beta_2;\\
  z_d^{1-\alpha} (1+{\bf 1}_{\{  1+\beta_2=\alpha  \}} |\log z_d|)
 &\text{ if } 1+\beta_2 \le \alpha.
\end{cases}
\end{align*}
Combining the estimates for $I$, $II$ and $III$, we get that 
\begin{align}
&L^\sB\psi(z)\le 
c_8\Phi(1/z_d)-C(\alpha, p, \sB)  |\wt{z}|^{\gamma}z_d^{-\alpha}  \nn\\
&\quad+ 
c_8|\wt{z}|^{\gamma-1}
\begin{cases}
\Phi(1/z_d) |\log z_d|
 &\text{ if } 1+\beta_1 \ge \alpha; \\
 [\Phi(1/z_d)  z_d^{1-\alpha+\beta_1}]\wedge z_d^{-\beta_2} 
 &\text{ if } 1+\beta_1 < \alpha<1+\beta_2;\\
  z_d^{1-\alpha}  |\log z_d|
 &\text{ if } 1+\beta_2 \le \alpha
\end{cases}
 \label{e:case1-1} \\
&= c_8 \Phi(1/z_d)\, -c_8|\widetilde{z}|^{\gamma-1}\nn\\
&\quad
 \times 
\begin{cases}
\Phi(1/z_d) \left( \frac{C(\alpha, p, \sB)}{{c_8}
 }\frac{|\wt{z}|}{\Phi(1/z_d)z_d^\alpha}-|\log z_d|\, \right) 
\qquad \qquad   \text{ if } 1+\beta_1 \ge  \alpha; \\
([\Phi(1/z_d)  z_d^{1-\alpha+\beta_1}]\wedge z_d^{-\beta_2} )\left(
\frac{C(\alpha, p, \sB)}{c_8}
\frac{ |\wt{z}| }{[\Phi(1/z_d)  z_d^{1+\beta_1}]\wedge z_d^{\alpha-\beta_2} 
}-1  \right)\\
\qquad \qquad\qquad \qquad\qquad \qquad\qquad \qquad\qquad \qquad \text{ if } 1+\beta_1 < \alpha<1+\beta_2;\\
{z_d^{1-\alpha}}\left(
 \frac{C(\alpha, p, \sB)}{c_8}
\frac{|\widetilde{z}|}{z_d}
- |\log z_d|\, \right)
\qquad \qquad\qquad \qquad   \text{ if } 1+\beta_2 \le \alpha.
\end{cases}
\label{e:case1-2}
\end{align}
We consider three cases separately:

\noindent
{\it Case $1+\beta_2 \le \alpha$:}
Let $\gamma=3$ and choose $\wt{\kappa} \in (0,1)$ so that 
$t^{-1/2}-
 \frac{c_8}{C(\alpha, p, \sB)}|\log t|>0$ 
for $t \in (0, \wt{\kappa}]$. 
When 
$|\wt{z}|\ge z_d^{1/2}$ 
and $z_d\le \wt{\kappa}$, it follows from \eqref{e:case1-2} and the choice of 
$\wt{\kappa}$ that 
$$
L^\sB \psi(z)\le  
c_8\Phi(1/z_d)\, -c_8\frac{|\widetilde{z}|^2}{z_d^{\alpha-1}}
\left(
 \frac{C(\alpha, p, \sB)}{c_8}{z_d^{-1/2}}
-|\log z_d|\, \right) 
 \le c_8  \Phi(1/z_d).
$$
In case 
when $|\wt{z}|\le z_d^{1/2}$ 
and $z_d\le \wt{\kappa}$, using the fact that
$2-\alpha >0$,
we estimate 
the last term in \eqref{e:case1-1} by
\begin{eqnarray*}
|\wt{z}|^2z_d^{1-\alpha}|\log z_d|\,  
\le  z_d^{2-\alpha}|\log z_d| 
\le c_9\le c_{10}\Phi(1/z_d).
\end{eqnarray*}
Thus, in this case we can disregard the middle term in \eqref{e:case1-1} and obtain again that 
$L^\sB\psi(z)\le c_{11}\Phi(1/z_d)$.
Finally, it follows from \eqref{e:case1-1} that for $z\in U$ with $z_d\ge \wt{\kappa}$ it holds that 
$L^\sB\psi(z)\le c_{12}\le c_{13}\Phi(1/z_d)$.
Therefore 
$L^\sB \psi(z)\le c_{14}\Phi(1/z_d) $ 
for all $z\in D(2^{-2}, 2^{-2})$.

\bigskip
\noindent
{\it Case $1+\beta_1 < \alpha<1+\beta_2$:}
Choose a $\gamma \ge 3$ such that $1 > (\gamma-1)/\gamma >\beta_2-\beta_1$.
Then, using  $1+\beta_1 < \alpha<1+\beta_2$, we have
$$
\gamma-1> \frac{\beta_2-\beta_1}{1+\beta_1-\beta_2} > \frac{\alpha-1-\beta_1}{1+\beta_1-\beta_2}>0.
$$
Let 
$$M:=\frac{(\gamma-1)(1+\beta_1-\beta_2)}{\alpha-1-\beta_1}>1.
$$
When 
$|\wt{z}|\ge 
 \frac{c_8}{C(\alpha, p, \sB)} 
([\Phi(1/z_d)  z_d^{1+\beta_1}]\wedge z_d^{\alpha-\beta_2} )^{1/M}$,
 it follows from \eqref{e:case1-2} 
 that 
\begin{align*}
&L^\sB \psi(z)\le 
c_8 \Phi(1/z_d)
\nn\\&\quad  -c_8|\wt{z}|^{\gamma-1}
([\Phi(1/z_d)  z_d^{1-\alpha+\beta_1}]\wedge z_d^{-\beta_2} )
\left(([\Phi(1/z_d)  z_d^{1+\beta_1}]\wedge z_d^{\alpha-\beta_2} )^{-(M-1)/M} 
-1 \right)\\
& \le c_8 \Phi(1/z_d)-c_8|\wt{z}|^{\gamma-1}
([\Phi(1/z_d)  z_d^{1-\alpha+\beta_1}]\wedge z_d^{-\beta_2} )
 \left(z_d^{-(\alpha-\beta_2)(M-1)/M} -1 \right)
\nn\\
&\le  c_8 \Phi(1/z_d).
\end{align*}
In case  when 
$|\wt{z}|\le \frac{c_8}{C(\alpha, p, \sB)} ([\Phi(1/z_d)  z_d^{1+\beta_1}]\wedge z_d^{\alpha-\beta_2} )^{1/M}$, 
using
$\Phi(1/z_d) \le c_{15}z_d^{-\beta_2}$,
we estimate 
the last term in \eqref{e:case1-1} by
\begin{align*}
&c_8|\wt{z}|^{\gamma-1} \Phi(1/z_d)  z_d^{1-\alpha+\beta_1} 
 \le   c_{16} (\Phi(1/z_d)z_d^{1+\beta_1})^{\frac{\gamma-1}{M}}
 z_d^{1-\alpha+\beta_1}  \Phi(1/z_d) \\
& =  c_{16} [\Phi(1/z_d)z_d^{1+\beta_1}
 z_d^{-1-\beta_1+\beta_2} ]^{\frac{\alpha-1-\beta_1}{1+\beta_1-\beta_2}} \Phi(1/z_d) 
\nn\\&=  c_{16} [\Phi(1/z_d)z_d^{\beta_2}
 ]^{\frac{\alpha-1-\beta_1}{1+\beta_1-\beta_2}} \Phi(1/z_d) 
\le  c_{17}\Phi(1/z_d).
\end{align*}
Thus, in this case we can disregard the middle term in \eqref{e:case1-1} and obtain 
 $L^\sB\psi(z)\le c_{18}\Phi(1/z_d)$.

\noindent
{\it Case $1+\beta_1 \ge  \alpha$:}
Let $\gamma=2$ and choose $\wt{\kappa} \in (0,1)$ so that 
$t^{-(\alpha-\beta_2)/2}-|\log t|>0$ 
for $t \in (0, \wt{\kappa}]$. 
When 
$|\wt{z}|\ge \frac{c_8c_{15}}{C(\alpha, p, \sB)}z_d^{(\alpha-\beta_2)/2}$ and $z_d\le \wt{\kappa}$,
it follows from \eqref{e:case1-2}  and $\Phi(1/z_d) \le c_{15}z_d^{-\beta_2}$
 and the choice of 
$\wt{\kappa}$ that 
\begin{align*}
L^\sB \psi(z)&\le 
c_8 \Phi(1/z_d)\,-c_8
|\widetilde{z}|\Phi(1/z_d) \left( 
 \frac{C(\alpha, p, \sB)}{c_8c_{15}}
\frac{|\wt{z}|}{z_d^{\alpha-\beta_2}}-|\log z_d|\, 
\right) 
\\
&\le 
c_8 \Phi(1/z_d)\, -c_8|\widetilde{z}|
\Phi(1/z_d)
 \left(z_d^{-(\alpha-\beta_2)/2}- |\log z_d|\, \right) 
\le c_{19}  \Phi(1/z_d).
\end{align*}
In case  when 
$|\wt{z}|\le  \frac{c_8c_{15}}{C(\alpha, p, \sB)}z_d^{(\alpha-\beta_2)/2}$, 
we estimate 
the last term in \eqref{e:case1-1} by
\begin{eqnarray*}
|\wt{z}|\Phi(1/z_d) |\log z_d|\,  
\le  c_{20}[z_d^{(\alpha-\beta_2)/2}
 |\log z_d|]\Phi(1/z_d)
 \le  c_{21}\Phi(1/z_d).
\end{eqnarray*}
Thus, in this case we can disregard the middle term in \eqref{e:case1-1} and obtain 
 $L^\sB\psi(z)\le c_{22}\Phi(1/z_d)$.
Finally, it follows from \eqref{e:case1-1} that for $z\in U$ with $z_d\ge \wt{\kappa}$ it holds that 
$L^\sB\psi(z)\le c_{23}\le c_{24}\Phi(1/z_d)$.
Therefore  
$L^\sB \psi(z)\le c_{25}\Phi(1/z_d) $ 
for all $z\in D(2^{-2}, 2^{-2})$.

\smallskip

\noindent
{\bf (a2)  $\alpha\in [1,2)$ and $z\in U\setminus D(2^{-2}, 2^{-2})$:}
We show that there exist constants $c_{26}>0$ and
$\kappa\in (0,1/4]$  such that (i) for $z_d\le \kappa$ and $|\wt{z}|\in [1/4,1/2)$   it holds that
 $L^\sB\psi(z)\le 0$; (ii) For $z_d\in [\kappa, 1/2)$, it holds that 
 $L^\sB\psi(z)\le c_{26}$.

We use \eqref{e:b} to get
\begin{align*}
L^\sB_{\alpha}\psi(z)&\le 
c_{27}\int_{D_{\wt z}(7,7)}\frac{\sB(y,z)}{|y-z|^{d+\alpha-2}}\, dy+c_{27}\int_{D_{\wt z}(7,7)}\frac{|1-\sB(y,z)|}{|y-z|^{d+\alpha-1}}\, dy +c_{27}
\int_{B(z, 6) \setminus B(z, z_d)}
\frac{dy}{|y-z|^{d+\alpha-1}}\, .
\end{align*}
Combining this with Lemmas 
\ref{l:1n}(b) and 
\ref{l:3}, we get that
there exists a positive constant $c_{28}>0$
such that
$$
L^\sB_{\alpha}\psi(z)\le c_{28} z_d^{-[(\alpha-1) \vee \beta_2]}
 |\log z_d|
\, , \quad z\in U\, .
$$
Thus
there exists $c_{26}=c_{26}(\kappa)$ such that $L^\sB\psi(z)\le L^\sB_{\alpha}\psi(z)\le c_{26}$ for all $z\in U$ with $z_d\ge \kappa$.

Finally, we assume that $|\wt{z}|\in (4^{-1}, 1)$. By the assumption on $\psi$, we have that 
$\psi(z)=\psi(\wt{z}, z_d)\ge 
4^{-\gamma}$. 
Since $a:=1 \wedge (\alpha-\beta_2)>0$, we have 
$\lim_{z_d\to 0}
z_d^{a} |\log z_d|  =0$ so  we can choose $\kappa >0$ so that
$$
z_d^{a}  |\log z_d|-
 \frac{C(\alpha, p, \sB)4^{-\gamma}}{c_{28}}\le 0
$$
for all $z_d\in (0,\kappa)$. Then,
\begin{eqnarray*}
L^\sB\psi(z)&=&L^\sB_{\alpha}\psi(z)-C(\alpha,p,\sB)z_d^{-\alpha}\psi(z)\le L^\sB_{\alpha}\psi(z)-C(\alpha,p,\sB)
4^{-\gamma}z_d^{-\alpha}\\
&\le &
  c_{28} z_d^{-\alpha}\left(z_d^{a}  |\log z_d| - \frac{C(\alpha, p, \sB)4^{-\gamma}}{c_{28}}\right)\le 0
\end{eqnarray*}
for all $z\in U\setminus D(2^{-2},2^{-2})$ with $|\wt{z}|\in  (4^{-1}, 1)$ and $z_d\in (0, \kappa)$.

 \noindent
{\bf (b)}
Recall that 
$h_{p, 1}(x)=x^p_d{\bf 1}_{D(1, 1)}(x)$.
Define $\phi:=h_{p, 1}-\psi$.
The function $\phi$ is obviously non-positive on $U^c$, hence 
Lemma \ref{l:key-lemma} (b2) holds true. Moreover, since $\psi((\wt{0},x_d))=0$, we have that $\phi((\wt{0}, x_d))=x_d^p$, for $(\wt{0}, x_d)\in U$, which is  Lemma \ref{l:key-lemma} (b1). 
Furthermore  Lemma \ref{l:key-lemma} (b3) 
follows from  Lemma \ref{l:estimate-of-L-hat-B}  and Lemma \ref{l:key-lemma} (a).
In fact, for $z\in U$
\begin{align*}
L^\sB \phi(z)=
L^\sB h_{p, 1}(z)
-L^\sB \psi(z)\ge 
-c_{29}\Phi(1/z_d)-c_{30}\Phi(1/z_d)
=- c_{31}\Phi(1/z_d).
\end{align*}
\qed

\section{Carleson estimates} \label{s:CE}
In this section we deal with the Carleson estimate.
The proof is similar to that of \cite[Theorem 1.2]{KSV20} and we provide only the part which requires some modification. 

\begin{thm}[Carleson estimate]\label{t:carleson}
Suppose $p\in (0, \alpha-\wt{\beta}_2) \cap [(\alpha-1)_+, \alpha-\wt{\beta}_2)$.
There exists a constant $C>0$  
such that for any $w \in\partial \R^d_+$, $r>0$, 
and any non-negative  function $f$ in $\R^d_+$ that is 
harmonic in $\R^d_+ \cap B(w, r)$ with respect to $Y$ and
vanishes continuously on $ \partial \R^d_+ \cap B(w, r)$, we have
\begin{equation}\label{e:carleson}
f(x)\le C f(x^{(0)}) \qquad \hbox{for all }  x\in \R^d_+\cap B(w,r/2),
\end{equation}
where 
 $x^{(0)}\in \R^d_+\cap B(w,4r/5)$  with 
$x^{(0)}_d\ge r/4$.
\end{thm}
\pf
We first choose $\beta_2$ such that
\eqref{e:Phi-infty} holds and $p\in (0, \alpha-\beta_2) \cap [(\alpha-1)_+, \alpha-\beta_2)$.
 By Lemma \ref{l:scaling} (a) and (b),
 it suffices to deal with the case $r=1$ and $\widetilde{w}=\widetilde 0$. 
 Moreover, by Theorem \ref{t:uhp}, we can assume that $x^{(0)}=(\wt 0, 4/5)$.
 
If $\kappa=0$, then  Corollary \ref{c:lemma4-1} (b) states that there is $n_0\ge 2$ such that for every $x\in \R^d_+$ it holds that $\P_x(\tau_{B(x, n_0 x_d)}=\zeta)\ge 1/2$. If $\kappa=C(\alpha, p, \sB)>0$, then
\begin{align*}
& \P_x(\tau_{B(x, n_0 x_d)}=\zeta)\ge \P_x(\tau_{B(x,  x_d/2)}=\zeta)=\E_x\int_0^{\infty} \1_{B(x,x_d/2)}(Y_s)\frac{C(\alpha, p,\sB)ds}{(Y^d_s)^{\alpha}}\\
& \ge C(\alpha,p, \sB)(x_d/2)^{-\alpha}\E_x \tau_{B(x,  x_d/2)}
\ge 
c \,C(\alpha, p, \sB),
\end{align*}
where the last inequality follows from Proposition \ref{p:exit-time-estimate} (a). Therefore, there exists a strictly positive constant $\delta_*$ depending on $\kappa$ such that for the corresponding process $Y$ (recall that we suppress dependence on $\kappa$ in the notation) it holds that
\begin{equation}\label{e:lemma-for carleson}
\P_x(\tau_{B(x, n_0 x_d)}=\zeta)\ge \delta_*, \quad \text{for all }x\in \R^d_+.
\end{equation}

Let $f$ be a non-negative function on $\R^d_+$ which is harmonic in $\R^d_+\cap B(0, 1)$ and vanishes continuously on 
$\partial \R^d_+\cap B(0, 1)$. 
By Theorem \ref{t:uhp} (b), it suffices to prove (\ref{e:carleson}) for 
$x\in \R^d_+\cap B(0,  1/(48 n_0))$.

Choose $\gamma\in (0,1/4)$ such that $0<\gamma < (\alpha-\beta_2)/(d+\alpha-2\beta_2)$ and $\gamma<\log 12/(\log n_0+\log 12)$ (the latter condition is equivalent to $n_0^{\gamma/(1-\gamma)}<12$). 
Recall that
 $x^{(0)}=(\wt 0, 4/5) \in \R^d_+\cap B(0, 1)$ and fix it.
For any $x\in \R^d_+\cap B(0, 1/(24 n_0))$, define
$$
B_0(x)=B(x, n_0 x_d)\, ,
\quad B_1(x)=B(x,x_d^{\gamma})\quad
\text{and}\quad
B_3=B(x^{(0)}, 4/15).
$$
Since $x\in B(0, 1/(24 n_0))$, we have $x_d<1/(12 n_0)$. 
By the choice of $\gamma$, we have that $B_0(x) \subset B_1(x)$. (Indeed, $(n_0x_d)/(x_d^{\gamma})=n_0 x_d^{1-\gamma}< n_0/(12 n_0)^{1-\gamma}=(n_0^{\gamma/(1-\gamma)}/12)^{1-\gamma}<1$.)

By \eqref{e:lemma-for carleson},
$\P_x(\tau_{B_0(x)}=\zeta)\ge \delta_*$ for 
$x\in \R^d_+\cap B(0, 1/(24 n_0)).$
By Theorem \ref{t:uhp} (b) there exists $\chi>0$ such that
$
f(x)<x_d^{-\chi} f(x^{(0)})$ for 
$x\in \R^d_+\cap B(0, 1/(24 n_0)).$
Since $f$ is harmonic in $\R^d_+\cap B(0,1)$,  for every 
$x\in \R^d_+\cap B(0, 1/(24 n_0)$,
\begin{align*}
f(x)&=\E_x\big[f\big(Y(\tau_{B_0(x)})\big); Y(\tau_{B_0(x)})\in B_1(x)\big]\nn\\
&\quad + \E_x\big[f\big(Y(\tau_{B_0(x)})\big);Y(\tau_{B_0(x)})\notin B_1(x)\big]. 
\end{align*}
Using \eqref{e:estimate-J}, \eqref{e:levy-system} and  Proposition \ref{p:exit-time-estimate}  
with
$
B_2:=B(x^{(0)}, 2/15)$, 
by the same arguments as in 
steps 1--2
of the proof of \cite[Theorem 1.2]{KSV}, we have that 
\begin{equation}\label{e:carleson-2}
f(x^{(0)})\ge c_1  \int_{\R^d_+\setminus B_3}J(x^{(0)},y) f(y)\, dy\, 
\end{equation}
and
\begin{align}\label{e:carleson-4}
&\E_x\big[f\big(Y(\tau_{B_0(x)})\big); Y(\tau_{B_0(x)}) \notin B_1(x)\big]\nn
\\&
\le c_{2} x_d^{\alpha} \left(\int_{(\R^d_+\setminus B_1(x))\cap B_3^c}
+ \int_{(\R^d_+\setminus B_1(x))\cap B_3}\right)J(x,y)f(y)dy=:c_{2} x_d^{\alpha}(I_1+I_2).
\end{align}

Suppose now that $|y-x|\ge x_d^{\gamma}$ and  $x\in \R^d_+\cap B(0, 1/(24 n_0))$. Then
$$
|y-x^{(0)}|\le |y-x|+2\le |y-x|+2x_d^{-\gamma}|y-x|\le 3x_d^{-\gamma}|y-x|.
$$
Moreover, since 
$(1/x_d)^{1-2\gamma}\ge (24 n_0)^{1-2\gamma} >4$, 
we have that $x^{(0)}_d=4/5 \ge x_d^{1-2\gamma}$. Therefore,
\begin{align*}
&J(x,y)\asymp \frac1{|x-y|^{d+\alpha} }
 \Phi\bigg(\frac{|y-x|^2}{x_d y_d}\bigg) \le c_{3}
  \frac{|x-y|^{-d-\alpha+2\beta_2}}{|x^{(0)}-y|^{2\beta_2}} (x_d^{(0)}/x_d)^{\beta_2}\Phi\bigg(\frac{|y-x^{(0)}|^2}{x_d^{(0)} y_d}\bigg)\\
&\le c_{4}
 \frac{(x_d^{\gamma}|y-x^{(0)}|)^{-d-\alpha+2\beta_2}}{|x^{(0)}-y|^{2\beta_2}} x_d^{-\beta_2}
 \Phi\bigg(\frac{|y-x^{(0)}|^2}{x_d^{(0)} y_d}\bigg)\le c_{5}  x_d^{-\gamma(d+\alpha-2\beta_2)-\beta_2} J(x^{(0)}, y).
\end{align*}
Now, using this and \eqref{e:carleson-2}, we get
\begin{align}\label{e:c:7}
I_1\le c_{5} x_d^{-\gamma(d+\alpha-2\beta_2)-\beta_2}\int_{(\R^d_+\setminus B_1(x))\cap B_3^c }J(x^{(0)}, y) f(y) \, dy\le  \frac{c_{6}  f(x^{(0)})}{x_d^{\gamma(d+\alpha-2\beta_2)+\beta_2}}.
\end{align}

If $y\in B_3$, 
then $y_d\asymp 1$ and
$$2\ge |x^{(0)}|+|x|+|y-x^{(0)}| \ge |y-x|\ge |x^{(0)}|-|x|-|y-x^{(0)}|>\frac45-\frac{1}{48n_0}-\frac14>\frac14.$$
Thus, for $y\in B_3$, 
$$
J(x,y)\le 
 \frac{c_{7}}{|x-y|^{d+\alpha} }\frac{|x-y|^{2\beta_2}}{x_d^{\beta_2}y_d^{\beta_2}}
\le c_8 x_d^{-\beta_2}.
$$
Moreover, by Theorem \ref{t:uhp}, there exists  
$c_9>0$ such that $f(y)\le c_9 f(x^{(0)})$ 
for all $y\in B_3$. 
Therefore, 
\begin{align}\label{e:c:8}
I_2 &\le 
 c_{9} f(x^{(0)}) \int_{(D\setminus B_1(x))\cap B_3} J(x, y) \, dy\nn\\
 &\le  c_{10} f(x^{(0)})\int_{2\ge |y-x|>1/4}   x_d^{-\beta_2}\, dy
\le c_{11}  x_d^{-\beta_2} f(x^{(0)}).
\end{align}

Combining \eqref{e:carleson-4}, \eqref{e:c:7} and \eqref{e:c:8} and using $\alpha-\beta_2>\gamma (d+\alpha-2\beta_2)>0$ (by the choice of $\gamma$),
we obtain
\begin{align}\label{e:c:9}
&\E_x[f(Y(\tau_{B_0(x)})); Y(\tau_{B_0(x)})\notin B_1(x)]\le c_{12}
f(x^{(0)}) x_d^{\alpha-\beta_2-\gamma(d+\alpha-2\beta_2)}.
\end{align}
 We choose 
  $\eta>0$ so that
$
 \eta^{\alpha-\beta_2-\gamma(d+\alpha-2\beta_2)}\le c_{12}^{-1}.
$
Then  for $x\in  \R^d_+ \cap B(0, 1/(24 n_0))$ with 
$x_d < \eta $,
we have by \eqref{e:c:9},
\begin{align*}
&\E_x\left[f(Y(\tau_{B_0(x)}));\, Y(\tau_{B_0(x)})\notin
B_1(x)\right]\nn\\ &\le  c_{12}\, f(x^{(0)})
\left(\eta^{-\gamma(d+\alpha-2\beta_2)-\beta_2+\alpha} +\eta^{-\beta_2+\alpha}\right)
\le f(x^{(0)})\, .
\end{align*}
The rest of the proof is the same as that of 
\cite[Theorem 1.2]{KSV}. 
Therefore we omit it. 
\qed

\section{Interior Green function estimates}\label{s:pgfe}
The goal of this section is to establish interior two-sided estimates of the Green function. We will distinguish between two cases: $d>\alpha$ and $d=1\le \alpha$. 
\subsection{Interior estimate: case $d>\alpha$}
In this subsection we establish the following interior two-sided estimates of the Green function in case $d>\alpha$.

\begin{prop}\label{p:green-int-est}
Suppose   $d>\alpha$. For any 
$a>0$, there exists  
$C=C(a)\ge 1$ 
such that for all $x,y\in \R^d_+$ satisfying  
$|x-y|\le a(x_d\wedge y_d)$, it holds that
$$
C^{-1}|x-y|^{-d+\alpha}\le G(x, y)\le C|x-y|^{-d+\alpha}.
$$
\end{prop}
\smallskip

We will first prove the upper bound which is more difficult. The idea of obtaining the upper bound of 
the Green function using the Hardy inequality originated from \cite{KSV20}.
The proof will be given through a number of auxiliary results. 

For $b>0$, 
let $\R^d_{b+} :=\{x\in \R^d_+:\, x_d \ge b\}$.  Define   
$$
Q(u,u):=
\int_{\R^d_{1+} }\int_{ \R^d_{1+} } 
(u(x)-u(y))^2 j(x, y)\, dx\, dy
$$
and
$\mathcal D(Q)=\{u\in 
L^2(\R^d_{1+}): 
Q(u, u)<\infty\}$.
Then $(Q, \mathcal D(Q))$ 
is a regular Dirichlet form and the corresponding symmetric Hunt process $X^{(1)}=(X^{(1)}_t)_{t\ge 0}$ is the reflected stable process 
on $\R^d_{1+}$, see \cite{BBC}.
Let  $p^{(1)}(t,x,y)$  be the transition density of  $X^{(1)}$. 
Using the estimates of $p^{(1)}(t,x,y)$ (see \cite{CK03}) 
we  get that for every $\gamma\in (0, (d/\alpha-1)\wedge2)$,
\begin{align*}
h(x,y)&:=\int_0^{\infty}t^{\gamma}  p^{(1)}(t,x,y)\, dt  \asymp
\int_0^{\infty}t^{\gamma}    \left( \bk t^{-d/\alpha}\wedge\frac{t}{|x-y|^{d+\alpha}} \right) \, dt
 \nn\\&\asymp \frac1{|x-y|^{d-(\gamma+1)\alpha},} \quad \quad \quad \quad x,y\in 
\R^d_{1+}
\end{align*}
and
$$
\overline{h}(x,y):=\int_0^{\infty}t^{\gamma-1}  p^{(1)}(t,x,y)\, dt \asymp \frac1{|x-y|^{d-\gamma\alpha}},\quad x,y\in 
\R^d_{1+}.
$$
Thus, 
$$
q(x):= \frac{\overline{h}(x,\e_d)}{h(x,\e_d)}\asymp \frac1{|x-\e_d|^{\alpha}}.
$$
It now follows from the Hardy inequality in
\cite[Theorem 2 and Corollary 3]{BDK} that 
there exists 
$c_1>0$ such that
\begin{equation}\label{e:estimate-Q}
Q(u,u)\ge c_1 
\int_{\R^d_{1+}}
u(x)^2 \frac{dx}
{|x-\e_d|^{\alpha}} 
\quad \textrm{for all }u\in 
L^2(\R^d_{1+}).
\end{equation}
 Using \eqref{e:estimate-Q},
we now follow the argument leading  to \cite[Corollary 4.4]{KSV20} line by line 
to get the following result.
\begin{prop}\label{c:extended-DF}
Suppose   $d>\alpha$. There exists 
$C>0$
such that for every non-negative Borel   function $f$ satisfying 
$\int_{\R^d_+}f(x) Gf(x)\, dx <\infty$ and every $z_b=(\wt{0},b)$ 
with $b\ge 0$, it holds that
$$
\int_{\R^d_{1+}} 
\frac{|Gf(x+z_b)|^2}
{|x-\e_d|^\alpha} \, dx
 \le C \int_{\R^d_+}f(x) Gf(x)\, dx.
$$
\end{prop}

\begin{prop}\label{p:l2norm} 
Suppose   $d>\alpha$. There exists $C>0$ such that for any $x\in \R^d_+$ with $x_d >6$,
it holds that
$
\int_{B(x,4)}  (G1_{B(x,4)}(y))^2\, dy \le C.
$
\end{prop}
\pf Without loss of generality we assume that 
$x=(\wt{0}, x_d)$. 
Set $B=B(x,4)$ and let $u=G \1_B$. 
It follows from \eqref{e:estimate-Q}  that for any $v\in C^\infty_c(\R^d_+)$,
\begin{align*}
 \int_B|v(y)|dy&\le |B|^{1/2}\left(\int_Bv^2(y)dy\right)^{1/2}\le 
c(x_d)
\left(\int_{\R^d_{1+}}v^2(y)\frac{dy}
{|y-\e_d|^\alpha}\right)^{1/2}\le c(x_d) (Q(v, v))^{1/2}.
\end{align*}
Thus $1_B(y)dy$ is of finite 0-order energy integral and $u\in \FF_e$,
where ${\mathcal F}_e$ is the extended Dirichlet space. 
By the definition of $\FF_e$,  there exists a $Q$-Cauchy sequence
$\{u_n\}\subset \FF$ with $u_n\to u$ almost everywhere.
Thus by \eqref{e:estimate-Q}  and Fatou's lemma,
\begin{align*}
\int_Bu^2(y)dy&\le  
c(x_d)\liminf_{n\to\infty} \int_Bu_n^2(y)dy\le \liminf_{n\to\infty} \int_{\R^d_{1+}}u_n^2(y)\frac{dy}{|y-\e_d|^\alpha}\\
&\le c(x_d)\liminf_{n\to\infty}  Q(u_n, u_n) =c(x_d) Q(u, u)<\infty.
\end{align*}
\normal

Let $z=(\wt{0}, x_d-6)$ and $\wt{B}=B((\wt{0},6),4)\subset 
\R^d_{2+}$.
By using the change of variables $w=x-z$ and the fact that 
$|w-\e_d|\asymp 1$ 
 for $w\in\wt{B}$ in the first line, and then Proposition \ref{c:extended-DF} and the Cauchy inequality in the third line below, we have
\begin{align*}
 \| u\|_{L^2(B)}^2 =& \int_{\wt{B}}|u(w+z)|^2\, dw \le 
c_1 
\int_{\R^d_{1+}} 
|G\1_B(w+z)|^2 \frac{dw}
{|w-\e_d|^\alpha}\\
\le & \,c_2 \int_{\R^d_+}\1_B(y) G\1_B(y)\, dy \le c_2 |B|^{
1/2} \| u\|_{L^2(B)}.
\end{align*}
Since $\| u\|_{L^2(B)}<\infty$, we have that 
$ \| u\|_{L^2(B)} \le c_2 |B|^{
1/2}$. This completes the proof. \qed

\noindent
\textbf{Proof of Proposition \ref{p:green-int-est}. Upper bound.} 
By \eqref{e:scaling-of-G} and Theorem  \ref{t:uhp}, it suffices to deal with 
$x, y\in \R^d_+$ with $|x-y|=1$ and 
$x_d\wedge y_d>10$.

We fix now two points $x^{(0)}$ and $y^{(0)}$ in $\R^d_+$ such that $|x^{(0)}-y^{(0)}|=1$,
$x^{(0)}_d \wedge y^{(0)}_d \ge10$ 
 and $\wt{x^{(0)}}=\wt{0}$. Let $E=B(x^{(0)},1/4)$,
$F=B(y^{(0)},1/4)$ and $D=B(x^{(0)},4)$. Let 
$f=G\1_E$ and $u=G\1_D$. 
Since $z \mapsto G(y^{(0)}, z )$ is harmonic in $B(x^{(0)},1/2)$ with respect to $Y$
and $f$ is harmonic in  $B(y^{(0)},1/2)$ with respect to $Y$, 
by applying 
Theorem  \ref{t:uhp}
to $f$ and $z \mapsto G(y^{(0)}, z )$, we get 
$$
f(y^{(0)}) =\int_E G(y^{(0)}, z ) dz \ge c_1 |E|G( y^{(0)}, x^{(0)}), \quad 
  \int_F  f(y)^2dy \ge c_2  |F| f(y^{(0)})^2.
$$
Thus, using the symmetry of $G$ and Proposition \ref{p:l2norm}, we obtain
$$
G(x^{(0)}, y^{(0)})\le \frac{1}{c_1 |E|} f(y^{(0)})\le \frac{1}{c_1|E|}\left(\frac{1}{c_2|F|} \int_F  f(y)^2dy\right)^{1/2}  \le 
\frac{c_3}{|E|^{3/2}}\|u\|_{L^2(D)},
$$
for  $c_3=c_1^{-1}c_2^{-1/2}>0$.

We have shown that there is a $c_4>0$ such that 
$G(z,w)\le c_4$ for all $z,w\in \R^d_+$ with $|z-w|=1$ and $z_d\wedge w_d\ge 10$.
By Theorem  \ref{t:uhp}, 
there exists $c_5=c_5(a)>0$ such that 
$G(z,w)\le c_5$ for all $z,w\in \R^d_+$ with $|z-w|=1$ and 
$z_d\wedge w_d>a^{-1}$.

Now let $x,y\in \R^d_+$ satisfy $|x-y|\le  a(x_d\wedge y_d)$ and set
$$
x^{(0)}=\frac{x}{|x-y|}, \quad y^{(0)}=\frac{y}{|x-y|}.
$$
Then $|x^{(0)}-y^{(0)}|=1$ and 
$x^{(0)}_d\wedge y^{(0)}_d> a^{-1}$ 
so that $G(x^{(0)},y^{(0)})\le c_5$. 
By scaling in Lemma \ref{l:scaling}(c),
$$
\qquad \qquad 
G(x,y)=G(x^{(0)},y^{(0)})|x-y|^{\alpha-d}\le \frac{c_5}{|x-y|^{d-\alpha}}. 
\qquad \qquad \qquad \quad \,  \Box
$$

We continue now by providing a proof of the lower bound and 
will use 
a well-known capacity argument to show
that there exists $c>0$ such that $G(x,y)\ge c$ 
for all $x,y\in \R^d_+$ satisfying $|x-y|=1$ 
and $x_d\wedge y_d\ge 10$.
 For such $x$ and $y$, let $D=B(x,5)$, $V=B(x,3)$ and $W_y=B(y,1/2)$. 
Recall that, for any $W\subset \R^d_+$, $T_W=\inf\{t>0: Y_t\in W\}$.
By Lemma \ref{l:krylov-safonov} (with $\epsilon=1/2$ and $r=5/2$),
there exists a constant $c_1>0$ such that
\begin{align}
\label{e:TW}
\P_x(T_{W_y}<\tau_D)\ge c_1\frac{|W_y|}{|D|}=c_2 >0\, .
\end{align}
Recall that $Y^D$ is the process $Y$ killed upon exiting $D$ and   
denote by $G_D(\cdot, \cdot)$  the Green function of $Y^D$. 
Let $\mu$ be the capacitary measure  of $W_y$ with respect to $Y^D$ (i.e., with respect to the corresponding Dirichlet form). 
Then $\mu$ is concentrated on $\overline{W_y}$, $\mu(D)=\mathrm{Cap}^{Y^D}(W_y)$ and $\P_x(T_{W_y}<\tau_D)=G_D \mu(x)$.  By
\eqref{e:TW} and applying
Theorem  \ref{t:uhp}  to the function $G(x,\cdot)$,  we get
\begin{align}
&c_2\le \P_x(T_{W_y}<\tau_D)=G_D\mu(x)=\int_D G_D(x,z)\mu (dz)\le \int_D G(x,z)\mu (dz)
\nonumber \\ 
&\le c_3  G(x,y) \mu(D)=
c_3 G(x,y)\mathrm{Cap}^{Y^D}(W_y)\, .\label{e:1}
\end{align}

Recall that $X$ denotes the isotropic $\alpha$-stable process in $\R^d$ 
and that  $X^D$ 
is the part of the process $X$ in $D$.  By 
Lemma \ref{l:df-comparison}
 and monotonicity of $\mathrm{Cap}^{X^D}$,
$$
\mathrm{Cap}^{Y^D}(W_y)\le c_4 \mathrm{Cap}^{X^D}(W_y)\le  c_4 \mathrm{Cap}^{X^D}(V)\, .
$$
The last term, $\mathrm{Cap}^{X^D}(V)$, 
is just a number, say $c_5$,  depending only on the radii of $V$ and $D$. Hence, $\mathrm{Cap}^{Y^D}(W_y)\le c_4 c_5$. 
Inserting in \eqref{e:1}, we get that
$$
G(x,y)\ge c_2c_3^{-1}c_4^{-1}c_5^{-1}.
$$

\noindent
\textbf{Proof of Proposition \ref{p:green-int-est}. Lower bound.} We have shown above that there is a $c_6>0$ such that 
$G(z,w)\ge c_6$ for all $z,w\in \R^d_+$ with $|z-w|=1$ and $z_d\wedge w_d\ge 10$. 
The rest of the proof is completed by the same argument as 
that
used for the upper bound. \qed

\subsection{Interior estimates: case $d=1\le \alpha$}\label{ss:int1d}
In this subsection we establish the following interior two-sided estimates of the Green function in case 
$d=1\le \alpha$.

\begin{prop}\label{p:green-int-est-2}
Suppose $d=1\le \alpha$. For any 
$a>0$, there exists a constant 
$C=C(a)\ge 1$ 
 such that for all $x,y\in \R^d_+$ satisfying  $|x-y|\le a(x_d\wedge y_d)$, 
it holds that
\begin{align*}
C^{-1} \left(x\vee y \vee |x-y|\right)^{\alpha-1} \le  &\,G(x,y)\le  C \left(x\vee y \vee |x-y|\right)^{\alpha-1}, 
& \alpha>1;\nn\\
C^{-1} \log \left(e+\frac{x\vee y}{|x-y|}\right) \le  &\,G(x,y)\le C  \log \left(e+\frac{x\vee y}{|x-y|}\right), 
& \alpha=1.
\end{align*}
\end{prop}

Again, we prove this result through a number of lemmas. 
The first one deals with the isotropic stable process killed upon exiting an interval.
This result might be known. Since we could not pinpoint a reference, we give a full proof.

\begin{lemma}\label{l:Cap1}
Suppose $d=1 \le \alpha$. 
There exists $C>1$ such that
for any $x_0 \in \R$ and $r \in (0, 3/4)$,
$$
C^{-1}
(
1+{\bf 1}_{\alpha=1}\log^{-1}(1/r)) \le \mathrm{Cap}^{X^{B(x_0,1)}}(\overline{B(x_0,r)}) \le C
(
1+{\bf 1}_{\alpha=1}\log^{-1}(1/r)).
$$
\end{lemma}
\pf
Without loss of generality, we assume that $x_0=0$.
Recall that $G^{X}_{B(0,1)}(x,y)$ is the Green function of the isotropic $\alpha$-stable process $X$ killed upon 
exiting 
$B(0,1)$. 
It is known that,
see e.g. \cite[Corollary 3]{BB01}, 
\begin{equation}\label{e-green-ball}
G^{X}_{B(0,1)}(x,y)\asymp
\begin{cases}
 \log\left(1+\frac{
(1-|x|)(1-|y|)}{|x-y|^2}\right),
&
\alpha=1;\\
[(1-|x|)^{(\alpha-1)/2}(1-|y|)^{(\alpha-1)/2}]  \wedge \frac{(1-|x|)^{\alpha/2}(1-|y|)^{\alpha/2}}{|x-y|}, &
\alpha>1.
\end{cases}
\end{equation}
Let $\mathcal{P}$ denote the family of all probability measures on
$\overline{B(0,r)}$. It follows from \cite[p.159]{Fu} that
\begin{equation}\label{e-formula-for-capacity}
\mathrm{Cap}^{X^{B(0,1)}}(
\overline{B(0,r)})=\Big(\inf\limits_{\mu\in
\mathcal{P}}\sup\limits_{x\in \mathrm{supp}(\mu)}
G^{X}_{B(0,1)}\mu(x)\Big)^{-1}\, .
\end{equation}
Let $m_r$ be the normalized Lebesgue measure on $\overline{B(0,r)}$.
By \eqref{e-formula-for-capacity},
\begin{equation}\label{e-estimate-for-capacity}
\mathrm{Cap}^{X^{B(0,1)}}(
\overline{B(0,r)})\ge \Big(\sup\limits_{x \in
\overline{B(0, r)}} G^{X}_{B(0,1)}m_r(x)\Big)^{-1}\, .
\end{equation}
Further, using
\eqref{e-green-ball} in the second line below, we have that for 
$\alpha=1$,
\begin{align*}
&\sup\limits_{x\in \overline{B(0,r)}} 
G^{X}_{B(0,1)}m_r(x)
=\sup\limits_{x\in \overline{B(0,r)}} \int_{B(0,r)}
G^{X}_{B(0,1)}(x,y)\, m_r(dy)\nn\\
&\le c \sup\limits_{x\in \overline{B(0,r)}} \int_{B(0,r)}\log\left(1+|x-y|^{-2}\right)\, m_r(dy)
\\&\le \frac{c}{ r} \int_{B(x,2r)}\log\left(1+|x-y|^{-2}\right)\, dy
\le \frac{c}{ r}
\int_{B(0,2r)}\log\frac{2}{|y|}\, dy\, 
\le\,
c \log\frac1r\, ,
\end{align*}
for some constant $c>0$.
Similarly, 
for $\alpha>1$,
\begin{equation*}
\sup\limits_{x\in \overline{B(0,r)}} 
G^{X}_{B(0,1)}m_r(x)
\le c \sup\limits_{x\in \overline{B(0,r)}} \int_{B(0,r)} m_r(dy) =c.
\end{equation*}
 This together with
\eqref{e-estimate-for-capacity} yields the desired 
lower bound.

For the upper bound we use that for any probability measure $\mu$ on $\overline{B(0,r)}$ it holds that
$$
\mathrm{Cap}^{X^{B(0,1)}}(\overline{B(0,r)})\le \Big(\inf_{x\in \overline{B(0,r)}}
G^{X}_{B(0,1)}\mu(x)\Big)^{-1},
$$
see \cite[Lemma 5.54]{SV09} (and note that there is a typo in the display -- the inf on the left-hand side should be taken over $x\in K$). 
Take $\mu=\delta_r$. Then in case 
$\alpha > 1$, 
for $x\in \overline{B(0,r)}$,
\begin{align*}
 G^{X}_{B(0,1)}\mu(x)&=
G^{X}_{B(0,1)}(x,r)
 \ge c\left( (1-r)^{\alpha-1}\wedge \frac{(1-r)^{\alpha}}{2r}\right)\\
&= c (1-r)^{\alpha-1}\left(1\wedge \frac{1-r}{2r}\right)\ge  c (1/4)^{\alpha-1}\left(1\wedge \frac{1}{8r}\right)\ge c .
\end{align*}
Hence, $\mathrm{Cap}^{X^{B(0,1)}}(\overline{B(0,r)})\le c$. 

In case $\alpha=1$, 
for $x\in \overline{B(0,r)}$,
\begin{align*}
& G^{X}_{B(0,1)}\mu(x)=
G^{X}_{B(0,1)}(x,r)\asymp \log\left(1+\frac{(1-|x|)(1-r)}{(x-r)^2}\right)\\
& \ge \log\left(1+\frac{(1-r)^2}{(2r)^2}\right)\ge \log\left(1+\frac{1}{4^3 r^2}\right)\ge c\log(1/r).
\end{align*}
Hence, $\mathrm{Cap}^{X^{B(0,1)}}(\overline{B(0,r)})\le c/\log(1/r)$.
\qed

Let
$B_n=(1-2^{-1}(1+2^{-n}), 1+2^{-1}(1+2^{-n}))$,
$n=1,2$.

\begin{lemma}\label{l:Cap2}
Suppose $d=1\le \alpha$.
There exists $C \ge 1$ such that for every $x \in \overline{B_2}$,
$$
C^{-1}(1+{\bf 1}_{\alpha=1} \log \bk \left(1/|x-1| \right)) \le G_{{B_1}}(x,1)\,\le\, C(1+{\bf 1}_{\alpha=1} \log \left(1/|x-1| \right)). 
$$
\end{lemma}
\pf Fix $x \in \overline{B_2}$ and let $r :=2^{-1} |x-1|$. Since
$\overline{B(1, r )}$ is a compact subset of ${B_1}$, there exists a
capacitary measure $\mu_{r}$ for $\overline{B(1, r )}$ with
respect to $Y^{{B_1}}$ such that
$$
\mbox{\rm Cap}^{Y^{{B_1}}} (\overline{B(1, r )}) = \mu_{r}
(\overline{B(1, r )})
$$
and $G_{B_1}\mu_r(x)=\P_x(T^{Y^{B_1}}_{\overline{B(1,r)}}<\infty)=  \P_x(T_{\overline{B(1,r)}}<\tau_{B_1})$ for $x\in B_1$
(see, for example, \cite[Section VI.4]{BG68} for details). 
Then by Theorem \ref{t:uhp} and using  \eqref{e:comp-cap} we have
\begin{align}
\int_{\overline{B(1, r )}}   
G_{B_1}(x,y)  \mu_{r} (dy) &\asymp  
G_{B_1}(x, 1) \mbox{\rm Cap}^{Y^{{B_1}}} (\overline{B(1, r )})
\nn \\
&\asymp 
G_{B_1}(x, 1) \mbox{\rm Cap}^{X^{{B_1}}}(\overline{B(1, r )}).
\label{e:ddd}
\end{align}

Moreover, $c\le \P_x(T_{\overline{B(1,r)}}<\tau_{B_1})\le 1$, where the left-hand side inequality follows from Lemma \ref{l:krylov-safonov} (with $\epsilon=1/10$ and the $r$ there equal to $5/8$). Therefore,
\begin{equation}\label{e:ddddd-2}
c\le \int_{\overline{B(1, r )}}  G_{B_1}(x,y)  \mu_{r}(dy) \le 1\, .
\end{equation}

Combining \eqref{e:ddd}-\eqref{e:ddddd-2} and applying Lemma \ref{l:Cap1}, we conclude that 
$$
G_{{B_1}}(x,1)\asymp \frac{1}{\mbox{\rm Cap}^{X^{{B_1}}}
\big(\overline{B(1, r)} \big)}
\asymp
\begin{cases}
  \log \left(1/r \right)\asymp \log
\left(1/|x-1| \right)&\text{if } 
\alpha=1;\\
1 &\text{if } 
\alpha>1.
\end{cases}
$$ 
\qed

\begin{lemma}
Suppose $d=1\le \alpha$.
There exists $C>0$ such that for all $x,y\in (0, 4/7)$ satisfying $|x-y|\le \frac58 (x\wedge y)$, 
\begin{align}
\label{e:G1d}
G_{(0, 1)} (x,y)\,\ge\, 
C\,
\begin{cases}
\log \left(e+\frac{x\vee y}{|x-y|}\right) 
&\text{if } \alpha=1;\\
(x \vee y \vee |x-y|)^{\alpha-1} &\text{if } \alpha>1.
\end{cases}
\end{align}
\end{lemma}
\pf
Note that by Lemma \eqref{l:scaling}(a), if $x<4/7$,
$$
G_{(0, 1)} (x,y)
\ge G_{(\frac14 x, \frac74 x)} (x,y)= x^{\alpha-1}G_{B_1}( 1, y/x).
$$
Thus, 
by Lemma \ref{l:Cap2}, for $x,y\in (0, 4/7)$ with
$|x-y|\le \frac58 (x\wedge y)$ we have 
$$
G_{(0, 1)} (x,y)
\ge  x^{\alpha-1}G_{B_1}( 1, y/x) \ge c x^{\alpha-1}
(1+{\bf 1}_{\alpha=1}\log (x/|y-x| ) ),
$$
so  \eqref{e:G1d} follows from this and 
the fact that $x \asymp x\vee y \vee |x-y|$. \qed

\noindent
\textbf{Proof of Proposition \ref{p:green-int-est-2}.}
Note that, if $|x-y|\le a(x\wedge y)$, then $x \asymp y$.
Without loss of generality, we assume $x \le y$.
We first consider the case 
$|x-y| \le \frac58 x $.
By \eqref{e:scaling-of-G},
$$G(x,y)= x^{\alpha-1}G(x/x, y/x)= x^{\alpha-1}G( 1, y/x).$$
Thus, 
 it suffices to show that  for $z \in \overline{B_2}$,
\begin{align}\label{e:Glog}
G( 1, z) \asymp 
1+{\bf 1}_{\alpha=1}\log (1/|z-1|).
\end{align}
By the strong Markov property, we have
\begin{equation}\label{e:sc2}
G(1, z)=
G_{{B_1}} (1, z) +\E_{1}\left[ G\big(X_{\tau_{{B_1}}},z\big) \right].
\end{equation}
Since $G(1, z) \ge G_{{B_1}} (1, z)$, the lower bound in \eqref{e:Glog}  follows from Lemma \ref{l:Cap2}.

For the upper bounds in  the proposition, define
$
h(v, w):=\E_v \left[ G\big (X_{\tau_{{B_1}}}, w\, \big) \right].
$
By Lemma \ref{l:Cap2}, to prove \eqref{e:Glog} we only need to show that 
\begin{equation}\label{e:sc211}
\sup_{z\in  \overline{B_2}}h(1, z) <\infty .
\end{equation}

For each fixed $v\in {B_1}$, the function $w \mapsto h(v, w)$ is harmonic in ${B_1}$ with respect to $Y$ and
for each fixed $w\in {B_1}$, 
$v\mapsto h(v, w)$ is  harmonic in
${B_1}$ with respect to the process $Y$. So it follows from
Theorem \ref{t:uhp} and the fact that $h(v, w) \le G(v, w)$ (see \eqref{e:sc2})\,
$$
\sup_{z \in  \overline{B_2}} h(1, z)
 \leq c \, \min_{v, w\in \overline{B_2}} h(v, w) \leq  c \,
\min_{v, w\in \overline{B_2}} G(v, w) 
 \le c\, G(1, 1/2)<\infty.
$$
We have shown that \eqref{e:sc211} and so  \eqref{e:Glog} hold. Thus, we have 
proved  the proposition  for   
$|x-y| \le \frac58 x $.
In particular, we have that 
$G(x,y) \asymp 1 $ for  
$\frac14 x <|x-y| \le \frac58 x $. 
Using  this and 
Theorem \ref{t:uhp}, we have 
$G(x,y) \asymp 1 $ for  $\frac14 x <|x-y| \le a x $.
The proof is complete.
\qed

\section{Preliminary upper bounds of Green function and Green potential}

The following result allows us to apply Theorem \ref{t:carleson} to get 
Proposition \ref{p:gfcnub} below, 
which is a key for obtaining the upper bound of  Green function.
In this section,  we always assume $p\in (0, \alpha-\wt{\beta}_2) \cap [(\alpha-1)_+, \alpha-\wt{\beta}_2)$.

\begin{thm}\label{t:green-function-decay}
For any $y\in \R^d_+$ and $w\in \partial \R^d_+$,  it holds that $\lim_{\R^d_+\ni x\to w}G(x,y)=0$.
\end{thm}
\pf By 
Lemma \ref{l:scaling}(b)
it suffices to show  $\lim_{|x|\to 0}G(x,y)=0$.  
We fix $y\in \R^d_+$ and consider $x\in \R^d_+$ with $|x|<2^{-10}y_d$. Let $B_1=B(y, y_d/2)$, $\overline B_1=\overline{B(y, y_d/2)}$
 and $B_2=B(y,y_d/4)$. For $z\in B_1$, we have that  $|z-x|\ge y_d/2-x_d\ge (7/16)y_d$. Thus, by the regular harmonicity of $G(\cdot, y)$  in $\R^d_+\setminus B(x, (7/16)y_d)$,
\begin{align}\label{e:green-function-decay}
G(x,y)&=\E_x[G(Y_{T_{\overline B_1}}, y), Y_{T_{\overline B_1}}\in B_1\setminus B_2]+\E_x[G(Y_{T_{\overline B_1}}, y), Y_{T_{\overline B_1}}\in B_2]\nn\\
&=:
I_1+I_2,
\end{align}
where, for any $V\subset \R^d_+$, $T_V:=\inf\{t>0: Y_t\in V\}$.
 For $z\in B_1$, $z_d>y_d/2$
 and so $|z-y|<y_d/2 \le z_d\wedge y_d$.
Thus, by Proposition \ref{p:green-int-est}, 
\begin{align}
\label{e:keyu}
G(z,y) \le c_1|z-y|^{-d+\alpha}, \quad z\in B_1.
\end{align}
Using \eqref{e:keyu}, we have
$$
I_1\le \sup_{z\in B_1\setminus B_2}G(z,y)\P_x(Y_{T_{\overline B_1}}\in B_1\setminus B_2)\le 
\frac{c_2}{y_d^{d-\alpha}}
\P_x(Y_{T_{\overline B_1}}\in B_1\setminus B_2).
$$
Further, it is easy to check that 
$J(w,z)\asymp J(w,y)$ for all $w\in \R^d_+\setminus B_1$ and $z\in B_2$. 
Moreover, by 
\eqref{e:keyu},
$$
\int_{B_2}G(y,z)\, dz \le c_1\int_{B(y,y_d/4)}|z-y|^{-d+\alpha} dz
=c_3 \int_0^{y_d/4} s^{\alpha-1} ds\le c_4 y_d^\alpha.
$$
Therefore, by \eqref{e:levy-system},
\begin{align*}
&I_2 =  \E_x \int_0^{T_{\overline B_1}} \int_{B_2} J (Y_t,z)G(z,y)\, dz\, dt 
\le c_5  \E_x \int_0^{T_{\overline B_1}} J(Y_t, y \normal)y_d^{\alpha} \, dt \\
&\le  c_6 y_d^{\alpha}\, \E_x \int_0^{T_{\overline B_1}} \left(\frac{1}{|B_2|}\int_{B_2} J (Y_t,z)\, dz \right)dt
= \frac{c_7}{y_d^{d-\alpha}}\P_x(Y_{T_{\overline B_1}}\in B_2).
\end{align*}
Inserting the estimates for $I_1$ and $I_2$ into \eqref{e:green-function-decay} and using 
Proposition \ref{p:allcomparable}
we get that
$$
G(x,y)\le 
 \frac{c_8}{y_d^{d-\alpha}}\P_x(Y_{T_{\overline B_1}}\in \R^d_+)\le \frac{c_8}{y_d^{d-\alpha}}\P_x(Y_{\tau_{U(y_d/4)}}\in \R^d_+)\le \frac{c_9}{y_d^{d-\alpha-p}}x_d^p,
$$
which implies the claim. \qed

Using Theorem \ref{t:green-function-decay}, 
we can combine Propositions \ref{p:green-int-est} and \ref{p:green-int-est-2}
with Theorem \ref{t:carleson} to get the following result.
\begin{prop}\label{p:gfcnub}
There exists a constant $C>0$ such that for all $x, y\in \R^d_+$,
\begin{align}
\label{e:upper}
G(x,y) \le C
\begin{cases}
|x-y|^{-d+\alpha} &\text{if } d>\alpha; \\
\log\left(e+ \frac{ x \vee y}{|x-y|}\right)
&\text{if } d=1=\alpha;\\
(x \vee y \vee |x-y|)^{\alpha-1}
&\text{if } d=1<\alpha.
\end{cases}
\end{align}
\end{prop}
\pf
When  $x_d \wedge y_d \ge |x-y|/8$, \eqref{e:upper} 
is proved in Propositions \ref{p:green-int-est} and \ref{p:green-int-est-2}.
In particular, for all $x,y$ such that $|x-y|=1$ and $2\ge x_d\wedge y_d\ge 1/8$, it holds that $G(x,y)\le c_1$ for some $c_1>0$.

By \eqref{e:scaling-of-G}, we only need to show that  for $x, y\in \R^d_+$ with $|x-y|=1$ and $x_d \wedge y_d \le 1/8$, 
\begin{align}
\label{e:upper09}
G(x,y) 
\le c_2
\begin{cases}
1 &\text{if } d>\alpha;\\
\log\left(e+ (x \vee y)\right)
&\text{if } d=1=\alpha;\\
(x \vee y \vee 1)^{\alpha-1}
&\text{if } d=1<\alpha
\end{cases}
\quad \asymp c_3.
\end{align}

Suppose that  $x, y\in \R^d_+$ with $|x-y|=1$, $x_d \le y_d$ and $x_d<1/8<y_d$.
Since $z \to G(z, y)$ is harmonic 
in $B((\wt x, 0), 1/4)$ with respect to $Y$ and vanishes on the boundary of $\R^d_+$ by Theorem \ref{t:green-function-decay}, we can use 
Theorem \ref{t:carleson} and see that there exists 
$c_4>0$ such that
\begin{align}
\label{e:cc1}
G(x,y)  \le c_4 G(x+(\wt 0, 1/8),y)  \le c_4 c_1.
\end{align}
Suppose that  $x, y\in \R^d_+$ with $|x-y|=1$, $x_d \le y_d$ and $y_d \le 1/8$ ($d \ge 2$). Then, 
since $z \to G(z, y)$ is harmonic 
in $B((\wt x, 0), 1/4)$ with respect to $Y$ and vanishes on the boundary of $\R^d_+$, 
by \eqref{e:cc1} and Theorem \ref{t:carleson}, we see that
$G(x,y)  \le c_4 G(x+(\wt 0, 1/8),y)  \le c_4^2 c_1$.
This finishes the proof of \eqref{e:upper09}.
\qed

\begin{lemma}\label{l:prelub} 
There exists $C>0$ such that for  $x,y \in \R^d_+$, 
$$
G(x, y) \le C
\left(\frac{x_d\wedge y_d}{|x-y|} \wedge 1\right)^p\times
\begin{cases} 
\frac{1}{|x-y|^{d-\alpha}} &\text{if } 
d>\alpha;\\
\log\big(e+ \frac{ x \vee  y}{|x-y|}\big)
&\text{if } d=1=\alpha;\\
(x \vee y \vee |x-y|)^{\alpha-1}
&\text{if }d=1<\alpha.
\end{cases}
$$
\end{lemma}

\pf 
We first choose $\beta_2$ such that
\eqref{e:Phi-infty} holds and $p\in (0, \alpha-\beta_2) \cap [(\alpha-1)_+, \alpha-\beta_2)$.
Suppose  $x, y\in \R^d_+$ satisfy $x_d \le 2^{-12}$ and $|x-y|=1$. Without loss of generality we assume that $\wt x =\wt 0$. 
Let $r=2^{-8}$.
For $z\in U(r)$ and $w\in \R^d_+\setminus D(r, r)$, we have $|w-z|\asymp |w|$. 
Moreover, by Proposition  \ref{p:gfcnub}, 
$G(w,y) \le c_1$ for $w \in  \R^d_+ \setminus B(y, r)$.
Thus, by using 
 Lemma \ref{l:new-two-estimates} with $q=0$ and 
 \eqref{e:B7},
\begin{align}
&\int_{\R^d_+\setminus D(r, r)}G(w,y)\frac{ \sB(z,w)}{|z-w|^{d+\alpha}}dw\nn\\ 
 \le &c_2
 \int_{\R^d_+ \cap B(y, r)}  G(w,y) \Phi\left(\frac{|w|^2}{z_d w_d}\right)\frac{dw}{ |w|^{d+\alpha}}
 +c_2
    \int_{\R^d_+\setminus (D(r, r) \cup B(y, r))} \Phi\left(\frac{|w|^2}{z_d w_d}\right) \frac{dw}{ |w|^{d+\alpha}} \nn\\\le     
&c_3\Phi\left(\frac{1}{z_d} \right)
\int_{\R^d_+ \cap B(y, r)} \frac{G(w,y)}{w_d^{\beta_2}|w|^{d+\alpha-2\beta_2}}dw+ 
c_3\Phi\left(\frac{1}{z_d} \right)=:c_3\Phi\left(\frac{1}{z_d} \right) (I+1).
\label{e:GG}
\end{align}
\noindent 
(i) We first estimate $I$ for $d > \alpha$:
Since $x\in U(r)$ and  $|y-x|=1$, we see that $|w| \asymp 1$ for $w\in B(y,r)$. If $r<y_d/2$, then $w_d\asymp y_d$ for $w\in B(y,r)$, and hence  by Proposition  \ref{p:gfcnub},
$$
 I\le c_4\int_{B(y,r)}G(w,y) dw \le c_5\int_{B(y,r)}\frac{dw}{|y-w|^{d-\alpha}}\le  c_6.
$$
If $r\ge y_d/2$, then 
$B(y, r)\cap \R^d_+\subset D_{\wt y}(3r, 3r)$
and thus by Proposition  \ref{p:gfcnub},
\begin{align*}
I&\le c_7\
\int_{B(y, y_d/2)}\frac{w_d^{-\beta_2}dw}{|y-w|^{d-\alpha}}
\nn\\
& +c_7
y_d^{\beta_2}
\int_{D_{\wt{y}}(3r, 3r) \cap \{|y-w| \ge y_d/2\}}
\left(\frac{|y-w|^{2}}{y_d w_d}\right)^{\beta_2} 
\frac{dw}{|y-w|^{d-\alpha+2\beta_2}} =:c_7(I_1+I_2).
\end{align*}
Clearly,
$$
I_1 \asymp y_d^{-\beta_2}\int_{B(y, y_d/2)}\frac{dw}{|y-w|^{d-\alpha}} 
 \le c_8 y_d^{\alpha-\beta_2} \le c_9.
$$
To estimate $I_2$ we use Lemma \ref{l:1n}(a) with $k=2\alpha-2\beta_2$
by taking $\Phi(t)=t^{\beta_2}$ and $\beta_1=\beta_2$ there.
Since $k+\beta_2 >\alpha$, by Lemma \ref{l:1n}(a), we get 
$I_2  \le c_{10} y_d^{\beta_2} y_d^{-\beta_2}  =c_{10}$.
Combining the estimates for $I_1$ and $I_2$, we get that $I\le c_{11}$.

\noindent 
(ii) We now estimate $I$ for  $d=1 \le \alpha$:
Since $1<y<1+2^{-9}$, 
we have  $w\asymp 1$ for $w \in (y-r, y+r)$. Thus,  
by  Proposition  \ref{p:gfcnub},
\begin{align*}
  I  &\le     \int_{y-r}^{y+r}   w^{\beta_2-1-\alpha}G(w,y) dw
\le    c _{12} \int_{y-r}^{y+r}    w^{\beta_2-1-\alpha}
   w^{\alpha-1} \log (e+ \frac{{\bf 1}_{\{\alpha=1\}}}{|w-y|})dw\nn\\
&=   c_{12}  \int_{y-r}^{y+r}    w^{-2+\beta_2}\log (e+ \frac{{\bf 1}_{\{\alpha=1\}}}{|w-y|})dw \le 
c_{13}  <\infty.
\end{align*}

\noindent 
(iii) By using \eqref{e:levy-system},
\eqref{e:GG} and the estimates for $I$ in (i) and (ii)  in the first inequality below, and Lemma \ref{l:upper-bound-for-integral} in the second, we get that 
\begin{align}\label{e:TAMSe6.41}
\E_x[G(Y_{\tau_{U(r)}}, y); Y_{\tau_{U(r)}}
\notin D(r, r)]
 \le c_{14}  \E_x\int^{\tau_{U(r)}}_0 
\Phi\left(\frac{r}{Y_t^d}\right) dt 
 \le c_{15}  x_d^p.
\end{align}

Let $x_0:=(\wt 0, r)$.
By Theorem \ref{t:carleson}, Propositions  \ref{p:gfcnub} and  \ref{p:allcomparable}, and scaling in Lemma \ref{l:scaling}(a),
we have
\begin{align}\label{e:TAMSe6.42}
\E_x[G(Y_{\tau_{U(r)}}, y); Y_{\tau_{U(r)}}
\in D(r, r)]\le 
c_{16} G(x_0, y)\P_x(Y_{\tau_{U(r)}}\in D(r, r))\le c_{17} x^p_d.
\end{align}
 Combining \eqref{e:TAMSe6.41} and \eqref{e:TAMSe6.42}, we get that 
 for $x, y\in \R^d_+$ satisfying 
 $x_d \le  2^{-12}$, $\wh x =\wt 0$ and $|x-y|=1$, 
\begin{align*}
G(x, y)&=
\E_x\left[G(Y_{\tau_{U(r)}}, y); Y_{\tau_{U(r)}}
\notin D(r, r)\right]\nn\\
&\quad + 
\E_x\left[G(Y_{\tau_{U(r)}}, y); Y_{\tau_{U(r)}}
\in D(r, r)\right]\le 
c_{18}x^p_d.
\end{align*}
Combining this with Proposition  \ref{p:gfcnub}, 
 \eqref{e:scaling-of-G} 
  and symmetry, we immediately get
the desired conclusion. \qed

As an application of Lemma \ref{l:prelub},
we get the following upper bound on Green potentials.

\begin{prop}\label{p:bound-for-integral-new}
{\rm (a)}  
Suppose $d> \alpha $.
There exists $C>0$ such that for any $\wt{w}\in \R^{d-1}$, $R>0$, any Borel set $D$ satisfying $D_{\wt{w}}(R/2,R/2) \subset D  \subset D_{\wt{w}}(R,R)$,
and any $x=(\wt{w}, x_d)$ with $0<x_d \le R/10$,
$$
\E_x \int_0^{\tau_D}(Y^{d}_t)^{ \gamma }\, dt 
= \int_D G_D(x,y)y_d^{ \gamma }\, dy  
\le C
\begin{cases} 
R^{\alpha+  \gamma   -p}x_d^p, & \gamma >p-\alpha;\\ 
x_d^p\log(R/x_d), &  \gamma =p-\alpha;\\ 
x_d^{\alpha+ \gamma }, &\hskip -.3in -p-1<  \gamma <p-\alpha.
\end{cases}
$$

\noindent 
{\rm (b)} Suppose  $d=1 \le \alpha$.
Let $\gamma>p-\alpha$. 	 
There exists $C>0$ such that for any 
$R>0$, any Borel set $D$ satisfying 
$(0, R/2) \subset D \subset (0, R)$, and any $0<x \le R/10$,
\begin{equation*}
\E_x \int_0^{\tau_D}(Y^{d}_t)^{\gamma }\, dt 
= \int_D G_D(x,y)y_d^{ \gamma }\, dy  \le 
CR^{\alpha+  \gamma   -p}x_d^p.
\end{equation*}
\end{prop}
\pf
(a) Using  Lemma \ref{l:prelub}, the proof is same as that of the upper bound of \cite[Proposition 6.10]{KSV20}.

\noindent 
(b)
Note that, by  scaling in Lemma \ref{l:scaling}(a), it suffices to show the lemma for $R=1$. We 
first consider the case 
$d=1<\alpha$.
By Lemma \ref{l:prelub},
\begin{align*}
G(x,y)\le c x^{\alpha-1} 
{\bf 1}_{0<y<x/2}+ 
x^{\alpha-1}c{\bf 1}_{x/2\le y<2x}
+ 
cx^{p}  y^{\alpha-p-1}{\bf 1}_{y\ge 2x}. 
\end{align*}
Thus, using 
$\gamma>p-\alpha$ and the fact that $D \subset (0,1)$, we have 
\begin{align*}
 &\int_D G_D(x,y)y^{ \gamma \normal}\, dy  \le   \int_0^{1} G(x,y)y^\gamma dy\nn\\
 & \le 
 cx^{\alpha-1}  \int_0^{x/2} y^{\gamma} dy+cx^{\gamma+\alpha-1} \int_{x/2}^{2x}dy+c x^{p}\int_{2x}^1   y^{\gamma+\alpha-1-p}dy\nn\\
  & \le 
 cx^{\alpha+ \gamma} +cx^{\alpha+ \gamma} +c x^{p}\int_{0}^1  y^{\gamma+\alpha-1-p}dy \le c x^{p}.
\end{align*}

We now consider the case $d=1=\alpha $.
By Lemma \ref{l:prelub},
\begin{align*}
G(x,y) \le c
{\bf 1}_{0<y<x/2}+ 
c \log\left(e+ \frac{ x }{|x-y|}\right) {\bf 1}_{x/2\le y<2x}
+ 
c (x/y)^p {\bf 1}_{y\ge 2x}
. 
\end{align*}
Thus, using 
$\gamma>p-1$ we have 
\begin{align*}
 &\int_D 
G_D(x,y)y^{ \gamma \normal}\, dy  \le   \int_0^{1} G(x,y)y^\gamma dy\nn\\
 & \le 
 c \int_0^{x/2} y^{\gamma} dy+cx^{\gamma} \int_{x/2}^{2x} \log\left(e+ \frac{ x }{|x-y|}\right) dy+c x^{p}\int_{2x}^1   y^{\gamma-p}dy\nn\\
   \qquad & \le 
 cx^{1+ \gamma} +cx^{1+ \gamma}\int_0^1\log\left(e+\frac{1}{t}\right)dt +c x^{p}\int_{0}^1  y^{\gamma-p}dy \le c x^{p}.
 \qquad \quad \  \Box
\end{align*}


\section{The proof of boundary Harnack principle and full Green function estimates} \label{s:fgfe}

\noindent
\textbf{Proof of Theorem \ref{t:BHP}.}
 By scaling in Lemma \ref{l:scaling}(a), it suffices to deal with the case $r=1$. Moreover, by
Theorem  \ref{t:uhp} (b), 
it suffices to prove \eqref{e:TAMSe1.8new}
for  $x, y\in D_{\wt w}(2^{-8}, 2^{-8})$.
Since $f$ is harmonic in $D_{\wt w}(2, 2)$ and vanishes continuously on $B(\wt w, 2)\cap \partial \R^d_+$,
it is regular harmonic in 
$D_{\wt w}(7/4, 7/4)$ and vanishes continuously on 
$B(\wt w, 7/4)\cap \partial \R^d_+$ (see \cite[Lemma 5.1]{KSV19} and its proof).
Throughout the remainder of this proof, we assume that 
$x\in D_{\wt w}(2^{-8}, 2^{-8})$. 
Without loss of generality we take $\widetilde{w}= \widetilde 0$.

Define $x_0=(\widetilde{x}, 
2^{-4})$.
By Theorem \ref{t:uhp}, 
Lemma \ref{l:POTAl7.4} and Proposition \ref{p:allcomparable}, 
we have for $x\in U$,
\begin{align}\label{e:TAMSe6.37}
f(x)&=
\E_x[f(Y_{\tau_{U}})]\ge \E_x[f(Y_{\tau_{U}}); Y_{\tau_{U}}
\in D(1/2, 1)\setminus D(1/2, 3/4)]
\nonumber\\
&
\ge c_1 f(x_0)\P_x(Y_{\tau_{U}}\in 
D(1/2, 1)\setminus D(1/2, 3/4))
\ge c_2 f(x_0)x_d^p.
\end{align}

Set $w_0=(\widetilde{0}, 2^{-7})$. 
Then, by \eqref{e:levy-system},
\begin{align}\label{e:POTAe7.27}
f(w_0)&\ge \E_{w_0}
\left[f(Y_{\tau_{U}}); Y_{\tau_{U}}\notin D(1, 1)
\right]
\ge \E_{w_0}\int^{\tau_{B({w_0}, 2^{-10})}}_0\int_{\R^d_+\setminus D(1, 1)}
J(Y_t, y)f(y)dydt\nonumber\\
&\ge c_{3}\E_{w_0}\tau_{B({w_0}, 2^{-10})}\int_{\R^d_+\setminus D(1, 1)}J({w_0}, y)f(y)dy=c_{4}\int_{\R^d_+\setminus D(1, 1)}J({w_0}, y)f(y)dy,
\end{align}
where in the last line we used Proposition \ref{p:exit-time-estimate} (a).

Note that $|z-y|\asymp |{w_0}-y|\asymp |y| \ge y_d \vee 1 \ge  z_d$ for any $z\in U$ and $y\in \R^d_+\setminus 
D(1, 1)$. Thus for $z\in U$ and $y\in \R^d_+\setminus 
D(1, 1)$,
\begin{equation}\label{e:new-estimate-for-J}
J(z,y)\le c_5 \frac{ 1}{|y|^{d+\alpha}}  
\Phi(   \frac{|y|^2}{z_dy_d}) 
\le c_6 z_d^{-\beta_2}
\Phi(   \frac{|y|^2}{y_d}) 
\frac{ 1}{|y|^{d+\alpha}}     \asymp
z_d^{-\beta_2} J({w_0},y).
\end{equation}

Combining \eqref{e:new-estimate-for-J} with  \eqref{e:POTAe7.27}
and using \eqref{e:levy-system} in the equality below and Proposition \ref{p:bound-for-integral-new} in the last inequality,  
 we now have
\begin{align}\label{e:POTAe7.29}
&\E_x\left[f(Y_{\tau_{U}}); Y_{\tau_{U}}\notin 
D(1, 1)\right]
=\E_x\int^{\tau_{U}}_0\int_{\R^d_+\setminus D(1, 1)}
J(Y_t, y)f(y)dydt\nonumber\\
&\le c_{7} \E_x\int_0^{\tau_{U}} 
(Y_t^{d})^{-\beta_2}dt \int_{\R^d_+\setminus 
D(1, 1)}J(w_0, y)f(y)dy\nn\\
&\le c_{8} f(w_0)\E_x\int_0^{\tau_{U}} 
(Y_t^{d})^{-\beta_2}dt \le c_9  f(w_0) x_d^p.
\end{align}

On the other hand, by Theorem \ref{t:uhp} and Theorem \ref{t:carleson} in the first inequality, 
and Proposition \ref{p:allcomparable} in the second, 
we have
\begin{align}\label{e:POTAe7.30}
\E_x\left[f(Y_{\tau_{U}}); Y_{\tau_{U}}\in D(1, 1)
\right]&\le c_{10}f(x_0)
\P_x\left(Y_{\tau_{U}}\in D(1, 1)
\right) \le c_{11} f(x_0) x_d^p.
\end{align}
Combining \eqref{e:POTAe7.29} and \eqref{e:POTAe7.30}, and using Theorem \ref{t:uhp}, we get
\begin{align*}
f(x)&=
\E_x\left[f(Y_{\tau_{U}}); Y_{\tau_{U}}\in D(1, 1)
\right]+
\E_x\left[f(Y_{\tau_{U}}); Y_{\tau_{U}}\notin D(1, 1)
\right]\nonumber\\
&\le c_{11}f(x_0)x_d^p+c_{9}f(w_0)x_d^p \le c_{12} f(x_0)x_d^p.
\end{align*}
 This with  
\eqref{e:TAMSe6.37} implies that 
$f(x) \asymp   f(x_0) x_d^p.
$
For any $y\in D(2^{-8}, 2^{-8})$, 
we have the same estimate with $f(y_0)$ instead of $f(x_0)$, 
where $y_0=(\widetilde{y}, 
2^{-4})$. 
By the Harnack inequality, we have $f(x_0)\asymp f(y_0)$. 
Thus,
$$
 \qquad \qquad \qquad \qquad \qquad \qquad
\frac{f(x)}{f(y)}\asymp \frac{x_d^{p}}{y_d^{p}}.
 \qquad \qquad \qquad \qquad \qquad \qquad \qquad \ \Box
$$

\begin{remark}\label{r:estimate0}
{\rm
Using \eqref{e:new-estimate-for-J}, 
one can follow the proofs \cite[Propositions 5.7 and 5.8]{KSV21} and show that any non-negative function which is regular harmonic near a portion of boundary vanishes continuously on that portion of boundary,
cf.~\cite[Remark 6.2]{BBC} and \cite[Lemma 3.2]{CK02}. Thus, the 
boundary Harnack principle also holds for regular harmonic functions.
We omit the details. 
}
\end{remark}

\noindent
\textbf{Proof of Theorem \ref{t:gfe}.}  
We first prove \eqref{e:FGE}.
Without loss of generality, we assume that $x_d \le y_d$ and $\wt x= \wt 0$.
By 
\eqref{e:scaling-of-G}, we can assume $|x-y|=1$ and just need to show that
\begin{align}
\label{e:Gfinal}
G(x,y)\asymp 
\begin{cases}
\left(x_d  \wedge 1 \right)^p\left(y_d \wedge 1 \right)^p 
&\text{if } d>\alpha\, ,\\
\left(x \wedge 1\right)^p\log\left(e+ y \right)
&\text{if } \alpha=1=d\, ,\\
\left(x \wedge 1\right)^p
( y \vee 1)^{\alpha-1}
&\text{if }\alpha>1=d\, .
\end{cases}
\end{align}
By \eqref{e:scaling-of-G}, Theorem \ref{t:uhp},  and Propositions 
\ref{p:green-int-est}
and \ref{p:green-int-est-2},
we only need to show \eqref{e:Gfinal} for $x_d \le 2^{-3}$ and $|x-y|=1$. In this case \eqref{e:Gfinal} reads
\begin{align}
\label{e:Gfinal2}
G(x,y)\asymp 
\begin{cases}
x_d^p  y_d ^p 
&\text{if } d \ge 2\, ,\\
x^p
&\text{if }d=1 \, .
\end{cases}
\end{align}
Thanks to 
Theorem \ref{t:green-function-decay},
\eqref{e:Gfinal2} is a direct consequence of Theorems \ref{t:BHP} and \ref{t:uhp},
see the proof of \cite[Theorem 1.2]{KSV21}.

From \eqref{e:FGE} it follows that
$\sup_{z\in \R^d_+ \setminus B(x, r)}G(x,z) <\infty$ for all $x \in \R^d_+$ and $r>0$.
The continuity of $y \mapsto G(x, y)$ on $\R^d_+\setminus \{x\}$ 
is a consequence of
 this observation and \cite[Proposition 6.3]{KSV22b}. 
 \qed

Using 
Propositions \ref{p:green-int-est}  and \ref{p:gfcnub} and Lemma  \ref{l:exit-probability-estimate}, 
the proof of the following lower bound is 
the same as that of \cite[Theorem 5.1]{KSV20}, hence we omit it.

\begin{thm}\label{t:GB}
Suppose $d> \alpha $, 
$p\in (0, \alpha-\wt{\beta}_2) \cap [(\alpha-1)_+, \alpha-\wt{\beta}_2)$.
For any $\eps \in (0, 1/4)$, there exists a constant $C>0$ such that for all 
$w \in \partial \R^d_+$,  $R>0$ and $x,y \in 
B(w, (1-\eps)R)\cap \R^d_+$, 
$$
G_{B(w, R)\cap \R^d_+}(x, y)\ge C
\left(\frac{x_d}{|x-y|}  \wedge 1 \right)^p\left(\frac{y_d}{|x-y|}  \wedge 1 \right)^p \frac{1}{|x-y|^{d-\alpha}}.
$$
\end{thm}

We now consider the lower bound in case $d=1 \le \alpha$.
\begin{thm}\label{t:GB1}
Let $d=1$. Suppose
$p\in (0, \alpha-\wt{\beta}_2) \cap [(\alpha-1)_+, \alpha-\wt{\beta}_2)$.
Then there exists a constant $c>0$ such that for all $R>0$ and all $x,y\in (0, R/2)$, 
$$
G_{(0,R)}(x,y)\ge c
\begin{cases}
\left(\frac{x\wedge y}{|x-y|} \wedge 1\right)^p\log\left(e+ \frac{ x \vee  y}{|x-y|}\right)
&\text{if } \alpha=1\, , \\
\left(\frac{x\wedge y}{|x-y|} \wedge 1\right)^p(x \vee y \vee |x-y|)^{\alpha-1}
&\text{if }\alpha>1\, .
\end{cases}
$$
\end{thm}
\pf By Lemma \ref{l:scaling}(a),
without loss of generality, we assume $R=1$ and $x \le y < 1/2$.
When $|x-y|\le \frac58 x$, the theorem follows from \eqref{e:G1d}.

Suppose $|x-y|=y-x> \frac58 x$. Then $y-x \le y=y-x+x <\frac{13}5(y-x)$.
Thus, we just need to show that 
$G_{(0,1)}(x,y)\ge c y^{\alpha-1} (x/y)^p$.

Since $y<1/2$, by Lemma \ref{l:scaling}(a), we have
$$
G_{(0, 1)} (x,y)\ge G_{(0, 2 y)} (x,y)= y^{\alpha-1}G_{(0, 2)}( 1, x/y).
$$
Since $x/y<8/(13)$, using Theorem \ref{t:BHP}, we get 
$$
G_{(0, 1)} (x,y)
\ge  y^{\alpha-1}G_{(0, 2)}( 1, x/y) \ge c G_{(0, 2)}( 1, 1/2)y^{\alpha-1} (x/y)^p
=c y^{\alpha-1} (x/y)^p\, .
$$
We have proved the theorem. 
\qed

As an application of Theorems  \ref{t:GB} and \ref{t:GB1}
we now get the full estimates of the following  Green potentials.

\begin{prop}\label{p:bound-for-integral-new-2} 
Suppose $d> \alpha $, 
$p\in (0, \alpha-\wt{\beta}_2) \cap [(\alpha-1)_+, \alpha-\wt{\beta}_2)$.
Then
for any 
$\wt{w} \in \R^{d-1}$, any Borel set $D$ satisfying $D_{\wt{w}}(R/2,R/2) \subset D  \subset D_{\wt{w}}(R,R)$ 
and any 
$x=(\wt{w}, x_d)$ with 
$0<x_d \le R/10$,
$$
\E_x \int_0^{\tau_D}(Y^{d}_t)^{ \gamma }dt 
= \int_D G_D(x,y)y_d^{\gamma }\, dy \asymp
\begin{cases} 
R^{\alpha+ \gamma   -p}x_d^p, &  \gamma >p-\alpha;\\ 
x_d^p\log(R/x_d), &   \gamma =p-\alpha, d> \alpha ;\\ 
x_d^{\alpha+\gamma }, & p-1<  \gamma <p-\alpha,  d> \alpha 
\end{cases}
$$
where the comparison constant  is
 independent of $\wt{w} \in \R^{d-1}$, $D$, $R$ and $x$.
\end{prop}
\pf
The upper bounds 
are given in Proposition \ref{p:bound-for-integral-new}.
Moreover, using  Theorem \ref{t:GB},
the proof for the lower bound for $d>\alpha$  is same as that 
of the lower bound of \cite[Proposition 6.10]{KSV20}.

Suppose $d=1 \le \alpha$ and  $\gamma>p-\alpha$.
By Lemma \ref{l:scaling}(a), it suffices to show the lemma for $R=1$. 
By Theorem \ref{t:GB1},
\begin{align*}
G_{(0, 1/2)}(x,y)y^\gamma \ge c
 x^{p}    y^{\gamma+\alpha-1-p} \quad \text{if }  1/4>y\ge 2x.
\end{align*}
Thus, using 
$\gamma>p-\alpha$ and the fact that $D \supset (0,1/2)$, we have that for $0<x<1/10$,
\begin{align*}
 &\int_D G_D(x,y)y^{ \gamma }\, dy  \ge   \int_{2x}^{1/4} 
G_{(0,1/2)} (x,y)y^\gamma dy\nn\\
\qquad \quad \  & \ge 
c x^{p}\int_{2x}^{1/4}   y^{\gamma+\alpha-1-p}dy \ge c x^{p}\int_{1/5}^{1/4}y^{\gamma+\alpha-1-p}dy = c x^{p}. 
\qquad \qquad \, \Box
\end{align*}

\section{Proof of Theorem \ref{t:estimate}}\label{s:proof-res-ker}
In this section we give a proof of Theorem \ref{t:estimate} for $d\ge 2$. The case $d=1$ has already been  treated in \cite{KSV22b}.
Let  $\wt e_1 =(1, 0,\dots,0)$ be the unit vector in the $x_1$ direction in $\R^{d-1}$.
For $a,b>0$, define 
\begin{align*}
h(a,b):=\int_{\R^{d-1}}\frac{\Psi\big({(|\wt u|+1+(1/a))^2 }b\big)}
{(|\wt u| + 1)^{d+\alpha}(|\wt u| +1+(1/a))^{d+\alpha}}
 d\wt u,
\end{align*}

\begin{align*}
\Upsilon(a,b,l):=\int_{\R^{d-1}}
\frac{\Psi\big((|\wt u-
(1/a) \wt e_1| + 1
)^2 b\big)}
{(|\wt u| + 1/l)^{d+\alpha}(|\wt u-
(1/a) \wt e_1| + 1
)^{d+\alpha}}
 d\wt u,\end{align*}
 and 
 $f(a,b):=\Upsilon(a,b,b)$,  $g(a,b):=\Upsilon(a,b,1)$.

\begin{lemma}\label{l:k}
For any $M>0$, 
$$
h(a,b) \asymp 
\begin{cases}
a^{d+\alpha} \Psi(b/a^2) & \text{if } a<M, b>0;\\
\Psi(b) & \text{if } a \ge M, b>0
\end{cases}
$$
with comparison constants depending on $M$.
\end{lemma}
\pf 
Using \eqref{e:Psi-infty0}, we have that for $a \ge M,$
\begin{align*}
h(a,b) &\asymp \int_{\R^{d-1}}\frac{\Psi\big({(|\wt u|+1)^2 }b\big)}
{(|\wt u| + 1)^{2d+2\alpha}}
 d\wt u \asymp
 {\Psi(b)}\int_{0}^{\infty} \frac{\Psi\big({(v+1)^2 }b\big)}{\Psi(b)(v+ 1)^{2d+2\alpha}}
 v^{d-2}
dv \\
&\le c  {\Psi(b)}\int_{0}^{\infty} \frac{ v^{d-2}}{(v+ 1)^{2d+2\alpha-2{\gamma_2}_+}}
dv
 \asymp
  \Psi(b).
\end{align*}
Similarly, using the lower bound in \eqref{e:Psi-infty0}, 
    \begin{align*}
h(a,b) & \ge c  {\Psi(b)}\int_{0}^{2} \frac{ v^{d-2}}{(v+ 1)^{2d+2\alpha+2{\gamma_1}_-
}}
dv\asymp    {\Psi(b)} \int_{0}^{2} 
 v^{d-2}
dv\asymp  \Psi(b).
\end{align*}

 For $a<M$,
 \begin{align*}
  h(a,b)& \ge 
  \int_{|\wt u| <1/a}\frac{\Psi\big({(|\wt u|+1+(1/a))^2 }b\big)}
{(|\wt u| + 1)^{d+\alpha}(|\wt u| +1+(1/a))^{d+\alpha}}
 d\wt u\asymp  a^{d+\alpha}
  \int_{|\wt u| <1/a}\frac{\Psi(b/a^2)}
{(|\wt u| + 1)^{d+\alpha}}
 d\wt u\\
&\ge  a^{d+\alpha}\Psi(b/a^2)
  \int_{|\wt u| <1/M}\frac{ d\wt u}
{(|\wt u| + 1)^{d+\alpha}}
 \asymp a^{d+\alpha}\Psi(b/a^2).
    \end{align*}

For the upper bound, we use \eqref{e:Psi-infty0}  
and get that for $a<1/M$,
    \begin{align*}
  h(a,b)& \asymp
  a^{d+\alpha}
  \int_{|\wt u| <1/a}\frac{\Psi(b/a^2)}
{(|\wt u| + 1)^{d+\alpha}}
 d\wt u
 +
 \int_{|\wt u|  \ge 1/a}\frac{\Psi\big({|\wt u|^2 }b\big)}
{|\wt u|^{2d+2\alpha}}
 d\wt u \nn\\
 & \le
a^{d+\alpha}\Psi(b/a^2)
  \int_{\R^{d-1}}\frac{ d\wt u}
{(|\wt u| + 1)^{d+\alpha}}
 +c\Psi(b/a^2)
 \int_{1/a}^{\infty}\frac{\Psi\big(bv^2\big)}
{\Psi(b/a^2)v^{d+2+2\alpha}}
 d v \nn\\
 & \le
  ca^{d+\alpha}\Psi(b/a^2)
 +c\Psi(b/a^2)
 \int_{1/a}^{\infty}\frac{ a^{2{\gamma_2}_+}   }
{v^{d+2+2\alpha-2{\gamma_2}_+}}
 d v 
\le c a^{d+\alpha}\Psi(b/a^2).  \end{align*}
\qed

  \begin{lemma}\label{l:g_2n} 
There exists a constant $C>0$ such that 
 $g (a,b)  \le C  \Psi(b) $ for all $a, b >0$.
 \end{lemma}
 \pf 
Since 
$d+\alpha-2{\gamma_2}_+>0$, 
$(|\wt u-
(1/a) \wt e_1| + 1)^{d+\alpha-2{\gamma_2}_+}\ge 1$. Thus,
  \begin{align*}
 g (a,b )
\le c  \Psi(b)\int_{\R^{d-1}}
\frac{(|\wt u| + 1)^{-d-\alpha}d\wt u}
{(|\wt u-
(1/a) \wt e_1| + 1
)^{d+\alpha-2{\gamma_2}_+}}
 \le c \Psi(b) \int_{\R^{d-1}}
\frac{d\wt u}{(|\wt u| + 1)^{d+\alpha}}  <c\Psi(b).
 \end{align*}
\qed

By the change of variables $1/a-u_1=v_1$ and $\wh u=\wh v$,
 and using that $|\wt v| =|(v_1, \wh{v})|=|(-v_1, \wh{v})|$, we see that for all $a,b,l>0$,
 \begin{align}\label{e:g111}
 &\int_{\R^{d-1}, u_1<\frac1{2a}}
\frac{\Psi\big((|\wt u-
(1/a) \wt e_1| + 1
)^2 b\big)}
{(|\wt u| + 1/l)^{d+\alpha}(|\wt u-
(1/a) \wt e_1| + 1
)^{d+\alpha}}
 d\wt u
\nn \\
= &\int_{\R^{d-1}, v_1>\frac1{2a}}
\frac{\Psi\big((|\wt v| + 1
)^2 b\big)}
{(|\wt v| +1)^{d+\alpha}(|\wt v-(1/a) \wt e_1|+1/l )^{d+\alpha}}
 d\wt v.
\end{align}
Thus, for all $a,b,l>0$,
 \begin{align}
\Upsilon(a,b,l)  &\le \int_{\R^{d-1}, u_1>\frac1{2a}}
\frac{\Psi\big((|\wt u-
(1/a) \wt e_1| + 1
)^2 b\big) }
{|\wt u|^{d+\alpha}(|\wt u-
(1/a) \wt e_1| + 1
)^{d+\alpha}}d\wt u\nn\\
&\qquad +
\int_{\R^{d-1}, \frac3{2a}\le u_1}
\frac{\Psi\big((|\wt u| + 1
)^2 b\big)}
{(|\wt u| +1)^{d+\alpha}|\wt u-(1/a) \wt e_1|^{d+\alpha}}
 d\wt u&\nn\\
 &\qquad+\int_{\R^{d-1}, \frac3{2a}>u_1>\frac1{2a}}
\frac{\Psi\big((|\wt u| + 1
)^2 b\big)d\wt u}
{(|\wt u| + 1
)^{d+\alpha}(|\wt u-(1/a) \wt e_1| +1/l)^{d+\alpha}} =:I+II+III(l).\label{e:f2n}
\end{align}
We use the above notations $I$, $II$ and $III(l)$ in the next two lemmas.
 \begin{lemma}\label{l:g_2}
 For any $M>0$, we have that 
 for all $a \in (0, M]$ and $b>0$, 
$
  g(a,b) \asymp a^{d+\alpha} \Psi(b/a^2),
$
with comparison constants depending on $M$.
\end{lemma}
\pf
For $\frac3{2a} \le |\wt  u|$,
we have 
 \begin{align}\label{e:g2}
\frac13|\wt u|=  |\wt u| -\frac23|\wt  u| \le  |\wt u-\frac1{a}\wt e_1|+\frac1{a}-\frac23 |\wt  u|   \le  |\wt u-\frac1{a}e_1|.
 \end{align}
 Thus, since $\Psi(t)t^{-(d+\alpha)/2}$ is almost decreasing because ${\gamma_2}_+-(d+\alpha)/2<0$, 
using the upper bound in \eqref{e:Psi-infty0}, we have that 
 for $\frac3{2a} \le |\wt  u|$,
  \begin{align}
  & \frac{\Psi\big((|\wt u-
(1/a) \wt e_1| + 1
)^2 b\big)}
{(|\wt u-
(1/a) \wt e_1| + 1
)^{d+\alpha}}
 =  \frac{\Psi\big((|\wt u-
(1/a) \wt e_1| + 1
)^2 b\big)b^{(d+\alpha)/2}}
{[(|\wt u-
(1/a) \wt e_1| + 1
)^2b]^{(d+\alpha)/2}}
 \le c  
 \frac{\Psi\big((|\wt u| + 1
)^2 b\big)b^{(d+\alpha)/2}}
{[(|\wt u| + 1
)^2b]^{(d+\alpha)/2}}\nn\\
&\le c  
 \frac{\Psi\big((|\wt u| + 1
)^2 b\big)}
{(|\wt u| + 1
)^{d+\alpha}}
= c  \Psi(b)
 \frac{\Psi\big((|\wt u| + 1
)^2 b\big)}
{\Psi(b)(|\wt u| + 1
)^{d+\alpha}}\le c 
 \frac{ \Psi(b)}{|\wt u|^{d+\alpha-2{\gamma_2}_+}}. \label{e:Line1}
   \end{align}

   Moreover, 
   \begin{align}
  \frac{\Psi\big((|\wt u-
(1/a) \wt e_1| + 1
)^2 b\big)}
{(|\wt u-
(1/a) \wt e_1| + 1
)^{d+\alpha}}
\le c
\frac{ \Psi(b)}{(|\wt u-(1/a) \wt e_1| +1)^{d+\alpha-2{\gamma_2}_+}}. \label{e:Line2}
  \end{align}
Using \eqref{e:g2}--\eqref{e:Line2}, for $a\le M$ with $\wt u=(u_1,  \wh u)$ for $d \ge 3$ 
(the case $d=2$ is simpler),
 \begin{align}
I &\le 
\int_{\R^{d-1}, \frac3{2a}>u_1>\frac1{2a}}
\frac{\Psi\big((|\wt u-
(1/a) \wt e_1| + 1
)^2 b\big)}
{u_1^{d+\alpha}(|\wt u-(1/a) \wt e_1| +1)^{d+\alpha}}
 d\wt u\nn\\
&\quad  +c\int_{\R^{d-1}, \frac3{2a} \le u_1}
\frac{\Psi\big((|\wt u-
(1/a) \wt e_1| + 1
)^2 b\big)}
{(u_1+|\wh u|)^{d+\alpha}(|\wt u-
(1/a) \wt e_1| + 1
)^{d+\alpha}
}
 du_1d\wh u
 \nn\\
&  \le c  a^{d+\alpha} \Psi(b)
\int_{\R^{d-1}}
\frac{ d\wt u}{(|\wt u-(1/a) \wt e_1| +1)^{d+\alpha-2{\gamma_2}_+}}
+c\Psi(b)
\int_{\R^{d-1}, \frac3{2a} \le u_1}
\frac{ du_1d\wh u }{(u_1+|\wh u|)^{2d+2\alpha-2{\gamma_2}_+}}
\nn\\ 
&\le c  a^{d+\alpha} \Psi(b)
\int_{\R^{d-1}}
\frac{ d\wt v}{(|\wt v| +1)^{d+\alpha-2{\gamma_2}_+}}
+c\Psi(b)
\int_{\frac3{2a}}^{\infty}
\frac{du_1}
 {u_1^{d+2+2\alpha-2{\gamma_2}_+}}
\int_{\R^{d-2}}
\frac{  d\wh w}{(1+|\wh w|)^{2d+2\alpha-2{\gamma_2}_+}}
\nn\\
&  \le    c \Psi(b) (a^{d+\alpha} 
+  a^{d+\alpha+(1+\alpha-2{\gamma_2}_+)} ) \le c a^{d+\alpha} \Psi(b) \label{e:neq1}
 \end{align}
 where in the last inequality we have used the facts $1+\alpha-2{\gamma_2}_+>0$ and $a\le M$ .

For
$3/(2a)>u_1>1/(2a)$, we have
$
 b/(4a^2)\le (|\wt u| + 1
)^2 b.
$
Thus, using the fact that $\Psi(t)t^{-(d+\alpha)/2}$ is almost decreasing, we have 
   for $3/(2a)>u_1>1/(2a)$,
   \begin{align}
&  \frac{\Psi\big((|\wt u| + 1
)^2 b\big)}
{(|\wt u| + 1
)^{d+\alpha}}
 =  
\frac{\Psi\big((|\wt u| + 1
)^2 b\big)b^{(d+\alpha)/2}}
{ 
[(|\wt u| +1)^2b]^{(d+\alpha)/2}}\le c
\frac{\Psi(b/a^2)b^{(d+\alpha)/2}}
{ 
[b/(4a^2)]^{(d+\alpha)/2}}
=ca^{d+\alpha}
\Psi(b/a^2).\label{e:Line3}
  \end{align}
Using the upper bound in \eqref{e:Psi-infty0}, we have 
   \begin{align}
& \frac{\Psi\big((|\wt u| + 1
)^2 b\big)}
{(|\wt u| + 1
)^{d+\alpha}}
=  \Psi(b)
 \frac{\Psi\big((|\wt u| + 1
)^2 b\big)}
{\Psi(b)(|\wt u| + 1
)^{d+\alpha}}
\le c 
 \frac{ \Psi(b)}{(|\wt u| + 1
)^{d+\alpha-2{\gamma_2}_+}}
\le c
 \frac{ \Psi(b)}{|\wt u|^{d+\alpha-2{\gamma_2}_+}}. \label{e:Line4}
   \end{align}

Using \eqref{e:g2} and \eqref{e:Line3}--\eqref{e:Line4},
we get
  \begin{align*}
III(1) \le ca^{d+\alpha}
\Psi(b/a^2) 
\int_{\R^{d-1}}
\frac{ d\wt v }
{(|\wt v| +1)^{d+\alpha}
}
 \le ca^{d+\alpha}
\Psi(b/a^2) 
 \end{align*}
 and, by the change of variables $\wt w=(u_1,  \wh u)=(u_1, u_1 \wh w)$ for $d \ge 3$, 
   \begin{align}
II&\le c  \Psi(b)
\int_{\R^{d-1}, \frac3{2a} \le u_1}
\frac{ d\wt u}{
|\wt u|^{2d+2\alpha-2{\gamma_2}_+}}\nn\\
 &\le c\Psi(b)
\int_{\frac3{2a}}^{\infty} 
\frac{du_1}
{u_1^{d+2+2\alpha-2{\gamma_2}_+}}
\int_{\R^{d-2}}
\frac{  d\wh w}{(1+|\wh w|)^{2d+2\alpha-2{\gamma_2}_+}}
 \le c\Psi(b)
 a^{d+\alpha+(1+\alpha-2{\gamma_2}_+)}. \label{e:g22}
 \end{align}
Therefore  using  \eqref{e:g111}, \eqref{e:neq1} and the fact $a\le M$, we have 
$$
g(a, b) \le  c  \left( a^{d+\alpha}\Psi(b)+  a^{d+\alpha}\Psi(b/a^2)
\right)\asymp a^{d+\alpha}\Psi(b/a^2), \quad a \in (0, M]. 
$$

We now show the lower bound: 
Note that, using the fact $a\le M$, we have that  for
$3/(2a)>u_1>1/(2a)$ and $ |\wh u|<1/(2a)$,
 \begin{align}\label{e:nba}
  |\wt u| + 1 
   \le ((3/2)^2+(1/2)^2)^{1/2}a^{-1} + 1  \le (\sqrt{10}+2M)/ (2a).
   \end{align}
Using \eqref{e:Psi-infty0}, 
\eqref{e:g111}, \eqref{e:nba} and the change of variable $(v_1, \wh v)=
(u_1-(1/{a}), \wh u)$, we have
  \begin{align*}
g(a,b)  \ge &\int_{\R^{d-1}, \frac3{2a}>u_1>\frac1{2a},  |\wh u|<\frac1{2a} }
 \frac{\Psi\big((|\wt u| + 1
)^2 b\big)}
{(|\wt u| + 1)^{d+\alpha}(|\wt u-(1/a) \wt e_1| +1)^{ d+\alpha}}
 d\wt u\nn\\
 \ge   & c  a^{d+\alpha}
\Psi(b/a^2)
 \int_{\R^{d-1}, \frac3{2a}>u_1>\frac1{2a},  |\wh u|<\frac1{2a} }
 \frac{ d\wt u}
{(|\wt u-(1/a) \wt e_1| +1)^{ d+\alpha}}
\nn\\
  \ge & c  a^{d+\alpha}\Psi(b/a^2)
\int_{\R^{d-1}, |v_1|<\frac{1}{M}, |\wh v|<\frac{1}{2M}}
\frac{ d\wt v} {(|\wt v| +1)^{d+\alpha}}
 = c a^{d+\alpha}\Psi(b/a^2).
 \end{align*}
Therefore  using  \eqref{e:g111}, we get
$g(a,b) \ge  
c a^{d+\alpha}\Psi(b/a^2)$
for all  $a \in (0, M]. $
\qed

 \begin{lemma}\label{l:f_2}
For any $M>0$, there exists  
$C=C(M)>0$ such that 
\begin{align}\label{e:k11_1f}
f(a,b) \le Ca^{d+\alpha} \Psi(b) 
 + Ca^{d+\alpha}b^{\alpha+1}\Psi
\big(\frac{b}{a^2}\big), \quad b>0 \,\mbox{ and } \,  a \in (0,  M (1\wedge b)].
 \end{align}
\end{lemma}

\pf 
By \eqref{e:neq1} and 
\eqref{e:g22}, we see that 
 \begin{align}
 I+II \le  ca^{d+\alpha} \Psi(b). \label{e:f3n}
 \end{align}

Since $a \le M$,  for
$\frac3{2a}>u_1>\frac1{2a}$, we have
      $
|\wh u| +1/a \asymp  |\wh u| +u_1 + 1
\asymp |\wt u| + 1 .
$
 Using this, by  the change of variables $v_1=u_1-1/a$ and  $\hat{v}=\hat{u}$, and then $v_1=t/a$, 
  \begin{align}
&III(b) \asymp 
\int_0^{\frac1{2a}}
\int_{\R^{d-2}}
\frac{\Psi\big(( |\wh v| +1/a)^2 b\big)}{( |\wh v| +1/a)^{d+\alpha}(|\wh v|+v_1 +1/b)^{d+\alpha}}
d\wh vdv_1\nn \\
 & =
   a^{-1}
\int_0^{1/2}
\int_{\R^{d-2}}
\frac{\Psi\big(( |\wh v| +1/a)^2 b\big)}{( |\wh v| +1/a)^{d+\alpha}(|\wh v| +(\frac{a}{b}+t)/a)^{d+\alpha}}
 d\wh v dt. \label{e:newlde}
  \end{align}
Using the change of variable $  \wh v=[(\frac{a}{b}+t)/a] \wh w$,  \eqref{e:newlde}
   is equal to   
\begin{align}
&
a^{\alpha+1}\int_0^{1/2}
(\frac{a}{b}+t)^{-2-\alpha}
a^{d+\alpha}
\int_{\R^{d-2}}
\frac{\Psi\big(((\frac{a}{b}+t) |\wh w| +1)^2 \frac{b}{a^2}\big)}{( (\frac{a}{b}+t) |\wh w| +1)^{d+\alpha}(|\wh w| +1)^{d+\alpha}}
 d\wh w dt\nn\\
&= c a^{d+2\alpha+1}
\int_{\frac{a}{b}}^{\frac{a}{b}+ 1/2}t^{-2-\alpha}
\int_0^\infty 
\frac{\Psi\big((ts +1)^2 \frac{b}{a^2}\big)s^{d-3}}{( ts +1)^{d+\alpha}(s+1)^{d+\alpha}}
 ds dt.\label{e:newlde1}
  \end{align}
Using the upper bound in \eqref{e:Psi-infty0},  
 for $a/b\le t\le a/b+1/2 \le M+1/2$,  
      \begin{align*}
  \int_0^\infty 
\frac{\Psi\big((ts +1)^2 \frac{b}{a^2}\big)s^{d-3}}{( ts +1)^{d+\alpha}(s+1)^{d+\alpha}}
 ds 
 & \le c\Psi\big(\frac{b}{a^2}\big) 
  \int_0^\infty 
\frac{(s+1)^{-d-\alpha}s^{d-3}ds}{  ( ts +1)^{d+\alpha-{\gamma_2}_+}}
 \le c  \Psi\big(\frac{b}{a^2}\big) \int_0^\infty 
\frac{s^{d-3}ds}{(s+1)^{d+\alpha}}.
   \end{align*}
   Therefore \eqref{e:newlde1} is less than or equal to 
         \begin{align*}
c a^{d+2\alpha+1}\Psi\big(\frac{b}{a^2}\big) 
\int_{\frac{a}{b}}^{\frac{a}{b}+ 1/2} t^{-2-\alpha}
 dt\asymp a^{d+2\alpha+1} (\frac{a}{b})^{-\alpha-1}\Psi\big(\frac{b}{a^2}\big) 
 \asymp a^{d+\alpha}b^{\alpha+1}\Psi\big(\frac{b}{a^2}\big) .
 \end{align*}
Therefore  using  this,  \eqref{e:f2n}  and \eqref{e:f3n}, we obtain \eqref{e:k11_1f}.
\qed

Note that for all $a>0$ and $p \in [(1\wedge\alpha)+1, \infty )$, 
 \begin{equation}\label{e:Psi-int}
 0 <\int_{a}^\infty
 v^{-p}\Psi(v)
 dv \le c\Psi(a)a^{-\gamma_2}\int_{a}^\infty
 v^{-p+\gamma_2}
 dv \le c(a, p) <\infty
 \end{equation}
 and 
 \begin{equation}\label{e:Psi-int2}
 \int_{0}^{1/a}
 u^{p-2}\Psi(1/u)
 du=\int_{a}^\infty
 v^{-p}\Psi(v)
 dv \asymp 1,
 \end{equation}
 with comparison constants depending on $a$.

For all $x, y \in \R^d_+$, let
\begin{align}\label{whI_211}
\Xi(x,y):= \int_{x_d}^\infty
z_d^{\alpha} 
\int_{\R^{d-1}}\frac{\Psi\big({(|\wt z|+1+z_d)^2 }/{(y_d z_d)}\big)}
{(|\wt z| + z_d)^{d+\alpha}(|\wt z| +1+ z_d)^{d+\alpha}}
 d\wt z
  dz_d.
\end{align}

\begin{lemma}\label{l:I2}
For all $x, y \in \R^d_+$  with $|x-y|=\sqrt 2$ and $y_d \ge x_d$, 
we have 
$$
\Xi(x,y) \asymp
\begin{cases}  
 x_d^{-d-\alpha} 
&\text{ for } x_d >1/4;\\[2pt]
\int_{\frac{1}{2y_d}}^{\frac{1}{x_dy_d}}
\Psi(v)\frac{dv}{v}
&\text{ for } x_d  \le1/4.
\end{cases}
$$
\end{lemma}
\pf
By the change of variables $\wt u=\wt z/z_d$ , 
we get
\begin{align*}
\Xi&=\int_{x_d}^\infty
z_d^{\alpha} 
\int_{\R^{d-1}}z_d^{-2(d+\alpha)} z_d^{d-1}
\frac{\Psi\big({(|\wt u|+1+(1/z_d))^2 }(z_d/y_d)\big)}{(|\wt u| + 1)^{d+\alpha}(|\wt u| + 1+ (1/z_d))^{d+\alpha}}
 d\wt u
  dz_d\nn \\
  &=\int_{x_d}^\infty
z_d^{-d-\alpha-1}\int_{\R^{d-1}}
\frac{\Psi\big({(|\wt u|+1+(1/z_d))^2 }(z_d/y_d)\big)}{(|\wt u| + 1)^{d+\alpha}(|\wt u| + 1+ (1/z_d))^{d+\alpha}}
 d\wt u
  dz_d.
\end{align*}

\noindent
{\it Case 1:} {\bf $x_d\ge 1/4$.}
In this case, $y_d \asymp x_d \ge 1/4$ 
so using \eqref{e:Psi-infty0}, 
we have that for $z_d \ge x_d \ge 1/4$,
\begin{align}
&\int_{\R^{d-1}}
\frac{\Psi\big({(|\wt u|+1+(1/z_d))^2 }(z_d/y_d)\big)}{(|\wt u| + 1)^{d+\alpha}(|\wt u| + 1+ (1/z_d))^{d+\alpha}}
 d\wt u= h(z_d, z_d/y_d)\asymp h(z_d, z_d/x_d).\nn
 \end{align}
Thus, by Lemma \ref{l:k} and 
\eqref{e:Psi-int}, 
for $x_d  \ge 1/4$ (so $x_d \asymp y_d$),
\begin{align*}
\Xi &
\asymp \int_{x_d}^\infty
z_d^{-d-\alpha-1}\Psi\big(\frac{z_d}{x_d}\big)
dz_d\asymp  x_d^{-d-\alpha}\int_{1}^\infty
\frac{\Psi(v)dv}{v^{d+\alpha+1}}\asymp  x_d^{-d-\alpha}.
\end{align*}

\noindent
{\it Case 2:} {\bf $x_d< 1/4$.}
In this case, 
by Lemma \ref{l:k}, 
\begin{align*}
\Xi = \int_{x_d}^\infty
z_d^{-d-\alpha-1}h(z_d,z_d/y_d) dz_d
\asymp \int_{2}^{\infty}
z_d^{-d-\alpha-1} \Psi\big(\frac{z_d}{y_d}\big) 
dz_d+ \int_{x_d}^{2}
z_d^{-1} \Psi\big(\frac{1}{z_dy_d}\big) 
 dz_d. 
\end{align*}
Note that 
 $$ 1 \asymp
 c\int_{2}^{\infty}
z_d^{-d-\alpha-1-{\gamma_1}_-} dz_d  \le 
 \int_{2}^{\infty}
z_d^{-d-\alpha-1} \frac{ \Psi\big(\frac{z_d}{y_d}\big) }
{ \Psi\big(\frac{1}{y_d}\big)  } dz_d \le 
c\int_{2}^{\infty}
z_d^{-d-\alpha+{\gamma_2}_+} dz_d 
\asymp 1,
$$
and
$$
\int_{\frac{1}{2y_d}}^{\frac{1}{x_dy_d}}
\Psi(v)\frac{dv}{v}
\ge \int_{\frac{1}{2y_d}}^{\frac{4}{y_d}}
\Psi(v)\frac{dv}{v} \asymp  \Psi\big(\frac{1}{y_d}\big) \int_{\frac{1}{2y_d}}^{\frac{4}{y_d}}
v^{-1} 
 dv  \asymp  \Psi\big(\frac{1}{y_d}\big).
$$
Thus, 
$\Xi \asymp  
\Psi\big(\frac{1}{y_d}\big) 
+ \int_{x_d}^{2}
z_d^{-1} \Psi\big(\frac{1}{z_dy_d}\big) 
 dz_d \asymp \int_{\frac{1}{2y_d}}^{\frac{1}{x_dy_d}}
\Psi(v)\frac{dv}{v}. 
 $
\qed

Suppose that $x, y \in \R_+^d$, 
$\wt x=\wt 0$, $|x-y|=\sqrt 2$, $y_d \ge x_d$ and
$y=(|\wt y| \wt e_1, y_d)$.
Let 
$$
I_1=I_1(x, y):=\int_0^{x_d}
z_d^{\alpha} 
\int_{\R^{d-1}}
\frac{
\Psi\big({(|\wt z-|\wt y| \wt e_1|+y_d)^2 }/{(y_d z_d)}\big)
}{(|\wt z| + x_d)^{d+\alpha}(|\wt z-|\wt y| \wt e_1| + y_d)^{d+\alpha}}
 d\wt z
  dz_d
  $$
and 
$$
I_2=I_2(x, y):=\int_{x_d}^\infty
z_d^{\alpha} 
\int_{\R^{d-1}}
\frac{\Psi\big({(|\wt z-|\wt y| \wt e_1|+y_d+z_d)^2 }/{(y_d z_d)}\big)}{
(|\wt z| + z_d)^{d+\alpha}
(|\wt z-|\wt y| \wt e_1| + y_d+ z_d)^{d+\alpha}}
 d\wt z
  dz_d.
$$
Since $x_d \asymp x_d+z_d$ and $y_d \asymp y_d+z_d$ if $z_d\le x_d$ and 
$z_d \asymp x_d+z_d$  if $z_d\ge x_d$, 
we see that 
\begin{align}\label{e:k10}
q(x,y) \asymp I_1(x, y)+I_2(x, y).
\end{align}

\begin{prop}\label{p:lower}
For all $x, y \in \R^d_+$  with $|x-y|=\sqrt 2$, 
we have 
\begin{align}\label{e:whks}
&q(x,y) \ge c\begin{cases} 
(x_d \wedge y_d)^{-d-\alpha}\asymp(x_d \vee y_d)^{-d-\alpha}
&\text{for } x_d \wedge y_d >1/4;\\[2pt]
\int_{1}^{\frac{1}{x_dy_d}}
\Psi(u)\frac{du}{u}
&\text{for }x_d \wedge y_d  \le1/4.
\end{cases}
\end{align}
\end{prop}
\pf
Suppose that $x, y \in \R_+^d$, 
$\wt x=\wt 0$, $|x-y|=\sqrt 2$ and $y_d \ge x_d$.
Without loss of generality we assume that 
$y=(|\wt y| \wt e_1, y_d)$.
Since 
\begin{align}
\label{e:g1000p}
 |\wt z-|\wt y| \wt e_1| + y_d \le |\wt z|+|\wt y|+ y_d-x_d +x_d \le  |\wt z|+ 2\sqrt 2+x_d,
\end{align}
we have that, for $z_d \ge  x_d$, 
\begin{align}
\label{e:g100p}
 |\wt z-|\wt y| \wt e_1| + y_d +z_d \le |\wt z|+ 2\sqrt 2+ 2z_d.
\end{align}

Since $t \to t^{-(d+\alpha)/2}\Psi\big(t)$ is almost decreasing, using \eqref{e:g100p} we have that for $z_d \ge  x_d$,
\begin{align}\label{whI_21}
&\frac{\Psi\big({(|\wt z-|\wt y| \wt e_1|+y_d+z_d)^2 }/{(y_d z_d)}\big)}{(|\wt z-|\wt y| \wt e_1| + y_d+ z_d)^{d+\alpha}}\nn\\
&=\frac{1}{(y_d z_d)^{(d+\alpha)/2}}\frac{\Psi\big({(|\wt z-|\wt y| \wt e_1|+y_d+z_d)^2 }/{(y_d z_d)}\big)}{[(|\wt z-|\wt y| \wt e_1| + y_d+ z_d)^2/(y_d z_d)]^{(d+\alpha)/2}}\nn \\
&\ge \frac{c}{(y_d z_d)^{(d+\alpha)/2}}\frac{\Psi\big({(|\wt z|+1+z_d)^2 }/{(y_d z_d)}\big)}
{[(|\wt z| +1+ z_d)^2/(y_d z_d)]^{(d+\alpha)/2}}=c \frac{\Psi\big({(|\wt z|+1+z_d)^2 }/{(y_d z_d)}\big)}
{(|\wt z| +1+ z_d)^{d+\alpha}}.
\end{align}
Thus, 
$I_2\ge c \, \Xi$,
where $\Xi=\Xi(x,y)$ is the function defined in \eqref{whI_211}.

If $|\wt y|\le 1/2$ and $x_d <1/4$, then 
\begin{align}
\label{e:g100}
y_d-x_d=\sqrt{2-|\wt y|^2}  \ge \sqrt{2-1/4}= \sqrt 7 /2.
\end{align} 
Thus $\sqrt 7 /2 \le y_d$. 
Applying  Lemma \ref{l:I2}, we get that
for $|\wt y|\le 1/2$ and $x_d <1/4$, 
 \begin{align*}
I_2  \ge c  \Xi \ge \int_{\frac{1}{2y_d}}^{\frac{1}{x_dy_d}}
\Psi(v)\frac{dv}{v}\ge c  \int_{\frac{1}{ \sqrt 7}}^{\frac{1}{x_dy_d}}
 \Psi(v)  \frac{dv}{v}\ge c  \int_{1}^{\frac{1}{x_dy_d}}
 \Psi(v) 
\frac{dv}{v}. 
\end{align*}

Since $x_d \le y_d$, 
  by the change of variables $\wt u=\wt z/z_d$, 
we get
\begin{align*}
 I_2 & \ge c \int_{y_d}^\infty
z_d^{\alpha} 
\int_{\R^{d-1}}
\frac{\Psi\big({
(|\wt z-|\wt y| \wt e_1|+z_d)^2 }/{(y_d z_d)
}\big)}{(|\wt z| + z_d)^{d+\alpha}(|\wt z-|\wt y| \wt e_1| + z_d)^{d+\alpha}}
 d\wt z
  dz_d\\\
  &\ge c\int_{y_d}^2
z_d^{-d-\alpha-1} 
g(z_d/|\wt y|, z_d/y_d)
  dz_d=: c I_3.
  \end{align*}

If $|\wt y| \ge 1/2$, then
$\sqrt 2= |x-y| \ge | \wt y|  \ge 1/2.  
$
 Thus,
 by Lemma \ref{l:g_2},   we get that for $|\wt y| \ge 1/2$ and $x_d <1/4$,
 \begin{align*}
&I_2 \ge c  I_3\ge c\int_{y_d}^2
z_d^{-d-\alpha-1} z_d^{d+\alpha} 
\Psi\big(\frac{z_d}{y_d}\big)
  dz_d 
   =
  \int_{y_d}^{2}
z_d^{-1} \Psi\big(\frac{z_d}{y_d}\big)  dz_d\asymp 
 \int_{1}^{\frac{2}{y_d}}
\Psi(u)\frac{du}{u}. 
 \end{align*}
Thus, combining this with Lemma \ref{l:I2} for $|\wt y| \ge 1/2$ and 
$x_d <1/4$,
\begin{align}
I_2& \ge c  (I_3+\Xi) \ge c \int_{1}^{\frac{2}{y_d}}
\Psi(u)\frac{du}{u}+ c\int_{\frac{1}{2y_d}}^{\frac{1}{x_dy_d}}
\Psi(v)\frac{dv}{v}  \nn\\
& \ge c \int_{1}^{\frac{2}{y_d}}
\Psi(u)\frac{du}{u}+ c\int_{\frac{2}{y_d}}^{\frac{1}{x_dy_d}}
\Psi(v)\frac{dv}{v} = c  \int_{1}^{\frac{1}{x_dy_d}}
\Psi(u)\frac{du}{u}.
 \label{e:l1111}
\end{align}
We now conclude from 
Lemma \ref{l:I2}
and \eqref{e:l1111}
that \eqref{e:whks} holds.
\qed

\begin{prop}\label{p:upper}
There exists a constant $C>0$ such that for all $x, y \in \R^d_+$  with $|x-y|=\sqrt 2$,
\begin{align}\label{e:whks_u}
&q(x,y) \le C\begin{cases} 
(x_d \wedge y_d)^{-d-\alpha}\asymp(x_d \vee y_d)^{-d-\alpha},
& x_d \wedge y_d >1/4;\\
\int_{1}^{\frac{1}{x_dy_d}}
\Psi(u)\frac{du}{u},
&x_d \wedge y_d  \le1/4.
\end{cases}
\end{align}
\end{prop}
\pf
Suppose that $x, y \in \R_+^d$ and $|x-y|=\sqrt 2$.
Without loss of generality we assume that 
$\wt x=\wt 0$, $y_d \ge x_d$ and 
$y=(|\wt y| \wt e_1, y_d)$.

\noindent
{\it Case 1,} {\bf $|\wt y| \le 1/2$:}
Suppose that $|\wt y| \le 1/2$. 
Then by \eqref{e:g100},
$\sqrt 2= |x-y| \ge y_d-x_d \ge \sqrt 7/2$ 
and
$y_d -x_d -|\wt y| \ge  (\sqrt 7-1)/2> 1/2$, 
and so, 
$$
|\wt z-|\wt y| \wt e_1| + y_d  \ge |\wt z|-|\wt y| + y_d 
  = |\wt z| +(y_d -x_d -|\wt y|) +x_d
\ge |\wt z|+1/2+x_d 
$$ 
and 
\begin{align}
\label{e:g10001p}
|\wt z-|\wt y| \wt e_1| + y_d+ z_d \ge |\wt z|+1/2+x_d + z_d  \ge |\wt z|+1/2+z_d.
\end{align}
Thus, using \eqref{e:g1000p} and  the change of variables $\wt u=\wt z/x_d$,
\begin{align*}
I_1&\asymp \int_0^{x_d}
z_d^{\alpha} 
\int_{\R^{d-1}}
\frac{
\Psi\big({(|\wt z| + 1+x_d)^2 }/{(y_d z_d)}\big)
}{
(|\wt z| + x_d)^{d+\alpha}(|\wt z| + 1+x_d)^{d+\alpha}}
 d\wt z
  dz_d \asymp 
 x_d^{-d-2\alpha-1}\int_0^{x_d}
z_d^{\alpha}   
h(x_d, x_d^2/(y_dz_d)) dz_d. 
\end{align*}
Since $\alpha>{\gamma_2}_+$, by \eqref{e:Psi-infty0},
$$
 \int_0^{x_d}
z_d^{\alpha}  
\Psi\big(\frac{1}{z_d}\big) dz_d
\le c 
 \Psi\big(\frac{1}{x_d}\big)x_d^{{\gamma_2}_+}
  \int_0^{x_d}
z_d^{\alpha-{\gamma_2}_+}  
 dz_d\le c 
\Psi\big(\frac{1}{x_d}\big)x_d^{\alpha+1}.
$$
Thus, 
by Lemma \ref{l:k}, \eqref{e:Psi-int2} and the fact that $x_d \asymp y_d$ if $x_d >1/4$, 
Thus, \begin{align}
  I_1 \le c
\begin{cases} 
\frac{1}{x_d^{d+2\alpha+1}}
\int_0^{x_d}
z_d^{\alpha}  
 \Psi\big(\frac{x_d}{z_d}\big) dz_d =
\frac1{x_d^{d+\alpha}}   \int_0^{1}
u^{\alpha}  
\Psi\big(\frac{1}{u}\big) du
\asymp \frac1{x_d^{d+\alpha}}, 
& x_d >1/4;\\
x_d^{1+\alpha}  
 \int_0^{x_d}
z_d^{\alpha}  
\Psi\big(\frac{1}{z_d}\big)dz_d \le c 
 \Psi\big(\frac{1}{x_d}\big),
&x_d  \le1/4.
\end{cases}\label{e:yI1}
\end{align}

On the other hand, 
by \eqref{e:g100p} 
and 
\eqref{e:g10001p}, we have 
$I_2\asymp \Xi$. 
Since $\sqrt 7 /2 \le y_d<7/4$ for  $x_d <1/4$ (because $|\wt y| \le 1/2$), by Lemma \ref{l:I2} for $x_d  \le1/4$,
$$I_2\asymp 
\int_{\frac{1}{2y_d}}^{\frac{1}{x_dy_d}}
\Psi(v)\frac{dv}{v} \ge \int_{\frac{1}{2x_dy_d}}^{\frac{1}{x_dy_d}}
\Psi(v)\frac{dv}{v}
\asymp \Psi\big(\frac{1}{x_dy_d}\big)\int_{\frac{1}{2x_dy_d}}^{\frac{1}{x_dy_d}}
\frac{dv}{v} \asymp \Psi\big(\frac{1}{x_d}\big) \ge c I_1
$$
and
$$
\int_{\frac{1}{2y_d}}^{1}
\Psi(v)\frac{dv}{v} \le 
\int_{\frac{2}{7}}^{1}
\Psi(v)\frac{dv}{v} \asymp 1 \asymp 
\int_{1}^{2}
\Psi(v)\frac{dv}{v}
\le \int_{1}^{\frac{1}{x_dy_d}}
\Psi(v)\frac{dv}{v}.
$$
Thus from these and \eqref{e:yI1}, we see that \eqref{e:whks_u} holds true 
 for $|\wt y| \le 1/2$.

\noindent
{\it Case 2,} {\bf $|\wt y| \ge 1/2$:}
Suppose that $|\wt y| \ge 1/2$. Then
\begin{align}
\label{e:g1}
\sqrt 2= |x-y| \ge | \wt y|  \ge 1/2.  
\end{align}

Since $\Psi(t)t^{-(d+\alpha)/2}$ is almost decreasing and $x_d \le y_d$, 
by the  argument in \eqref{whI_21} we have 
\begin{align}
I_1\le c\int_0^{x_d}
z_d^{\alpha} 
\int_{\R^{d-1}}
\frac{\Psi\big({(|\wt z-|\wt y| \wt e_1|+x_d)^2 }/{(y_d z_d)}\big)}{(|\wt z| + x_d)^{d+\alpha}(|\wt z-|\wt y| \wt e_1| + x_d)^{d+\alpha}}
 d\wt z
  dz_d=:cJ_1
\end{align}
and
\begin{align*}
I_2\le c 
\int_{x_d}^{\infty}z_d^{\alpha} 
\int_{\R^{d-1}}
\frac{\Psi\big({(|\wt z-|\wt y| \wt e_1|+z_d)^2 }/{(y_d z_d)}\big)}
{(|\wt z|+z_d)^{d+\alpha}(|\wt z-|\wt y| \wt e_1| +  z_d)^{d+\alpha}}
d\wt z
  dz_d=:c\wh J_2.
\end{align*}
By the change of variables $\wt u=\wt z/x_d$ in $J_1$,  the change of variables $\wt u=\wt z/z_d$  in $
\wh J_2$, 
we get
\begin{align*}
J_1
=  \frac{1}{x_d^{d+2\alpha+1}}\int_0^{x_d}
z_d^{\alpha} 
g(x_d/|\wt y|, x_d^2/(y_dz_d)) dz_d, 
\quad
\wh J_2
  = \int_{x_d}^\infty
\frac{
g(z_d/|\wt y|, z_d/y_d)}
{z_d^{d+\alpha+1}} 
  dz_d.
  \end{align*}

Since $x_d \asymp y_d$ for $x_d  >1/4$, 
by \eqref{e:Psi-infty0}, \eqref{e:Psi-int} and \eqref{e:Psi-int2},  
we have that for $x_d  >1/4$, 
 \begin{align*}
 & \int_{x_d}^\infty
\frac{
\Psi({z_d}/{y_d})}
{z_d^{d+\alpha+1}}dz_d
 \asymp \int_{x_d}^\infty
\frac{\Psi({z_d}/{x_d})}
{z_d^{d+\alpha+1}}
dz_d
=x_d^{-d-\alpha} \int_{1}^\infty
\Psi(u)\frac{du}  {u^{d+\alpha+1}}
\asymp x_d^{-d-\alpha}
   \end{align*}
     and
$$
  \int_0^{x_d}
z_d^{\alpha} 
\Psi\big(\frac{x_d^2}{z_dy_d}\big) dz_d
 \asymp 
  \int_0^{x_d}
z_d^{\alpha} 
\Psi\big(\frac{x_d}{z_d}\big) dz_d\asymp  x_d^{\alpha+1}   \int_0^{1}
u^{\alpha} 
\Psi(1/u) du \asymp x_d^{\alpha+1}. 
$$
  Thus by Lemma \ref{l:g_2n}, for $x_d  >1/4$, we get 
\begin{align*}
q(x,y) \le c(J_1+\wh J_2)
 \le \frac{c}
{x_d^{d+2\alpha+1}}\int_0^{x_d}
z_d^{\alpha} 
\Psi\big(\frac{x_d^2}{z_dy_d}\big) dz_d+ \int_{x_d}^\infty
\Psi\big(\frac{z_d}{y_d}\big)\frac{dz_d}{z_d^{d+\alpha+1}}
\asymp x_d^{-d-\alpha}.
   \end{align*}
   
For the remainder of the proof, we assume  $x_d  \le1/4$.
Clearly, 
\begin{align*}
I_2 \asymp &
\int_{x_d}^{y_d}
z_d^{\alpha} 
\int_{\R^{d-1}}
\frac{\Psi\big({(|\wt z-|\wt y| \wt e_1|+y_d)^2 }/{(y_d z_d)}\big)}{
(|\wt z| + z_d)^{d+\alpha}
(|\wt z-|\wt y| \wt e_1| + y_d)^{d+\alpha}}
 d\wt z
  dz_d\\
  +& 
 \int_{y_d}^\infty
z_d^{\alpha} 
\int_{\R^{d-1}}
\frac{\Psi\big({(|\wt z-|\wt y| \wt e_1|+z_d)^2 }/{(y_d z_d)}\big)}{(|\wt z| + z_d)^{d+\alpha}(|\wt z-|\wt y| \wt e_1| +  z_d)^{d+\alpha}}
 d\wt z
  dz_d=:J_3+J_2.
\end{align*}
By   the change of variables $\wt u=\wt z/z_d$  in $J_2$ and the change of variables $\wt u=\wt z/y_d$  in $J_3$, 
  \begin{align*}
J_2= \int_{y_d}^\infty
\frac{
g(z_d/|\wt y|, z_d/y_d)}
{z_d^{d+\alpha+1}} 
  dz_d
\quad\text{and}\quad
  J_3=  \frac{1}{y_d^{d+2\alpha+1}}  \int_{x_d}^{y_d}  z_d^{\alpha} 
  f(y_d/|\wt y|, y_d/z_d)dz_d.
  \end{align*}
  Since $x_d  \le1/4$,   we get
     \begin{align*}
x_d^{-\alpha-1}\int_0^{x_d}
z_d^{\alpha} 
\Psi\big(\frac{1}{z_dy_d}\big)dz_d
\le c x_d^{-\alpha-1} \Psi( \frac{1}{x_dy_d}) 
 x_d^{{\gamma_2}_+}\int_0^{x_d}
z_d^{\alpha-{\gamma_2}_+} 
 dz_d \asymp
\Psi( \frac{1}{x_dy_d}) .
     \end{align*}     
   
     Thus, by Lemma \ref{l:g_2} and the fact $|\wt y| \asymp 1$ by \eqref{e:g1} (and recalling 
$\Psi(t)  \equiv \Psi(2)>0$ on $[0, 2)$),
  \begin{align}
   J_1&\asymp x_d^{-\alpha-1}\int_0^{x_d}
z_d^{\alpha} 
\Psi\big(\frac{1}{z_dy_d}\big)  dz_d\le c  \Psi( \frac{1}{x_dy_d}).
 \label{e:betge0}
       \end{align}    
Note that, by \eqref{e:Psi-int},
   \begin{align*}
  \int_{2}^\infty
v^{-d-\alpha-1}\Psi\big(\frac{v}{y_d}\big)dv
\le c \Psi\big(\frac{1}{y_d}\big)
  \int_{2}^\infty 
v^{-d-\alpha-1+{\gamma_2}_+}dv
\asymp \Psi\big(\frac{1}{y_d}\big).
     \end{align*}
     Thus,
 by   Lemmas \ref{l:g_2n} and \ref{l:g_2} with the fact $|\wt y| \asymp 1$, 
  \begin{align}\label{e:k15p}
  J_2\le c \Psi\big(\frac{1}{y_d}\big) + c  \int_{y_d}^{2}
z_d^{-1} 
\Psi\big(\frac{1}{z_dy_d}\big)
  dz_d\le 
 c \Psi\big(\frac{1}{y_d}\big) +  c \int_{\frac{1}{2y_d}}^{\frac{1}{x_dy_d}}
\Psi(u)\frac{du}{u}.
  \end{align}
  For $J_3$, we use Lemma \ref{l:f_2} with the fact $|\wt y| \asymp 1$ (so that  $(z_d/|\wt y|) \le (y_d/|\wt y|) \le c$ for $z_d \le y_d$),
    \begin{align*}
J_3&\le c
\int_{x_d}^{y_d}  z_d^{-1} \Psi\big(\frac{1}{y_dz_d}\big)
dz_d+c
y_d^{-\alpha-1} 
\int_{x_d}^{y_d}  z_d^{\alpha} 
\Psi\big(\frac{y_d}{z_d}\big)
dz_d=
c\int_{\frac{1}{y_d^2}}^{\frac{1}{x_dy_d}}   \Psi(u)\frac{du}{u}
+c
\int_{1}^{\frac{y_d}{x_d}} 
\Psi(u)\frac{du}{ u^{2+\alpha}}
\nn \\
&
 \le c \int_{\frac{1}{4}}^{\frac{1}{x_dy_d}}   \Psi(u)\frac{du}{u}
+
c\int_{1}^{\frac{y_d}{x_d}}  
\Psi(u)\frac{du}{u}
\asymp 
 \int_{\frac{1}{4}}^{\frac{1}{x_dy_d}}  \Psi(u)\frac{du}{u}
+
\int_{\frac14}^{\frac{y_d}{4x_d}}  
\Psi(u)\frac{du}{u} 
.
\end{align*}
 Since   \begin{align*} 
 \int_{\frac{1}{4}}^{\frac{1}{x_dy_d}}  \Psi(u)\frac{du}{u}
+
\int_{\frac14}^{\frac{y_d}{4x_d}}  
\Psi(u)\frac{du}{u} 
\le
 2\int_{\frac{1}{4}}^{1}  \Psi(u)\frac{du}{u}+2\int_{1}^{\frac{1}{x_dy_d}}   \Psi(u)\frac{du}{u},
\end{align*}
 we obtain 
 we obtain 
    \begin{align}\label{e:k15pp}
    J_3&  \le c+c\int_{1}^{\frac{1}{x_dy_d}}   \Psi(u)\frac{du}{u}.
    \end{align}
   Using $x_dy_d \le (\sqrt 2 +1/4)/4 <7/16$, we get
    \begin{align} \label{e:fin1}
&\int_{1}^{\frac{1}{x_dy_d}}
\Psi(u)\frac{du}{u} \ge \frac13 \left(\int_{\frac{1}{2}\frac1{x_dy_d}}^{\frac1{x_dy_d}} +
  \int_{\frac{2}{y_d}}^{\frac{4}{y_d}}
+ \int_{1}^{2} \right)
\Psi(u)\frac{du}{u} \nn \\
    &\asymp \Psi\big(\frac{1}{x_dy_d}\big) 
  \int_{\frac{1}{2}\frac1{x_dy_d}}^{\frac1{x_dy_d}}\frac{du}{u}
  + \Psi\big(\frac{1}{y_d}\big)
  \int_{\frac{2}{y_d}}^{\frac{4}{y_d}}\frac{du}{u}
+
 1\asymp \Psi\big(\frac{1}{x_dy_d}\big) + \Psi\big(\frac{1}{y_d}\big) +
 1.
\end{align}
Therefore, we conclude from
\eqref{e:betge0}--\eqref{e:fin1}
 that 
   \begin{align*} q(x,y) & \le
   c(J_1+J_2+J_3)\le
     c+ c\Psi( \frac{1}{x_dy_d}) +  c \Psi\big(\frac{1}{y_d}\big) + c \int_{1}^{\frac{1}{x_dy_d}}   \Psi(u)\frac{du}{u}\asymp \int_{1}^{\frac{1}{x_dy_d}}
\Psi(v)\frac{dv}{v}.
    \end{align*}
    \qed

Recall that $J(x,y)= j(x, y)+q(x,y)$ and 
$\sB(x,y)=J(x,y)/j(x, y)=1+q(x,y)/j(x, y)$,
so that $\sB(x,y)-1=q(x,y)/j(x, y)$.

\smallskip
\noindent
\textbf{Proof of Theorem \ref{t:estimate}:}
Using \eqref{e:whsc}, \eqref{e:whks1}--\eqref{e:whks3} follow from Propositions \ref{p:lower} and \ref{p:upper}.
The assertions \eqref{e:whks4}-\eqref{e:whks5} follow from \eqref{e:whks3} and Lemma \ref{l:phsnew}(a)-(b).
The proof is now complete.
\qed

\textbf{Declaration of competing interest:}
The authors declare that they have no known competing financial interests or personal
relationships that could have appeared to influence the work reported in this paper.

Data sharing not applicable to this article as no datasets were generated or analysed during the current study.

\end{document}